~
\vskip 0.45cm
\font\we=cmb10 at 14.4truept
\font\li=cmb10 at 12truept
\noindent
\centerline {\we A Program for Geometric Arithmetic}
\vskip 1.0cm
\centerline {\li Lin WENG}
\vskip 1.0cm
In this article, we originate a program for what I call Geometric 
Arithmetic. Such a
program would consist of four parts, if I were able to properly understand
the essentials now. Namely,  (1) Non-Abelian Class Field Theory; (2) Geo-Ari
Cohomology Theory; (3) New Non-Abelian Zeta Functions;  and (4)  Riemann
Hypothesis. However, here I could only provide the reader with ${{1+1(={1\over
2}+{1\over 2})+1}\over 4}$ of them. To be
more precise, discussed in this article are the following particulars;
\vskip 0.30cm
\noindent
(A) Representation of Galois Group, Stability and Tannakian Category;

\noindent
(B) Moduli Spaces, Riemann-Roch, and New Non-Abelian Zeta Function; and

\noindent
(C) Explicit Formula, Functional Equation and Geo-Ari Intersection.
\vskip 0.30cm
So what are these ABC of the Geometric Arithmetic?!
\vskip 0.30cm
As stated above,  (A) is aimed at establishing a Non-Abelian Class Field 
Theory. The
starting point here is the following classical result: Over a compact 
Riemann surface,
a line bundle is of degree zero if and only if it is flat, i.e.,  induced
from a representation of fundamental group of the Riemann surface. Clearly, 
being a
bridge connecting divisor classes and fundamental groups, this result may be
viewed as and is indeed a central piece of the classical (abelian) class field
theory. (See e.g., [Hilbert] and [Weil].) Thus it is  then only natural to 
give a
non-abelian generalization  of it in order to offer a non-abelian class
field theory. This  was first done by Weil. In his fundamental paper on
generalization of abelian functions [Weil1], Weil showed that over a compact
Riemann surface, a vector bundle is of degree zero if and only if it is 
induced from a
representation of fundamental group of the surface.

Thus far, two new aspects naturally emerge. That is, unitary
representations and non-compact Riemann surfaces,  reflecting
finite quotients of Galois groups and  ramifications  in
Class Fields Theory, CFT for short, respectively: In a (complex) 
representation class
of a finite group, there always exists a unitary one, while a discussion 
for compact
Riemann surfaces  results only unramified CFT. Thus mathematics demands new
results to couple with them. As it is well-known that to this end we then 
have (i)
Mumford's stability of vector bundles in terms of intersection; (ii)
Narasimhan-Seshadri's correspondence; and (iii) Seshadri's parabolic analog 
of (i)
and (ii). That is to say, now the above result of Weil is further refined 
to the follows:
Over (punctured) Riemann surfaces, (Seshadri) equivalence classes  of 
semi-stable
parabolic bundles of parabolic degree zero correspond naturally in one-to-one
to  isomorphism classes of unitary representations of
fundamental groups.
\vskip 0.30cm
On the other hand, the above results, while central, do only  parts of the 
CFT -- at
its best, the Weil-Narasimhan-Seshadri correspondence reflects a micro 
reciprocity
law. What CFT really stands should not be a relation between a single 
representation
and an isolated bundle, instead, CFT should expose  Galois groups 
intrinsically in
terms of bundles globally. Thus an integration process  aiming at 
constructing a
global theory becomes a great necessity.

It is at this point where the theory of Tannakian category enters into the 
picture.
Recall that the existing theory of Tannakian category takes the following 
forms: (i)
groups may be reconstructed from their associated  categories of
representations; (ii) Fiber functors equipped Tannakian categories
are clone categories of (i), i.e., are equivalent to the categories of 
representations; and (iii)
original groups may be recovered from the automorphism groups of  fiber
functors.
\vskip 0.30cm
At it turns out, with this strongest form of the standard theory of Tannakian
category, we have little hope to match it perfectly with the CFT we are
looking for. Fortunately, there are still room to manoeuvre, since in CFT 
we only
care about finite quotients of the associated groups,  and  in terms of
representations finite quotients correspond to what we call finitely completed
Tannakian subcategories according to Tannaka duality and van Kampen
completeness theorem. In this way, we finally  establish a
non-abelian CFT for Riemann surfaces, or better, for function fields over 
complex
numbers successfully. Main results include the Existence Theorem, the 
Conductor Theorem
and the Reciprocity Law. See e.g. Theorem A.2.4.2.
\vskip 0.30cm
By establishing a CFT for Riemann surfaces as above, possibly, we may
give the reader an impression that everything works smoothly.
No, practically, it is not the case. For example, we do not need all unitary
representations. Or put this in another way, all semi-stable parabolic 
bundles of
parabolic degree zero lead us to nowhere. Consequently, we must carefully 
select
among these semi-stable objects a handful portion so that (i) the standard
theory of Tannakian category could be applied; and (ii) there are still 
rooms for us to
manoeuvre, along the line of Tannaka duality and van Kampen
completeness theorem. This then leads to  what we call geo-ari
representations and geo-ari bundles.
\vskip 0.30cm
\noindent
{\it Remark.} By definition, as a direct consequence of the Narasimhan-Seshadri
correspondence,  the correspondence between geo-ari representation and geo-ari
bundles for function fields over complex numbers holds more or less trivially.
However, the situation changes dramatically for global fields. For example,
for curves over finite fields, we need to introduce a new principal, called the
Harder-Narasimhan correspondence, to tackle this.
\vskip 0.30cm
The experienced reader here naturally would ask how we overcome the difficulty
about tensor products of geo-ari bundles, since, generally speaking, to 
show the
tensor operation is closed is the key to apply the theory of Tannakian 
category. Here
for Riemann surfaces, two approches are available. For one, we use the
Narasimhan-Seshadri correspondence, as easily one sees that tensors of unitary
representations are again unitary. But this analytic approach is not a 
genuine one,
since a micro reciprocity law, i.e., the Weil-Narasimhan-Seshadri 
correspondence is
used. Thus a purely algebraic proof should be pursued. This then leads to 
the works
of  Kempf and Ramanan-Ramanthan on instability flags, which I call the
$KR^2$-trick. Moreover, as the original $KR^2$-trick only works for bundles 
without
parabolic structures, so to stylize the non-abelian CFT (for Riemann
surfaces), we ask for a parabolic version of the $RK^2$-trick. To
achieve this, we follow a supplementary work of Faltings  and Tataro: 
First, as in
[Fa], rewrite any geo-ari subbundle in terms of filtrations over certain 
points on the
surface,  disjoint from parabolic points; then use the GIT stability to 
check whether the
associated point for a subbundle of the tensor is semi-stable. If so, we 
are done by definition.
If not, by the instability flag of Mumford-Kempf, we obtain a modified GIT 
stable point
according to Ramanan-Ramanathan, from which,  the  intersection stability for
the tensor may finally be proved by using the intersection stability of all the
components in the tensor product as in [To].
\vskip 0.30cm
Motivated by such a success in non-abelian CFT for function fields over
complex numbers, we anticipate that in principal, the non-abelian CFT for
local and global fields  works  similarly. So the building blocks of our
program for a non-abelian CFT then are the follows:
 
(1) there should be a suitable type of representations of Galois groups, 
which we call
geometric representations and a suitable type of intersection stability for 
bundles
which we call geometric parabolic bundles such that an
analog of the Weil-Narasimhan-Seshadri Correspondence holds;

(2) there should exist subclasses of geometric representations and
geometric parabolic bundles, which we call geo-ari representations and geo-ari
bundles, respectively,  such that (i) these classes form naturally two abelian
categories,  (ii) an analog of Hader-Narasimhan Correspondence holds; and 
(iii) an
analog of $KR^2$-trick works. Thus in particular, by (1) and (2), we obtain two
equivalent (generalized) Tannakian categories together with natural fiber 
functors;

(3) The (generalized) Tannakian categories contain systems of the so-called 
finitely
completed Tannakian subcategories, so that via an analog of Tannaka Duality and
van Kampen completeness theorem,
we  obtain the so-called fundamental theorem of non-abelian CFT such as the
existence theorem, the conductor theorem and the reciprocity law.
\vskip 0.30cm
To end this brief discussion on Part (A), we would like to point out that
to realize the above mentioned 123 for our non-abelian CFT, standard theories
on GIT, Tannakian category and representations of Galois groups are far 
from being
enough. For examples,   to achieve (1), we require  (i)  a Geometric Invariant
Theory over integral bases in the spirit of Arakelov; (ii) a deformation 
theory for
geometric representations of Galois groups; and (iii) a suitable 
completeness  for
representation and stability along the line of Fountaine and Langton, 
respectively;
and to achieve (3), we require a theory of Tannakian category over
integral bases.
\vskip 0.30cm
Our next main scheme is devoted to non-abelian zeta functions.
As stated at the very beginning, Part (B) is a combination of
our partial understanding of our  new non-abelian zeta functions and
what we call geo-ari cohomology. This part to a large extent is practical 
rather than
theoretical, due to the fact that not only all studies here are based on
practical constructions, but we have not yet  understood the mathematics
involved theoretically.

Unlike for the classical  Weil zeta functions, instead of working on
general algebraic varieties and counting their rational points (over finite 
fields) in a
very primitive way, for our non-abelian zeta functions, we  concentrate  our
attentions to  moduli spaces of semi-stable bundles and count their 
rational  points
from  moduli point of view, in a similar way as what Shimura does for  Shimura
varieties.

To be more precise, consider function fields over finite fields first.
Then, for each fixed natural number $r$, we, by using a work of
Mumford-Seshadri, obtain the associated moduli spaces of
rank $r$ semi-stable bundles. In particular, with the so-called 
Harder-Narasimhan
correspondence, which claims that the rationalities of bundles and moduli 
points
coincide, we could then introduce a new type of zeta functions by considering
rational points of the moduli space as moduli points  associated with rational
semi-stable bundles.
\vskip 0.30cm
This approach, while different from that of Weil, is indeed a natural
generalization of that of Artin: When $r=1$,
our construction  recovers the classical Artin zeta functions. Moreover,
just like classical abelian zeta functions, our non-abelian version 
satisfies rationality
and a standard type of function equation as well. Since we even can give 
uniform
bounds for the coefficients of these (local) zeta functions, so via an 
Euler product, we
further introduce a more global non-abelian zeta functions for curves defined
over number fields. Needless to say,  when $r=1$, these global zeta 
functions are
nothing but the classical  Hasse-Weil zeta functions for curves. So non-abelian
arithmetic aspect of  curves is supposed to be reflected by these new zeta 
functions.
\vskip 0.30cm
Well, while this latest general statement should finally lead us to a 
mathematics
wonderland, we  have no yet found our theoretical feet. For this purpose, I 
then
turn my attention to some concrete examples. This directs us to
the study of what I call the refined Brill-Noether locus and their 
intersections:
Beyond the classical consideration, the refined Brill-Noether locus 
measures  how
automorphisms of the associated bundles change too. As a direct consequence, we
obtain a concrete reciprocity law for elliptic curves in  ranks 2 and 3.
\vskip 0.30cm
It now becomes quite apparent that key points for our construction of 
non-abelian
zeta functions are the follows: (i) moduli spaces of semi-stable bundles
admit naturally algebraic variety structures;
(ii) there exists a well-established cohomology theory, in which the 
so-called Serre
Duality and Riemann-Roch hold. So to construct non-abelian zeta functions of
number fields, we should carry out some basic researches since  both
(i) and (ii) above seem to be virgin lands in number theory.
 
However,  we cannot offer the reader  a very satisfied
GIT and a completed cohomology theory over number fields now. Fortunately, in
this article,  we manage successfully to obtain some practical items which are
sufficient to the construction of our new non-abelian zeta functions for number
fields. To say the truth, the outcome turns out to be equally nice: Not 
only the
non-abelian zeta functions for number fields could be defined as a natural
generalization of the classical Dedekind zeta functions,  these new zeta 
functions
are as canonical as they should be -- they
satisfy the  functional equation, and the residues of them at simple poles are
nothing but the volumes of what I call the Tamagawa measures of the associated
moduli spaces of semi-stable bundles (over number fields).
In particular, when rank is one, our work essentially recovers
Iwasawa's ICM talk at MIT about Dedekind zeta functions.

By saying this, I have no intention to claim that we are satisfied
with what we have achieved. Far from being it, we have little understanding of
these  new zeta functions. For examples,

(1) we have no idea now on how the non-abelian reciprocity law, which, by 
(A), are
supposed to hold naturally, could be read from our non-abelian zeta functions;

(2) generally speaking, we are less sure about the meaning of special 
values of our
non-abelian zeta functions -- We meet essential difficulties when trying to
explian our non-abelian zeta functions in terms of the existing motivic 
language.

(Recall that all classical abelian zeta functions are supposed to be 
motivic in the sense
of Grothendieck, Shimura, Deligne and Langlands, and hence that the associated
special values of classical zeta functions are supposed to have motivic
interpretations as conjectured by Beilinson, and Bloch-Kato, based on the
fundamental works of Euler, Riemann, Borel, Quillen and Tate, among others.)
\vskip 0.30cm
As to (i) and (ii) above for number fields, what is accomplished here is
the  introduction of  intersection stability, a construction of the 
corresponding
moduli spaces, and a practical formation of one dimensional geo-ari cohomology
for which the duality and Riemann-Roch hold.  In fact, the intersection 
stability may
be dated back earlier from the works of Stuhler and Grayson, despite the fact
that we work out this independently.

Note that also for the construction of non-abelian zeta functions for number
fields, main properties for moduli spaces we need are the compactness and the
existence of natural measures. So, based on Arakelov intersection theory and
Chevalley-Weil's adeles, we may easily generalize (i) to number fields.
In comparison, (ii), the key to the convergence, the functional equation, 
remains
very challenging. However,  in this article,   based on the earlier works in
particular, that of  Tate [Tate] and Schoof-van der Geer, (see also Lang 
[Lang],
Arakelov,  Szpiro, Neukirch), we are able to offer a practical definition 
of one
dimensional geo-ari cohomology in terms of Chevalley idelic language via 
Fourier
analysis. In our definition of geo-ari $h^i$'s, we make a clear
distinction between algebraic and arithmetic aspects -- algebracially,
cohomology groups are finite generated abelian groups, while
arithmetically, geo-ari cohomology is a
finite definite quantities measuring geo-arithmetical
complexities by counting all and hence infinitely many elements in
the about algebraic cohomology groups (with the help of Fourier analysis).
Thus, it would be extremely interesting to compare our sheaf theoretic approach
with Deninger's Betti cohomology approach
in which  infinite dimensional spaces are used with the help of the regularized
determinant formalism.
\vskip 0.30cm
Our practical cohomology works only in dimension one. However,
based   on linear compacity of Chevalley, as given in Iwasawa's Princeton 
lectures
notes, and Parshin's approach to duality and residue in dimension two,  we 
at the
end of  Part (C)  provide with the reader an program for what I call a
half-theoretical geo-ari cohomology in lower dimensions, which should  play 
a key
role in establishing the Hodge Index Theorem for our geo-ari intersection
introduced in (C).

Part (C) of the program is designed to give a geometric justification of 
the formal
summation $\sum_{s:\xi(s)=0}\rho^s$ for $\rho\in {\bf R}$. For this purpose,
we propose a new two dimensional geo-ari  intersection theory. This geo-ari
intersection turns out to be very interesting, since the Riemann Hypothesis 
may be
naturally studied within the framework of this model along with the line of 
Weil's
original proof of the so-called Hasse-Weil Theorem, or better, the Riemann
Hypothesis for Artin zeta functions. Due to the facts that  the Cramer 
formula is
behind the above summation for zeros of Riemann zeta and that the Explicit 
Formula
of Weil is behind the above proposed geo-ari intersection, our approach to the
Riemann Hypothesis is in appearence different from but in essence related to
that of Deninger and Quillen.

More precisely, to introduce our model on a two dimensional geo-ari 
intersection,
first we assume that there exist two dimensional mathematics sites,
which I call geo-ari surfaces; Then, as in geometry, we assume that on
these geo-ari surfaces, there are naturally divisors, which I call  micro
divisors; (Unlike in the geometric case, these micro divisors are assumed to be
parametrized by {\bf R}.) With this, motivated by the standard properties of
intersections, the Riemann-Roch in dimension one, the adjunction formula, the
global functional equation,  and the Weil explicit formula, we introduce  six
simple axioms for the intersections, consisting of one for
(permutation) symmetry, one for mirror symmetry, two for fixed points,  one
for micro explicit formula, and one for normalization.

The advantage of having this mathematics model on geo-ari intersection is 
that, as
in geometry, then the Riemann Hypothesis may be deduced from an analog of the
Hodge index theorem. Thus, motivated by what happens in geometry, we should 
also
search for a good geo-ari cohomology in dimension two. As stated above, such a
program  is proposed at the end of (C).
\vfill
\eject
\vskip 0.50cm
\centerline {\li Contents}
\vskip 0.50cm
\centerline {\li A. Representation of Galois Group, Stability and Tannakian
Category}
\vskip 0.30cm
\noindent
{\bf A.1. Summary}
\vskip 0.30cm
\noindent
{\bf A.2. Non-Abelian CFT for Function Fields over C}
\vskip 0.30cm
\noindent
A.2.1. Weil-Narasimhan-Seshadri Correspondence
\vskip 0.30cm
A.2.1.1. Unitary Representations of Fundamental Groups
\vskip 0.30cm
A.2.1.2. Semi-Stable Parabolic Bundles
\vskip 0.30cm
A.2.1.3. Weil-Narasimhan-Seshadri Correspondence: A Micro Reciprocity Law
\vskip 0.30cm
\noindent
A.2.2. Rationality: Geo-Ari Representations and Geo-Ari Bundles
\vskip 0.30cm
A.2.2.1. Branched Coverings of Riemann Surfaces
\vskip 0.30cm
A.2.2.2. Geo-Ari Representations and Geo-Ari Bundles
\vskip 0.30cm
\noindent
A.2.3.  $KR^2$-Trick and  Completed Tannakian Categories
\vskip 0.30cm
A.2.3.1. Completed Tannakian Category and van Kampen Completeness Theorem
\vskip 0.30cm
A.2.3.2.  $KR^2$-Trick
\vskip 0.30cm
\noindent
A.2.4. Non-Abelian CFT for Function Fields over Complex Numbers
\vskip 0.30cm
A.2.4.1. Micro Reciprocity Law, Tannakian Duality and the Reciprocity Map
\vskip 0.30cm
A.2.4.2.  Non-Abelian CFT
\vskip 0.30cm
\noindent
A.2.5. Classical (abelian) CFT: An Example of Kwada-Tata and Kawada
\vskip 0.30cm
A.2.5.1. Class Formation
\vskip 0.30cm
A.2.5.2. The Work of Kawada and Tate
\vskip 0.30cm
A.2.5.3. Abelian CFT for Riemann Surfaces In Terms of Geo-Ari Bundles
\vskip 0.30cm
\noindent
{\bf A.3. Towards Non-Abelian CFT for Global Fields}
\vskip 0.30cm
\noindent
A.3.1. Weil-Narasimhan-Seshadri Type Correspondence
\vskip 0.3cm
A.3.1.1. Grometric Representations
\vskip 0.30cm
A.3.1.2. Semi-Stability in terms of Intersection
\vskip 0.30cm
A.3.1.3. Weil-Narasimhan-Seshadri Type Correspondence
\vskip 0.30cm
\noindent
A.3.2. Harder-Narasimhan Correspondence
\vskip 0.30cm
\noindent
A.3.3. $KR^2$-Trick
  \vskip 0.30cm
\noindent
A.3.4. Tannakian Category Theory over Arbitary Bases
\vskip 0.30cm
\noindent
A.3.5. Non-Abelian CFT for Global Fields
\eject\vskip 0.30cm
\centerline {\li B. Moduli Spaces, Riemann-Roch, and New Non-Abelian Zeta
functions}
\vskip 0.30cm
\noindent
{\bf B.1. New Local and Global Non-Abelian Zeta Functions for Curves}
\vskip 0.30cm
\noindent
B.1.1. Local Non-Abelian Zeta Functions
\vskip 0.30cm
B.1.1.1. Artin Zeta Functions for Curves
\vskip 0.30cm
B.1.1.2. Too Different Generalizations: Weil's Zeta Functions and A New 
Approach
\vskip 0.30cm
B.1.1.3. Moduli Spaces of Semi-Stable Bundles
\vskip 0.30cm
B.1.1.4. New Local Non-Abelian Zeta Functions
\vskip 0.30cm
B.1.1.5. Basic Properties for  Non-Abelian Zeta Functions
\vskip 0.30cm
\noindent
B.1.2. Global Non-Abelian Zeta Functions for Curves
\vskip 0.30cm
B.1.2.1. Preparations
\vskip 0.30cm
B.1.2.2. Global Non-Abelian Zeta Functions for Curves
\vskip 0.30cm
B.1.2.3. Working Hypothesis
\vskip 0.30cm
\noindent
B.1.3. Refined Brill-Noether Locus for Elliptic Curves: Towards A 
Reciprocity Law
\vskip 0.30cm
B.1.3.1 Results of Atiyah
\vskip 0.30cm
B.1.3.2. Refined Brill-Noether Locus
\vskip 0.30cm
B.1.3.3. Towards A Reciprocity Law: Measuring
Refined Brill-Noether Locus Arithmetically
\vskip 0.30cm
B.1.3.4.  Examples In Ranks Two and Three: A Precise Reciprocity Law
\vskip 0.30cm
B.1.3.5. Why  Use only Semi-Stable Bundles
\vskip 0.5cm
\noindent
{\bf Appendix to B.1:  Weierstrass  Groups}
\vskip 0.3cm
\noindent
{1. Weierstrass Divisors}
\vskip 0.30cm
\noindent
{2. K-Groups}
\vskip 0.45cm
\noindent
{3. Generalized Jacobians}
\vskip 0.45cm
\noindent
{4. Galois Cohomology Groups}
\vskip 0.45cm
\noindent
{5.  Deligne-Beilinson Cohomology}
\vskip 0.5cm
\noindent
{\bf  B.2. New Non-Abelian Zeta Functions for Number Fields}
\vskip 0.30cm
\noindent
B.2.1. Iwasawa's ICM Talk on Dedekind Zeta Functions
\vskip 0.30cm
\noindent
B.2.2. Intersection Stability
\vskip 0.30cm
B.2.2.1.  Classification of Unimodular Lattices: A Global Approach
\vskip 0.30cm
B.2.2.2. Semi-Stable  Bundles over Number Fields
\vskip 0.30cm
B.2.2.3. Adelic Moduli and Its Associated Tamagawa Measure
\vskip 0.30cm
\noindent
B.2.3. Geo-Ari Duality and Riemann-Roch: A Practical Geo-Ari Cohomology 
following
Tate
\vskip 0.30cm
B.2.3.1. An Example
\vskip 0.30cm
B.2.3.2. Canonical Divisors and  Space of Different Forms
\vskip 0.30cm
B.2.3.3. Algebraic Cohomology for Matrix Divisors
\vskip 0.30cm
B.2.3.4. Geo-Arit Cohomology and Its Associated
Riemann-Roch
\vskip 0.30cm
\noindent
B.2.4. Non-Abelian Zeta Function For Number Fields
\vskip 0.30cm
B.2.4.1. The Construction
\vskip 0.30cm
B.2.4.2.  Basic Properties
\vskip 1.0cm
\centerline{\li C. Explicit Formula, Functional Equation and Geo-Ari 
Intersection}
\vskip 0.45cm
\noindent
{\bf C.1. The Riemann Hypothesis for Curves}
\vskip 0.45cm
\noindent
C.1.1. Weil's Explicit Formula: the Reciprocity Law
\vskip 0.30cm
\noindent
C.1.2. Geometric Version of Explicit Formula
\vskip 0.30cm
\noindent
C.1.3. Riemann Hypothesis for Function Fields
\vskip 0.45cm
\noindent
{\bf  C.2. Geo-Ari Intersection in Dimension Two: A Mathematics Model}
\vskip 0.45cm
\noindent
C.2.1. Motivation from Cram\'er's Formula
\vskip 0.30cm
\noindent
C.2.2. Micro Divisors
\vskip 0.30cm
\noindent
C.2.3. Global Divisors and Their Intersections: Geometric Reciprocity Law
\vskip 0.30cm
\noindent
C.2.4.  The Riemann Hypothesis
\vskip 0.30cm
\noindent
C.2.5. Not so serious Convergence Problem
\vskip 0.30cm
\noindent
C.2.6. Weil's Explicit Formula and Two Dimensional Geometric Arithmetic
Intersections
\vskip 0.30cm
\noindent
{\bf C.3. Towards A  Geo-Ari Cohomology in Lower Dimensions}
\vskip 0.30cm
\noindent
C.3.1. Classical Approach in Diemnsion One
\vskip 0.30cm
\noindent
C.3.2. Chevalley's Linear Compacity
\vskip 0.30cm
\noindent
C.3.3. Adelic Approach in Geometric Dimension Two
\vfill
\eject
\vskip 0.45cm
\centerline {\we A. Representation of Galois Group, Stability and Tannakian
Category}
\vskip 0.30cm
\noindent
{\li A.1. Summary}
\vskip 0.30cm
Our Program for a non-abelian class CFT may be roughly summarized in the
following table.
$$\vbox{\tabskip=0pt \offinterlineskip
\def\tablerule{\noalign{\hrule}}
\halign to450pt{\strut#& \vrule#\tabskip=1em plus 2em&
\hfil#&\vrule#& \hfil#\hfil&\vrule#&
\hfil#& \vrule#\tabskip=0pt\cr\tablerule
&&\multispan5\hfil Non-Abelian Class Field Theory: A
Program\hfil&\cr\tablerule
&& Galois Aspect \quad  && Principal
&&Bundle Aspect\quad\quad&\cr\tablerule
&&Geometric Reps\quad&&Narasimhan-Seshadri Correspondence&&
S. Stable Parabolic Bundles&\cr\tablerule
&&$\Downarrow$ Rationality\qquad&&Vanishing of Brauer
Groups$\Downarrow$
&&Rationality$\Downarrow$\qquad\quad&\cr\tablerule
&&Geo-Ari Reps\qquad&&Harder-Narasimhan Correspondence&&Geo-Ari
Bundles\qquad&\cr\tablerule
&&$\Downarrow\otimes$\qquad\qquad&&$KR^2$
Trick&&$\otimes\Downarrow$\qquad\qquad&\cr\tablerule
&& Tannakian Category&&Tannaka Duality, van Kampen Cplt Th&&Clone
Tannakian Category&\cr\tablerule
&&Galois Group\qquad&&Reciprocity
Map&&$Aut^\otimes$\qquad\quad\ \ &\cr\tablerule
&&Finite Quotient\quad&&Existence Theorem, Reciprocity
Law&&Finitely Completed Module&\cr\tablerule \noalign{\smallskip}
&\multispan7\hfil\cr}}$$
\vskip 0.30cm
\noindent
{\li A.2. Non-Abelian CFT for Function Fields over C}
\vskip 0.30cm
\noindent
{\bf A.2.1. Weil-Narasimhan-Seshadri Correspondence}
\vskip 0.30cm
\noindent
{\bf A.2.1.1. Unitary Representations of Fundamental Groups}
\vskip 0.30cm
Let $M^0$ be a punctured Riemann surface of signature $(g,N)$ with $M$  the
smooth compactification. Then,
$M^0=M\backslash\{P_1,\dots,P_N\}$, and  $M$ is  of genus $g$ with
$P_1,\dots,P_N$ pairwise distinct points on $M$. Suppose that $2g-2+N>0$. 
 From the
uniformization theorem, $M^0$ can be represented as a quotient
$\Gamma\backslash {\cal H}$ of the upper half plane ${\cal H}=\{z\in {\bf
C}:{\rm Im}z>0\}$ modulo an action of a torsion-free  Fuchsian
group $\Gamma\in PSL_2({\bf R})$, generated by
$2g$ hyperbolic transformations $A_1,B_1,\dots,A_g,B_g$ and $N$ parabolic
transformations $S_1,\dots,S_N$ satisfying a single relation
$$A_1B_1A_1^{-1}B_1^{-1}\dots A_gB_gA_g^{-1}B_g^{-1}S_1\dots
S_N=1.$$ Denote the fixed points, the so-called cusps, of the parabolic 
elements
$S_1,\dots,S_N$  by $z_1,\dots,z_N$ respectively. Then, images of the cusps
$z_1,\dots,z_N\in {\bf R}\cup\{\infty\}$ under the projection $p:{\cal
H}^*:={\cal H}\cup{\bf R}\cup\{\infty\}\to\Gamma\backslash{\cal H}^*=M$
result the punctures $P_1,\dots,P_N\in M$. For each $i=1,\dots,N$, denote by
$\Gamma_i=\Gamma_{z_i}$ the stablizer of $z_i$ in $\Gamma$. Then $\Gamma_i$ is
a  cyclic subgroup in $\Gamma$ generated by
$S_i$. Moreover, for an element
$\sigma_i\in PSL_2({\bf R})$ such that $\sigma_i\infty=z_i$, we have
$\sigma_i^{-1}S_i\sigma_i=
\left(\matrix {1&\pm 1\cr 0&1\cr}\right)$, and hence $<
\sigma_i^{-1}S_i\sigma_i>=\Gamma_\infty$. (For simplicity, from now on, 
when only a
local discussion is involved, we always assume that $z_i=\infty$.)

For a representation $\rho:\pi_1(M^0)\simeq \Gamma\to GL(n,{\bf C})$  of
$\Gamma$ into a complex vector space $V$,  the vector bundle ${\bf V}:={\cal
H}\times V$ on ${\cal H}$ admits a natural $\Gamma$-vector bundle structure
via $\gamma(z,v)=(\gamma(z),\rho(\gamma)v)$ for
$\gamma\in
\Gamma, z\in {\cal H}$ and $v\in V$. The quotient of ${\cal H}\times V$ modulo
the action of $\Gamma$ is then a vector bundle of rank $n$ over $M^0$.
Moreover, since the same representation  $\rho$  defines also a $\Gamma$-vector
bundle structure on ${\cal H}^*\times V$, we obtain  a vector bundle
$V_\rho$ on $M=\Gamma\backslash{\cal H}^*$ as well.

Next assume that
$\rho$ is unitary. Then with respect to a suitable basis of
$V$,  $$\rho(S_i)={\rm diag}\Big(\exp(2\pi i\alpha_{i1}),\dots,\exp(2\pi
i\alpha_{i,n})\Big)$$ where
$\alpha_{ij}\in [0,1)$ for all $i=1,\dots,N, j=1,\dots,n$. Hence,  $V_\rho$,
or better, the associated sheaf $p_*^\Gamma({\bf V})$ of sections  may
be interpreted as follows: On $M^0$, it corresponds to the
$\Gamma$-invariant sections of ${\bf V}$, while near parabolic punctures
$P=P_i\in M$, over a neighbourhood $U$ of $P$ of the form
${\cal H}_\delta/\Gamma_\infty$ where ${\cal H}_\delta:=\{z=x+iy:y>
\delta>0\}$, the sections are all bounded
$\Gamma_\infty$-invariant sections of {\bf V} on ${\cal H}_\delta$.  Thus, 
as an
${\cal O}_{M,p}$-module, a basis of $p_*^\Gamma({\bf V})$ at $P$ is given by
the $\Gamma_\infty$ sections $\theta_j:z\mapsto \exp(2\pi i\alpha_jz)e_j$
where $\{e_1,\dots,e_n\}$ is a basis of $V$ such that
$S_i(e_j)=\exp(2\pi i\alpha_{ij})e_j$.

As a direct consequence, in addition to the associated bundles $V_\rho$ on $M$,
there exist as well the following  structures on the fibers of $V_\rho$ at 
punctures
$P_1,\dots,P_N$: Over $P=P_i$, we obtain real numbers
$$\alpha_{i1}=\alpha_{i2}=\dots=\alpha_{ik_1}<\alpha_{i,k_1+1}=\alpha_{i,k_1+2}
=\dots=\alpha_{ik_2}<\dots=\alpha_{i,k_{r_i}}$$ and a decreasing flag of
$V_\rho|_P$ defined by

\noindent
(i) $F_1(V_\rho|_P):=V_\rho|_P$;

\noindent
(ii) $F_2(V_\rho|_P)$  the subspace spanned by
$\theta_{k_1+1}\dots,\theta_n$;

\noindent
(iii) $F_3(V_\rho|_P)$  the subspace spanned by
$\theta_{k_2+1}\dots,\theta_n$, etc.

\noindent
Clearly $k_1=\dim F_1(V_\rho|_P)-\dim
F_2(V_\rho|_P),\dots, k_r=\dim F_r(V_\rho|_P)$, and
these additional structures are indeed determined by 
$\alpha_{ij}':=\alpha_{ik_j},
j=1,\dots, r_i$,  and $k_j$'s.
\vskip 0.30cm
\noindent
{\bf Proposition.} (Seshadri) {\it With the same notation as above, we
have
${\rm deg}(V_\rho)=-\sum_{i,j=1}^{N,r_i}k_{ij}\alpha_{ij}'$.}
\vskip 0.30cm
\noindent
{\bf A.2.1.2. Semi-Stable Parabolic Bundles}
\vskip 0.30cm
Following Seshadri, by definition, a parabolic structure on a vector bundle 
$E$ over a
compact Riemann surface is given by the following data:

\noindent
(1) a finite collection of points $P_1,\dots,P_N\in M$; and for each $P=P_i$,

\noindent
(2) a flag $E_P=F_1E_P\supset F_2E_P\dots\supset F_rE_P$; and

\noindent
(3) a collection of  parabolic weights $\alpha_1,\dots,\alpha_r$ attached to
$F_1E_P,\dots,F_rE_P$ such that $0\leq
\alpha_1<\alpha_2<\dots<\alpha_r<1$.

\noindent
Often $k_1=\dim F_1E_P-\dim
F_2E_P,\dots, k_r=\dim F_rE_P$ are called the multiplicities of
$\alpha_1,\dots,\alpha_r$; and a bundle $E$ together with a paraboluc structure
$$\big(P=P_i; E_P=F_1E_P\supset F_2E_P\dots\supset
F_rE_P;\alpha_1=\alpha_{i1},\dots,\alpha_r=\alpha_{ir_i}\big)_{i=1}^N$$ is 
called
a parabolic bundle and is written as
$$\Sigma(E):=\Sigma:=\Big(E; \big(P=P_i; E_P=F_1E_P\supset F_2E_P\dots\supset
F_rE_P;\alpha_1=\alpha_{i1},\dots,\alpha_r=\alpha_{ir_i}\big)_{i=1}^N\Big).$$
Trivially, if
  $W$ is a subbundle of $E$, then  $\Sigma$  induces a natural parabolic 
structure
$\Sigma(W)$ on $W$.

For parabolic bundles, its associated parabolic degree is defined  to be $${\rm
para\,deg}(\Sigma):={\rm 
deg}E+\sum_{i=1}^N\big(\sum_{j=1}^{r_i}k_{ij}\alpha_{ij}\big).$$
So, in particular,  if $\Sigma$ is induced from a unitary
representation of fundamental group of a Riemann surface,  its
associated para degree is zero. By
definition, a parabolic bundle
$\Sigma$ is called (Mumford-Seshadri) semi-stable (resp. stable) if for any 
subbundle
$W$ of $E$,
$${{{\rm para\,deg}(\Sigma(W))}\over {{\rm rank}(W)}}
\qquad \leq \ \ ({\rm resp}.\ <)\qquad {{{\rm
para\,deg}(\Sigma(E))}\over {{\rm rank}(E)}}.$$

\noindent
{\bf Proposition.} (Seshadri) {\it With the same notation as above, if $\Sigma$
is induced from a unitary representation of the corresponding
fundamental group, then $\Sigma$ is  semi-stable. Moreover, if the
representation is irreducible, then $\Sigma$ is in fact stable.}
\vskip 0.30cm
\noindent
{\bf A.2.1.3. Weil-Narasimhan-Seshadri Correspondence: A Micro Reciprocity Law}
\vskip 0.30cm
The real suprising result is the inverse of Proposition 2.1.2.
The starting point for all this is the following classical result on line 
bundles:
Over a compact Riemann surface, a line bundle is of degree zero if and only
if it is flat, i.e., it is induced from a representation of fundamental group.
It is Weil who first generalized this to vector bundles in his fundamental
paper on Generalisation des fonctions abeliennes dated in 1938: Over a compact
Riemann surface, a vector bundle is of degree zero if and only if it is 
induced from a
representation of fundamental group. (The reader may find a modern proof in
Gunning's Princeton lecture notes.)

It is said that Weil's primitive motivation is to develop a non-abelian CFT 
for Riemann
surfaces. Granting this, clearly, the next step is to study what happens
for general Riemann surfaces, which need not to be compact. This then leads 
to Weil
and Toyama's theory on matrix divisors.
\vskip 0.30cm
While all this seems to be essentially in the right direction, still many 
crucial points
are missing in these earlier studies.

Recall that  the reciprocity law in CFT is essentially the one for finite 
field extensions  and
that for a finite group, in any equivalence class of (finite dimensional 
complex)
representations there always exists a unitary one. Thus naturally we should 
strengthen
Weil's theorem from any representation to that of unitary representation. 
For doing so, the
first difficulty appears in algebraic side. That is, how to give a 
corresponding algebraic
condition?! This is solved by Mumford with his famous intersection 
stability. In fact, not only
Mumford introduces  the intersection stability, he also
studies the associated deformation theory via his fundamental work on GIT 
stability.

On the other hand, for Riemann surfaces, fundamental
groups may be described very precisely, and hence deformations of
the associated unitary representations may be quantitatively studied. All 
this, together
with certain completeness for both unitary representations and Mumford's
intersection stability, the so-called Langton's Principal, then leads 
Narasimhan and Seshadri
to prove that over a compact Riemann surface,  Mumford semi-stable
vector bundles of degree zero are naturally associated with  {\it unitary} 
representations of
the fundamental group of the surface. Later on, Seshadri first generalizes 
this result
to $\pi$-bundles, and  then to parabolic bundles.
\vskip 0.30cm
\noindent
{\bf Theorem.} (Seshadri) {\it There is a natural one-to-one
correspondence between  isomorphism classes of unitray
representations of fundamental groups of punctured Riemann surfaces
and equivalence classes of semi-stable parabolic bundles of parabolic degree
zero.}
\vskip 0.30cm
We would like call this result the Weil-Narasimhan-Seshadri correspondence.
Due to the fact that it reveals an intrinsic relation
between fundamental groups and vector bundles, often
we  call it a micro reciprocity law as well. (Over higher dimensional 
compact K\"ahler
manifold, similar correspondence is named the Kobayashi-Hitchin 
correspondence, which is
established by Ulenberk-Yau. See also Donaldson for projective manifolds.)

Note also that, as stated above, a proof  is based on

\noindent
(a) Geometric Invariant Theory;

\noindent
(b) Deformation of Reresentations; and

\noindent
(c) completeness  of semi-stable parabolic bundles --the Langton principal.

\noindent
Consequently, all this is  supposed to play a crucial role in our program 
for non-abelian CFT,
the details of which will be  discussed later.
\vskip 0.30cm
While the Narasimhan-Seshadri correspondence is a kind of Reciprocity Law, 
it is only a
micro one. Thus to find the genuine one, we need to study it globally. 
Thus, the following
result for finite groups naturally enters into the picture: Any finite 
group is determined by
its characters and vice versa. Hence, if we  have  a similar result for 
more general
groups, we could establish our non-abelian CFT for Riemann surfaces. So
naturally, we are led to  the so-called Tannakian category. But before that 
we still
need to make sure that coverings and parabolic bundles are under control. 
It is for this
purpose, we introduce the Rationality discussion in our program.
\vskip 0.30cm
\noindent
{\bf A.2.2. Rationality: Geo-Ari Representations and Geo-Ari Bundles}
\vskip 0.30cm
\noindent
{\bf A.2.2.1. Branched Coverings of Riemann Surfaces}
\vskip 0.30cm
The advantage in developing a non-abelian CFT for Riemann surfaces is that
all the time we have concrete geometric models ready to use. As there is no 
additional
cost and also for later discussion on non-abelian CFT for higher dimensional
function fields over complex numbers, we next recall some basics for 
branched coverings
of complex manifolds.

Let $M$ be an $n$-dimensional connected complex manifold. A branched
covering of $M$ is by defintion an $n$-dimensional irreducible normal
complex space $X$ together with a surjective holomorphic mapping
$\pi:X\to M$ such that

\noindent
(1) every fiber of $\pi$ is discrete in $X$;

\noindent
(2) $R_\pi:=\{q\in X:\pi^*:{\cal O}_{\pi(q),M}\to {\cal O}_{q,X}\ {\rm is\ not\
isomorphic}\}$ and $B_\pi:=\pi(R_\pi)$ are hypersurfaces of $X$ and
$M$ respectively. As uausl, we call $R_\pi$ and $B_\pi$ ramification locus 
and branch
locus respectively;

\noindent
(3) $\pi:X\backslash \pi^{-1}(B_\pi)\to M\backslash B_\pi$ is an
unramified covering; and

\noindent
(4) for any $p\in M$, there is a connected open
neighbourhood $W$ of $p$ in $M$ such that for every connected component
$U$ of $\pi^{-1}(W)$, (i) $\pi^{-1}(p)\cap U=\{q\}$ is one point and (ii)
$\pi_U:=\pi|_U:U\to W$ is surjective and proper.
Thus in particular, the induced map
$\pi:X\backslash \pi^{-1}(B_\pi)\to M\backslash B_\pi$ is a topological
covering and $\pi_U$ is finite.

For example, if $\pi:X\to M$ is a surjective
proper finite holomorphic map,  $\pi$ is a finite branched covering of
$M$.
\vskip 0.30cm
Two branch coverings $\pi:X\to M$ and $\pi':X'\to M$ are said to be
equivalent if there is a biholomorphic map $\phi:X\to X'$ such that
$\pi=\pi'\circ \phi$. In this case we write $\pi\geq \pi'$ or $\pi'\leq \pi$.
The set of all automorphisms of $\pi$ forms a group
$G_\pi$ naturally. One chacks easily that if
we denote by $\pi_1$ the restriction of $\pi$ to $X\backslash
\pi^{-1}(B_\pi)$, then $G_\pi$ is canonically isomorphic to $G_{\pi'}$.
By definition, $\pi$ is called a Galois covering  if $G_\pi$
acts transitively on every fibre of $\pi$; and $\pi$ is called abelian if
$\pi$ is Galois and $G_\pi$ is abelian.
\vskip 0.30cm
\noindent
{\bf Theorem.} {\it If $\pi:X\to M$ is a Galois covering, then

\noindent
(1) for every subgroup $H$ of
$G_\pi$, there is a branched covering $\pi_H:X/H\to M$ such that
$\pi_H\leq\pi$;

\noindent
(2) the correspondence $H\to\pi'=\pi_H$ gives a bijection
between subgroups $H$ and equivalence classes of branched coverings
$\pi'$ of $M$ such that $\pi'\leq \pi$; and

\noindent
(3) $H$ is normal if and only if
$\pi_H$ is a Galois covering, for which $G_{\pi_H}$ is siomorphic to 
$G_\pi/H$.}
\vskip 0.30cm
Note that for any branched covering $\pi:X\to M$, if $p,q$ are points of
$B_\pi$ and $\pi^{-1}(B_\pi)$ respectively such that

\noindent
(i) $B_\pi$ is normally
crossing at $p$, $q\in \pi^{-1}(p)$;

\noindent
(ii) $X$ is smooth at $q$; and

\noindent
(iii) $\pi^{-1}(B_\pi)$ is normally crossing at $q$,

\noindent
then, for a sufficiently small
connected open neighbourhood of $p$ with  a coordinate system
$(w_1,\dots,w_n)$ such that $p=0$ and $B_\pi\cap
W=\{(w_1,\dots,w_n):w_k\dots w_n=0\}$ for some $k$,  there is a
coordinate system $(z_1,\dots,z_n)$ in the connected component $U$ of
$\pi^{-1}(W)$ with $q\in U$ such that $q=0$, $\pi^{-1}(B_\pi)\cap
U=\{(z_1,\dots,z_n)\in U:z_k\dots z_n=0\}$ and
$\pi_U(z_1,\dots,z_n)=(z_1,\dots,z_{k-1}, z_k^{e_k},\dots,z_n^{e_n})$. 
Often we call
$e_j$, $j=k,\dots,n$ the ramification index of the irreducible $C_j$ of
$\pi^{-1}(B_\pi)$ such that $C_j\cap U=\{(z_1,\dots,z_n):z_j=0\}$.

Suppose now that $B_\pi=D_1\cup\dots\cup D_N$ is the irreducible
decomposition of $B_\pi$. Let $D=\sum_{i=1}^N e_iD_i$ be an effective 
divisor on $M$.
By definition, the branched covering $\pi:X\to M$ is called branched at $D$
(resp. at most at $D$) if for every irreducible component $C$ of
$\pi^{-1}(B)$ with $\pi(C)=D_j$, the ramification index of $\pi$ at $C$ is 
$e_j$
(resp. divides $e_j$). In particular,  a branched at $D$ Galois covering 
$\pi:X\to
M$  is called maximal if for any  branched covering $\pi':X'\to M$ which 
branches
at most at $D$, $\pi\geq \pi'$.
\vskip 0.30cm
While, in general, it is  complicated to describe  maximal branched 
covering for
higher dimensional complex manifolds (see however [Kato] and [Namba]), for
Riemann surfaces, it may be simply stated as follows.
 
A result of Bundgaard-Nielsen-Fox says that there is no finite
Galois covering $\pi:X\to M$ branched at $D=\sum_{j=1}^Ne_jP_j$ if and only
if either (i) $g=0$ and $N=1$ or (ii) $g=0$ $N=2$ and $e_1\not= e_2$. Here
we write $D_j$ as $P_j$, $j=1,\dots,N$. Thus from now on, we always assume that
we are not in these exceptions. Also we assume that $e_j\geq 2$. (Otherwise, we
may omit it from the beginning as there is no ramification for the 
corresponding points
$P_j$ then.)

Let $J(D)$ be the smallest normal subgroup of
$\pi_1(M^0)$ which contains $S_1^{e_1},\dots,S_N^{e_N}$. Then (as used in
the proof of the above result of Bundgaard-Nielsen-Fox,) the normal
subgroup $J(D)$ satisfies the following condition:

\centerline{Condition (*): \hskip 1cm If $S_j^d\in J(D)$, then
$d|e_j$ for $j=1,\dots,N$.}

\noindent
As a direct consequence, we have the following
\vskip 0.30cm
\noindent
{\bf Theorem.} (Bundgaard-Nielsen-Fox) {\it With the same notation as above,

\noindent
(1) there is a maximal covering $\tilde\pi:\tilde
M(D)\to M$ which branches at $D$. In particular, $\tilde M(D)$ is simply
connected;

\noindent
(2) there is a canonical one-to-one correspondence between
subgroups (resp. normal subgroups) $H$ of the quotient group
$\pi_1(M^0)/J(D)$ and equivalence classes of branched coverings (resp. 
Galois coverings)
$\pi:X\to M$ which branch at most at $D$;

\noindent
(3) $\pi:X\to M$ branches at
$D$ if and only if the following condition for $K$ is satisfied (here 
H=K/J(D)):
if $S_j^d\in K$, then $d|e_j, i=1,\dots,N$.}
\vskip 0.30cm
\noindent
{\bf A.2.2.2. Geo-Ari Representations and Geo-Ari Bundles}
\vskip 0.30cm
Motivated by the above discussion, we now introduce what we call
geo-ari representations of fundamental groups.

As before, let $M$ be a compact Riemann surface of genus $g$ with  marked 
points
$P_1,\dots,P_N$. Set $M^0=M\backslash\{P_1,\dots, P_N\}$. Then, $\pi_1(M^0)$ is
generated by hyperbolic elements $A_1, B_1,\dots, A_g,B_g$ and parabolic 
elements
$S_1,\dots, S_N$ such that
$[A_1,B_1]\dots [A_g,B_g]S_1\dots S_N=1$. Fix an effective divisor
$D=\sum_{i=1}^Ne_jP_j$ on $M$. By definition, a geometric arithmetic 
representation, a
geo-ari representation for short, of the fundamental group $\pi_1(M^0)$ 
along with
$D$ is a unitary representation $\rho:\pi_1(M^0)\to U(l)$ such that

\noindent
(i) $\rho(S_i)={\rm diag}\Big(\exp(2\pi i\beta_{i1}),\dots,\exp(2\pi 
i\beta_{il})\Big)$ for all
$i=1,\dots,N$; and

\noindent
(ii) there exist  integers $\gamma_{ij}>0$ and $\delta_{ij}\geq 0$ such that

\noindent
(a) $\gamma_{ij}|e_j$;

\noindent
(b) $(\gamma_{ij},\delta_{ij})=1$; and

\noindent
(c) $\beta_{ij}={{\delta_{ij}}\over {\gamma_{ij}}}$ for all $i=1,\dots,N$ and
$j=1,\dots,l$.
\vskip 0.30cm
Parallelly, we define a geo-ari bundle on $M$ along $D$ to be a parabolic
degree zero  semi-stable parabolic bundles
$$\Sigma=\Big(E; \big(P=P_i; E_P=F_1E_P\supset F_2E_P\dots\supset
F_rE_P;\alpha_1=\alpha_{i1},\dots,\alpha_r=\alpha_{ir_i}\big)_{i=1}^N\Big)$$
such that the following conditions are satisfied:
  there exist  integers $\gamma_{ij}>0$ and $\delta_{ij}\geq 0$ such that

\noindent
(a) $\gamma_{ij}|e_j$;

\noindent
(b) $(\gamma_{ij},\delta_{ij})=1$; and

\noindent
(c) $\alpha_{ij}={{\delta_{ij}}\over {\gamma_{ij}}}$ for all $i=1,\dots,N$ and
$j=1,\dots,r_i$.

Obviously, by applying Theorem 2.1.3, we obtain the following
\vskip 0.30cm
\noindent
{\bf Theorem}$'$. {\it With the same notation as
above,  there exists a natural one to one correspondence between
isomorphic classes of geo-ari representations of $\pi_1(M^0)$ along with $D$
and (Seshadri) equivalence classes of geo-ari bundles along $D$ over
$M$.}
\vskip 0.30cm
\noindent
{\it Remark.} We call this result  the Harder-Narashimhan Correspondence, 
despite the
fact that in the situation now, i.e., over complex
numbers,  the Harder-Narasimhan correspondence is simply the direct 
consequence of
Narasimhan-Seshadri correspondence. Later we will see that when the 
constant field is
finite, such a correspondence, first established by Harder-Narasimhan, is 
based on the
vanishing of the related Brauer groups.
\vskip 0.30cm
Clearly, if $D'=\sum_{j=1}e_j'P_j$ is an effective divisor such that
$e_j|e_j'$ for all $j=1,\dots,N$, then geo-ari representations of 
$\pi_1(M^0)$ (resp.
geo-ari bundles) along with $D$ are also geo-ari representations of 
$\pi_1(M^0)$
(resp. geo-ari bundles) along with $D'$.  (Usually, we write $D|D'$.)
Thus if we denote $U(M;D)$ the category of equivalences classes of  geo-ari
representations of $\pi_1(M^0)$ along with $D$, (see, e.g., 2.3.1 below for 
a brief discussion
on categories,) and ${\cal M}(M;D)$ the category of  geo-ari bundles along 
$D$ over $M$, then
by using a result of Mehta-Seshadri, see e.g., Prop. 1.15 of [MS], we have 
the following
\vskip 0.30cm
\noindent
{\bf Proposition.} {\it With the same notation as above,  $U(M;D)$ and  ${\cal
M}(M;D)$  are equivalent abelian categories. Moreover if $D|D'$, then 
$U(M;D)$ and
${\cal M}(M;D)$ are abelian subcategories of  $U(M;D')$ and ${\cal 
M}(M;D')$, respectively.}
\vskip 0.30cm
\noindent
{\bf A.2.3.  $KR^2$-Trick and  Completed Tannakian Categories}
\vskip 0.30cm
\noindent
{\bf A.2.3.1. Completed Tannakian Category and van Kampen Completeness Theorem}
\vskip 0.30cm
Recall that a category {\bf A} consists of two sets Obj and Arr of objects and
morphisms with two associates dom and cod from morphisms to objects. Often an
arrow $f$ with dom$(f)=x$ and cod$(f)=y$ is written as $f:x\to y$. There is 
also a map called
composition $Arr\times Arr\to Arr$ defined for $(f,g)$ when 
dom($g$)=cod($f$) such
that dom($g\circ f$)=dom($f$) and cod($g\circ f$)=cod($g$),
$f\circ(g\circ h)=(f\circ g)\circ h$, and for every object $x$ there is a 
morphism $1_x:x\to
x$ such that $f\circ 1_x=1_x\circ f=f$. Thus we may form the set $Hom(x,y)$ 
by taking
the collection of all $f:x\to y$. By definition, if moreover the hom sets 
are all abelian
groups such that compositions are bilinear, we call it a preadditive
category.

Among two categories, a functor  $T:A\to B$ is defined to be a pair of
maps $Obj(A)\to Obj(B)$ and $Arr(A)\to Arr(B)$ such that if
$f:c\to c'$ is an morphism in $A$, then $T(f)$ is a  morphism $T(f):T(c)\to 
T(c')$
in $B$, and that $T(1_c)=1_{T(c)}$, $T(g\circ f)=T(g)\circ T(f)$; and a natural
transformation $\tau$ bewteen two functors $S,T:A\to B$ is defined to be  a
collection of morphisms  $\tau_c:T(c)\to S(c)$  in $B$ such that for any 
$f:c\to
c'$ in $A$, $\tau_c T(f)=S(f)\tau_{c'}$; if moreover all $\tau_c$ have
inverses, then the natural transformation is called a functorial
isomorphism.

Among categories, abelian categories are of special importance.
By definition, an abelian category is a preadditive category such that

\noindent
(1) there is a unique object called zero object such that it is the initial 
as well as the final
object of the category;

\noindent
(2) direct product construction exists; and

\noindent
(3) associated to any
morphism are kernel and cokernel which are objects of the
category as well; moreover, every monomorphism is the kernel of its 
cokernel, every
epimorphism is the cokernel of its kernel; and

\noindent
(4) every morphism can be factored into an epimorphism followed by a 
monomorphism.

To facilitate ensuing discussion, we introduce some new concepts in 
category theory.
By definition,  an object $x$
in an abelian category is called decomposible if there exist objects, $y,z$ 
different from
zero, such that $x=y\oplus z$; and  $x$ is called irreducible if it is not 
decomposible.
Moreover, an abelian subcategory of an abelian category if called  completed if

\noindent
(1) for any object $x$, there is a unique finite decomposition $x=\oplus x_i$
with $x_i$ irreducible,  the irreducible components of $x$;
 
\noindent
(2) the subcategory contains all of its irreducible components of its objects.
\vskip 0.30cm
\noindent
{\bf Proposition.} {\it With the same notation as above, the
categories $U(M;D)$ and ${\cal M}(M;D)$ are  completed
abelian categories.}
\vskip 0.30cm
The reader may prove this proposition from the following facts:

\noindent
(1) Among two stable parabolic bundles of the same parabolic
degree,  homomorphisms are either zero or an isomorphism;

\noindent
(2) There exist Jordan-H\"older fitrations for parabolic
semi-stable bundles;

\noindent
(3) The associated Jordan-H\"older graded
parabolic bundles in (2) is unique.
\vskip 0.30cm
\noindent
{\it Remark.} Motivated by this proposition, the reader may
give a more abstract criterion to check when an  abelian subcategory is
completed.
\vskip 0.30cm
Next, let us recall what a tensor category should be. By definition, a 
cetagory is called a
tansor category if there is an operation, called tensor product $A\otimes 
B$ for any two
objects $A,B$ of the category, such that the following conditions are 
satisfied:

\noindent
(1) There are natural isomorphisms $S:A\otimes B\to B\otimes
A$ and $T:(A\otimes B)\otimes C\to A\otimes (B\otimes C)$;

\noindent
(2) $S$ and $T$ satisfy the so-called pentagon and hexagon axioms;

\noindent
(3) There is a unique identity object $1$ such that $A\simeq
A\otimes 1\simeq 1\otimes A$ for all object $A$.

Clearly, for tensor categories, we may interduce the so-called tensor
operation for any finite number of objects. Usually, there
are many ways to do so, but  the above conditions
for tensor category implies that these different ways are all the same.
Moreover, if for any
objects $x,y$ of $A$, $Hom(x,y)$ is again an object in $A$ satisfies the 
following
conditions, we  call $A$ a rigid tensor category:

\noindent
(4) there exists a morphism, the evaluation map,
$ev_{x,y}:Hom(x,y)\otimes x\to y$ such that for any object
$t$ and any morphism $g:t\otimes x\to y$, there exists a
unique morphism $f:t\to Hom(x,y)$ such that the following
commutative diagram commutes
  $$\matrix {t\otimes x&\buildrel f\otimes 1_x\over\to&
Hom(x,y)\otimes x\cr
\downarrow g&\downarrow ev_{x,y}&\cr
y&&\cr}$$
\noindent
(5) There exists natural isomorphism $$Hom(x_1,y_1)\otimes
Hom(x_2,y_2)\simeq Hom(x_1\otimes x_2,y_1\otimes y_2)$$
which is compactible with (4);

\noindent
(6) For any object $x$, by (5), if we set
$x^\vee:=Hom(x,1)$, then there exists a natural isomorphism
$x\simeq ({x^\vee})^\vee$.
\vskip 0.30cm
Now we are ready to introduce Tannakian categories.
By definition, a Tannakian category is a category which is both an abelian 
category and a
rigid tensor category such that the tensor operation is bilinear. Clearly 
then in such
categories, $Hom(x,y)\simeq x\otimes y^\vee$. For example, the
category $Vec_k$ of finite dimensional vector spaces over a field $k$ is a
Tannakian category. By definition, a functor $\omega:A\to Vec_k$ from a
Tannakian category $A$ to the category $Vec_k$ is called a
fiber functor, if $\omega$ is an exact faithful tensor
functor. (Recall that faithful means there is  a natural
injection $\omega(Hom_A(x,y))\hookrightarrow
Hom_{Vec_k}(\omega(x),\omega(y))$ which is induced from
$\omega$.)

Associated to a fiber functor $\omega$ is naturally its automorphic group
$Aut^\otimes\omega$. In particular, then we have the following fundamental 
theorem of
Tannakian category:
\vskip 0.30cm
\noindent
{\bf Theorem.} (Tannaka, Grothendieck,  Saavedra Rivano) {\it With the same 
notation as
above, assume that $(A,\omega:A\to Ver_k)$ consists of a Tannakian category
and a fiber functor such that the field $k$
is canonically isomorphic to $1_A^\vee$, then $A$ is
equivalent to the category of representations of the
group $Aut^\otimes\omega$.}
\vskip 0.30cm
A Proof of this theorem may be deduced from the following facts:

\noindent
(1) Category of all representations forms naturally a
Tannakian category together with a fiber functor;

\noindent
(2) A knowledge of the representations of a group is
equivalent to a knowledge of the group; and

\noindent
(3) Any Tannakian category is in fact a clone of (1), via the
so-called Tannaka Duality Principal.
\vskip 0.30cm
Now we are ready to introduce our own Tannakian categories which are 
completed and
equipped with natural fiber functors. There are two of them, i.e., the one 
for geo-ari
representations and the one for geo-ari bundles. Unlike what we did before,
for convinence, here we make no clear distinctions between them.
 
By Proposition 2.2.2,  it is enough to introduce the fiber functors and 
show that
geo-ari representations and geo-ari bundles are closed under the tensor 
operation. But all
this is quite straightforward:  By definition,  tensors of unitary 
representations are again
uintary, while for geo-ari bundles, the fiber functor may be defined by 
taking  special
fibers for the bundles at any point which is not a marked one. Note that
morphisms between geo-ari bundles are either zero or isomorphisms, so the 
latest defined
functor is faithful. Therefore, we obtain the following
\vskip 0.30cm
\noindent
{\bf Key Proposition.} {\it With the same notation as above, $U(M;D)$ and hence
${\cal M}(M;D)$ are completed  Tannakian categories equipped with natural 
fiber functors
to $Ver_{\bf C}$.}
\vskip 0.30cm
\noindent
{\it Remark.} While it is enough for the purpose to develop a non-abelian 
CFT for
Riemann surfaces to use the above standard theory of Tannakian categories,
in general, such a theory for Tannakian category (over base fields) is
not adequate. See ??? for details.
\vskip 0.30cm
We end this brief discussion on Tannakian category by the so-called
van Kampen completeness theorem, which will be used in the proof of the 
fundamental
theorem of CFT for Riemann surfaces later.
 
For any group $G$, denote its associated Tannakian category of equivalence 
classes
$[\rho]$  of unitary representations $\rho:G\to {\rm Aut}(V_\rho)$ by $U(G)$.
Fix for all classes $[\rho]$ a representative $\rho$ once and for all. A 
subset $Z$
of $U(G)$ is said to contain sufficiently many representations if for any 
two distinct
elements $g_1,g_2$ of $G$, there exists $[\rho]$ in $Z$ such that 
$\rho(g_1)\not=\rho(g_2)$.
\vskip 0.30cm
\noindent
{\bf van Kampen Completeness Theorem.} {\it With the same notation as above, if
$G$ is compact, then, as a completed Tannakian category, $U(G)$ may be 
generated by any
collection of objects which contains sufficiently many representations.}
\vskip 0.30cm
\noindent
{\bf A.2.3.2.  $KR^2$-Trick}
\vskip 0.30cm
While it is quite nice to have a proof of Key Proposition 2.3.1, we are by 
no means
satisfied with it. The proposition consists of two aspects: the algebraic 
one and the
analytic one. So, what we should do is the follows;
 
\noindent
(1) From analytic point of view, to prove that

  (a)   the tensor operation is
closed;  and

(b)  the so-called forgetful functor is faithful; and

\noindent
(2) From algebraic point of view, to show that

(a)  the tensor operation is
closed;  and

(b) the functor introduced in 2.3.1 is faithful.

\noindent
However in our proof outlined above, only
(1.a) and (2.b) are shown. That is to say,
with the help of the so-called Narasimhan-Seshadri
correspondence, a micro reciprocity law, we make no distinction between 
algebraic and
analytic aspects.  Thus, a purely analytic proof for (1.a)
and a purely  algebraic proof for (2.a) should be pursued.
 
The proof of
(1.a) is a simple one since  unitary representations for our fundamental 
groups, or
better the unitary representations of quotients of fundamental groups by 
(normal)
subgroups generated by weighted parabolic generators have a 
semi-simplification. We
leave the details for the reader. Thus the real challenge
here is an algebraic proof of (2.a), i.e., an proof for that tensor products of
geo-ari bundles are again geo-ari bundles. It is here we should use another 
central
concept in GIT, the so-called instability flag of Mumford-Kempf (with the
refined version given by Ramanan and Ramanathan).  This goes as follows.

In his study of GIT stabilities, Mumford conjectures that if a  point is 
not semi-stable,
then there should exist a  parabolic subgroup which takes reponsibility, in the
view of the so-called Hilbert-Mumford criterion. This is  confirmed by 
Kempf. (In
Kempf's study, as suggested by Mumford, the rationality problem is also treated
successfully, at least when constant fields are perfect.) Kempf's result 
then motivates
Ramanan and Ramanathan  to show that even though, for the original
action, the corresponding point is not semi-stable, but if a certain type 
well-controlled
modification is allowed (with the aim to cancell the instability 
contribution in the original
action),  a new (yet well-associated)  point could be constructed  such 
that with respect to
the natural induced action this new point becomes  semi-stable. As a direct 
consequence, in
the case for semi-stable vector bundles without parabolic structures, since 
the new
well-associated action may be associated to the intersection stability 
condition for bundles
naturally,  we may then obtain an algebraic proof of (2.a) for bundles.
\vskip 0.30cm
\noindent
{\it Remark.} Over complex numbers, the rationality is not a serious 
problem. However,
over finite fields, this turns to be a  difficult one. As a matter of fact, 
we then need to use
the Frobenius to tackle the rationality problem. For details, see ???.
\vskip 0.30cm
Therefore, to develop a non-abelian CFT for function fields over complex 
numbers,  we
need to find a more general version, what we call the $KR^2$-trick, of the 
above result of
Kempf, and Ramanan and Ramanathan, which works for parabolic bundles. For 
this purpose,
we may follow a more down-to-earth approach of Faltings and Tataro. That is 
to say,
following Faltings, we first write  any  parabolic sub-bundle of the tensor 
product  in terms
of  filtrations over some points, disjoint from parabolic ones; then 
introduce a certain GIT
stability for general filtrations to check whether  the induced filtrations 
from subbundles
and the  original parabolic filtration  are GIT stable. If so, we by 
definition are done.
Otherwise, following Totaro, from the associated instability flag of 
Mumford-Kempf, we can
complete the proof by using the intersection stablility of each of the 
components of the
tensor product.
\vskip 0.30cm
\noindent
{\bf A.2.4. Non-Abelian CFT for Function Fields over Complex Numbers}
\vskip 0.30cm
\noindent
{\bf A.2.4.1. Micro Reciprocity Law, Tannakian Duality and the Reciprocity Map}
\vskip 0.30cm
Let $M$ be a compact Riemann surface of genus $g$ with marks $P_1,\dots,P_N$.
Set $M^0:=M\backslash\{P_1,\dots,P_N\}$. Then, the fundamental group
$\pi_1(M^0)$ is generated by $2g$ hyperboilic generators 
$A_1,B_1,\dots,A_g,B_g$ and $N$
parabolic generators $S_1,\dots,S_N$ which satisfy one single relation
$\prod_{i=1}^g[A_i,B_i]\cdot\prod_{j=1}^N S_j=1.$
Moreover, with respect to an effective divisor
$D=\sum_{i=1}^Ne_iP_i$,  we have a Tannakian category
$U(M;D)$ consisting of equivalent classes of unitary representations
$[\rho:\pi_1(M^0)\to {\rm Aut}(V_\rho)]$ of $\pi_1(M^0)$ such that

\noindent
(i) $\rho(S_i)={\rm diag}\Big(\exp(2\pi i\beta_{i1}),\dots,\exp(2\pi 
i\beta_{il})\Big)$ for all
$i=1,\dots,N$; and

\noindent
(ii) there exist  integers $\gamma_{ij}>0$ and $\delta_{ij}\geq 0$ such that

\noindent
(a) $\gamma_{ij}|e_j$;

\noindent
(b) $(\gamma_{ij},\delta_{ij})=1$; and

\noindent
(c) $\beta_{ij}={{\delta_{ij}}\over {\gamma_{ij}}}$ for all $i=1,\dots,N$ and
$j=1,\dots,l$.
 
\noindent
Now in each equivalence class $[\rho]$, fix a unitary representation, 
denoted also by
$\rho$, once and hence for all. Clearly, $\rho$ induces a unitary 
representation  of the
group
$\pi_1(M^0)/J(D)$, where $J(D)$ denotes the normal subgroup of
$\pi_1(M^0)$ generated by $S_1^{e_1},\dots,S_N^{e_N}$.
Call this latest representation $\rho$ as well. Then, for any element $g\in 
\pi_1(M^0)/J(D)$,
$\rho(g)$ induces for each object $[\rho:\pi_1(M^0)\to {\rm Aut}(V_\rho)]$ 
in $U(M;D)$ an
automorphism of $V_\rho$. As a direct consequence, we obtain a natural 
group morphism
from
$\pi_1(M^0)/J(D)$ to the automorphism group of the corresponding fiber 
functor $U(M;D)\to
Vec_{\bf C}$.

On the other hand, for the  Tannakian category ${\cal M}(M;D)$
of geo-ari bundles on $M$ along $D$ together with the fiber functor 
$\omega(M;D):{\cal
M}(M;D)\to Ver_{\bf C}$, by Tannakian Duality, we conclude that
$\omega(M;D):{\cal M}(M;D)\to Ver_{\bf C}$ is equivalent to the Tannakian 
category of the
representations of ${\rm Aut}^{\otimes}\omega$. Therefore, by the
Narasimhan-Seshadri correspondence (and the Harder-Narasimhan correspondence),
the so-called micro reciprocity law, we obtain a canonical group morphism
$$\Omega(D):\pi_1(M^0)/J(D)\to {\rm Aut}^{\otimes}(\omega(M;D)).$$
We will call $\Omega(D)$  the reciprocity map associated with $(M,D)$.
\vskip 0.30cm
\noindent
{\bf A.2.4.2.  Non-Abelian CFT}
\vskip 0.30cm
By definition, a subcategory $S$ of a Tannakian category with respect to a 
fiber functor
$\omega$ is called  a {\it finitely completed
Tannakian subcategory}, if

\noindent
(1) it is  a completed Tannakian subcategory;

\noindent
(2) there exist finitely many objects which generated $S$ as an abelian 
tensor subcategory;

\noindent
(3) ${\rm Aut}^\otimes(\omega|_S)$ is a finite group.

With this,   by using the
van Kampen completeness theorem and  the above reciprocity map, we then can 
manage
to obtain the  following fundamental theorem on non-abelian class field 
theory for
function fields over complex numbers.
\vskip 0.30cm
\noindent
{\bf Fundamental Theorem in Non-Abelian Class Field Theory for Riemann 
Surfaces.}
\vskip 0.30cm
\noindent
(1) (Existence and Conductor Theorem) {\it There is a natural one-to-one 
correspondence
$\omega_{M,D}$ between
$$\{{\bf S}:{\rm finitely\ completed\ Tannakian\ subcategory\ of} \ {\cal
M}(M;D)\}$$ and $$
\{\pi:M'\to M:{\rm finite\ Galois\ covering\ branched\ at\ most\ at}\
D\};$$}

\noindent
(2) (Reciprocity Law) {\it There is a natural group isomorphism}
$${\rm Aut}^\otimes (\omega(M;D)\big|_{\bf S})\simeq {\rm 
Gal}\,(\omega_{M,D}({\bf
S})).$$

We end this discussion by pointing out that, as an application, one may use 
this
fundamental theorem to solve the geometric inverse Galois problem.
\vskip 0.30cm
\noindent
{\bf A.2.5. Classical (abelian) CFT: An Example of Kwada-Tata and Kawada}
\vskip 0.30cm
\noindent
{\bf A.2.5.1. Class Formation}
\vskip 0.30cm
Classical abelian class field theory was first formulated axiomatically in 
Artin-Tate
seminar in terms of cohomology of groups. As an example of this 
formulation, later
Kawada-Tate and Kawada studied function fields over complex numbers. To make a
comparison between this classical approach of CFT and what we outlined 
above, next we
recall the works in Kawada-Tate and   Kawada.
\vskip 0.30cm
Let $k_0$ be a given ground field and $\Omega$ a fixed infinite separable 
normal
algebraic extension of $k_0$. Let ${\cal R}$ be the set of all finite 
extensions of
$k_0$ in $\Omega$. By definition we call a collection $\{E(K):K\in {\cal 
R}\}$ of
abelian groups $E(K)$ a formulation if the following conditions are satisfied:

\noindent
F1. If $k\subset K$ then there is an injective morphism
$\phi_{k/K}:E(k)\hookrightarrow E(K)$;

\noindent
F2. If $k\subset l\subset K$, then $\phi_{l/K}\circ\phi_{k/l}=\phi_{k/K}$;

\noindent
F3. If $K/k$ is normal and $G={\rm Gal}(K/k)$ is its Galois group, then $G$ 
acts
on $E(K)$ and $\phi_{k/K}(E(k))=E(K)^G$.

\noindent
F4. If $k\subset L\subset K$  and $L/k,K/k$ are both normal, then the
Galois group $F={\rm Gal}(L/k)$ is a quotient group of $G={\rm Gal}(K/k)$. 
Denote
by $\lambda_{G/F}:G\to F$ the canonical quotient map. Then for any
$\sigma\in G$ and $f\in E(L)$, $\sigma\circ \phi_{L/K}(f)=\phi_{L/K}\circ
(\lambda_{G/F}\sigma)(f)$.

Moreover if a formation is called a class formation if it satisfies the 
following additional
conditions on group cohomology:

\noindent
C1: $H^1(G,E(K))=0$; and

\noindent
C2: $H^2(G,E(K))\simeq [K:k_0]{\bf Z}$.
\vskip 0.30cm
By a theorem of Tate,  C1 and C2 imply and hence are equivalent to
the following stronger condition

\noindent
C. In a class formation, for all $r\in {\bf Z}$, $H^r(G,E(K))\simeq 
H^{r-2}(G,{\bf
Z})$ for all $r$.  In particular, if a 2-cocycle $\alpha(K)$ generates the 
cyclic group
$H^2(G,E(K))$, then the cup-product $g\mapsto\alpha(K)\cup g$
induces the isomorphism.
\vskip 0.30cm
Note that by definition,  $H^{-2}(G,{\bf Z})=G^{\rm ab}:=G/[G,G]$, and
the 0-th cohomology group $H^0(G,E(K))$ is nothing
but $E(K)^G/T_GE(K)$, (where $T_G a=\sum_{\sigma\in G}\sigma(a)$). So for a 
class
formation, we obtain then the reciprocity law
$$E(k)/\phi^{-1}(T_GE(K))\simeq G^{\rm ab}.$$

Furthermore, for a class formation, with respect to the so-called
res (restriction), infl (inflation) and ver (Verlagerung) operations of 
group cohomology,
we have the follows;

\noindent
Case 1. For $k\subset l\subset K$ with $K/k$ normal,  $G:={\rm Gal}(K/k)$ and
$H:={\rm Gal}(K/l)$,  in $H^2(G,E(K)), ver_{H/G}\circ res_{G/H}=[G:H]\cdot 1.$
So, $res_{G/H}$ is surjective and $ver_{H/G}$ is injective;

\noindent
Case 2. For $k\subset L\subset K$ with $L/k$, $K/k$ normal, and $G:={\rm
Gal}(K/k)$, $H:={\rm Gal}(K/L)$,  then we have the exact
sequences
$$0\to H^2(F,E(L))\buildrel infl\over\to H^2(G,E(K))\buildrel res\over\to
H^2(H,E(K)).$$

Therefore, there exist 2-cocycles $\{\alpha(K)\}$ of $G(K/k)$ over $E(K)$ 
such that

\noindent
D1. $\alpha(K/k)$ are generators of the cyclic groups $H^2(G,E(K))$;

\noindent
D2. In Case 1, $$res_{G/H}(\alpha(K/k))\sim f(K/l);
\qquad ver_{H/G}(\alpha(K/k)\sim [G:H]\cdot \alpha(K/k);$$

\noindent
D3. In Case 2, $infl_{F/G}\alpha(L/k)\sim [K:L]\cdot \alpha(K/k).$ Here $\sim$
means cohomologous.

Often, we call such a system $\{\alpha(K/k)\}$ the canonical 2-cocycle of
$G(K/k)$ over $E(K)$ which may be used to write down the reciprocity law as
follows: Introduce $({{K/k}\over\sigma})\in E(k)$,
$\sigma\in G:={\rm Gal}(K/k)$ for normal extension $K/k$ by $$({{K/k}\over
\sigma})=\phi^{-1}_{k/K}\Big(\prod_{\rho\in
G}\alpha_{K/k}(\rho,\sigma)\Big)\quad{\rm mod}\ \phi_{k/K}^{-1}(T_GE(K/k)).$$
Then by a result of Nakayama, the symbol $(a,K/k)\in G^{\rm ab}$,
$a\in E(k)$ defined by
$$(a,K/k):=\sigma\ {\rm mod}\, [G,G]\qquad{\rm when}\qquad ({{K/k}\over
\sigma})=a\quad{\rm mod}\,\phi_{k/K}^{-1}(T_GE(K)),$$   satisfies all the 
properties of the
norm residue symbol in number theory.
\vskip 0.30cm
Let $k\in {\cal R}$ be fixed and $\Omega^a(k)$ be the maximal abelian extension
of $k$ in $\Omega$. Set ${\cal R}^a(k)$ be the collection of finite abelian 
extensions of
$k$ (in $\Omega^a$). Then for $K\in {\cal R}^a(k)$ define a subgroup 
$A(K/k)$ of
$E(k)$ to be $\phi_{k/K}^{-1}(T_{\rm Gal}(K/k)E(K))$. By definition, a subgroup
$F$ of $E(k)$ is called admissible if $F=A(K/k)$ for some $K\in {\cal R}^a(k)$.
Denote the set of all admissible subgroups of $E(k)$ by ${\cal U}(E(k))$. 
Then one
checks, by the properties of norm residue symbol, that

\noindent
(1) (Combination Theorem) $A(K_1\cdot K_2/k)=A(K_1/k)\cap A(K_2/k)$, $A(K_1\cap
K_2/k)=A(K_1/k)\cdot A(K_2/k)$;

\noindent
(2) (Ordering Theorem) $A(K_1/k)\supset A(K_2/k)$ if and only if $K_1\subset
K_2$;

\noindent
(3) (Uniqueness Theorem) $A(K_1/k)=A(K_2/k)$ if and only if $K_1=K_2$
for $K_1,K_2\in {\cal R}^a(k)$.
\vskip 0.30cm
Also, if  let $\Gamma(k)$ be the compact Galois group $\Omega^a(k)/k$, then
$\Gamma(k)$ is the inverse limit group of $\{G(K/k):K\in {\cal R}^a(k)\}$. 
Moreover, for $a\in
E(k)$, the limit, called the generalized norm residue symbol, 
$(a,k):=\lim_{K\in {\cal
R}^a(k)}(a,K/k)\in \Gamma(k)$ exists. Set ${\cal T}(k):=(E(k),k)\subset
\Gamma(k)$ and ${\cal R}(k):=\{(a\in E(k):(a,k)=1\}\subset E(k)$, then the
mapping $a\mapsto (a,k)$ induces an isomorphism $$E(k)/{\cal R}(k)\simeq {\cal
T}(k)\subset \Gamma(k),$$ where ${\cal T}(k)$ is dense in $\Gamma(k)$.
\vskip 0.30cm
For examples, classical (abelian) CFT for global fields are all class 
formations.
More precisely,

\noindent
(A) For a number field $k$,  we may take $E(k)$ to be the idele class group
$C_k$  and  $\phi_{k/K}$ the natural inclusion.  In particular, ${\cal 
R}(k)$ is then
the connected component of the unity of $C_k$, $\Gamma(k)={\cal T}(k)$, and
${\cal U}(E(k))$ is the set of all open subgroup of $E(k)$ of finite index;

\noindent
(B)  For a function field $k$ of one variable over a finite field, we  may take
$E(k)$ to be the idele class group $C_k$ and $\phi_{k/K}$ the natural 
inclusion.
In particular, ${\cal R}(k)=1$, $\Gamma(k)/{\cal T}(k)$ is a uniquely 
divisible group
isomorphic to $\hat {\bf Z}/{\bf Z}$, and ${\cal U}(E(k))$ is the
set of all open subgroup of $E(k)$ of finite index.
\vskip 0.30cm
\noindent
{\bf A.2.5.2. The Work of Kawada and Tate}
\vskip 0.30cm
Choose $k_0$ to be an algebraic function field of one variable over complex
numbers with $\Omega$ the algebraic closure of $k_0$. Let $D(k)$ be the 
group of
all fractional divisors $\prod_vP_v^{r_v}$, where $r_v\in {\bf Q}$, and
$r_v=0$  for almost all $v$, and $P(k)$ the group of all principal 
divisors. Let ${\bf
D}(k)=D(k)/P(k)$.

Define   $E(k):={\bf D}(k)^\vee={\rm Hom}({\bf D}(k),{\bf
R}/{\bf Z})$, the group of characters of ${\bf D}(k)$, then  $\{E(k)\}$ 
with the conorm map
$\phi_{k/K}$ is a class formation such that the norm residue map 
$\Phi_k:E(k)\to A(k)$ is
surjective but $F(k):=Ker\Phi_k\not=0$.
As a direct consequence, if $T{\bf D}(k)$ denotes the torsion subgroup of 
${\bf D}(k)$,
$\{E(k)^*=E(k)/F(k)=(T{\bf D}(k))^\vee\}$ gives also a class formation.
\vskip 0.30cm
On the other hand, easily,

\noindent
(i) the character of $k_0$ is zero;

\noindent
(ii) $k_0$ contains  all the roots of unity; and

\noindent
(iii) for any finite normal extension $K/k$ with $k\supset k_0$ and $K/k$ 
finite, $N_{K/k}K=k$.

\noindent
So, by applying the so-called Kummer theory,  $E(k):=(k^*\otimes {\bf 
Q}/{\bf Z})^\vee$,
  and the conorm
$\phi_{k/K}:E(k)\to E(K)$, i.e., $\phi_{k/K}(\chi)(A)=\chi(N_{K/k}A)$
for $\chi\in E(k),A\in K^*\otimes {\bf Q}/{\bf Z}$, give a class formation with
$E(k)=A(k)$.

Thus, in particular, by comparing these class formations, we obtain a
canonical isomorphism $k^*\otimes {\bf Q}/{\bf Z}\simeq T{\bf D}(k).$
Moreover, if $\Omega_\phi$ is a maximal unramified extension of $k_0$,
$t:=t(\Omega_\phi/k)$ and $E(k)_\phi:=({\bf D}_\phi(k))^\vee$ with ${\bf 
D}_\phi(k)$ the
divisor class group of the usual sense, i.e., with integral coefficients, 
then $\{E_\phi(k)\}$
gives a class formation for $t_\phi$.  On the other hand, for a fixed 
finite set
$S\not=\emptyset$ of prime divisors of $k_0$, let
$\Omega_S$ be the maximal $S$-ramified extension of $k_0$. Put
$t_S:=t(\Omega_S/k_0)$ and $E(k)_S:=({\bf D}_S(k))^\vee$ with ${\bf 
D}_S(k)$ the
$S$-fractional divisor class group of $k$. Then $\{E_S(k)\}$ is a class 
formation for
$t_S$, such that $Im\,\Phi_k=A(k)$ and $Ker\,\Phi_k$ is the connected
component of $E(k)_S$ which is not zero. Similarly, if we take 
$E(k)_S^*:=(T{\bf
D}_S(k))^\vee$, then $\{E(k)_S^*\}$ forms also a class formation for $t_S$ 
such that
$E(k)_S^*\simeq A(k)$.
\vskip 0.30cm
\noindent
{\bf A.2.5.3. Abelian CFT for Riemann Surfaces In Terms of Geo-Ari Bundles}
\vskip 0.30cm
However, it is via the third approach of a class formation for Riemann 
surfaces due to
Kawada, these classical results are related with our approach to the
CFT.  It goes as follows. For $k_0\subset k$, let $S(k)$ denote the set of 
prime divisors
of
$k$ which are extensions of a prime divisor of $k_0$ contained in $S$. Let 
$R(k)$ be the
Riemann surface of $k$ and $R_S(k)=R(k)\backslash S(k)$. Let $E(k)$ be the one
dimensional integral homology group
$H_1(R_S(k),{\bf Z})$ of $R_S(k)$, and define $\phi_{k/K}:E(k)\to E(K)$  by
$\gamma\mapsto V\gamma$ where $V\gamma$ is the covering path of
$\gamma\in R_S(k)$ on the unramified covering surface $R_S(K)$ of $R_S(k)$.
Then  $\{E(k)\}$ forms again a class formation for $t_S$ with
$Im\,\Phi_k$ dense in $A(k)$ and $Ker\,\Phi_k=0$. Indeed, one may obtain 
this  by
looking at the canonical pairing, the micro reciprocity law in this context,
$H_1(R_S(k),{\bf Z})\times T{\bf D}_S(k)\to {\bf S}^1\subset {\bf C}$ 
defined by
$(\eta,A)_S\mapsto \exp(\int_\eta d\log A)$ for 1-cycle $\eta$ on $R_S(k)$ and
$A\in T{\bf D}_S(k)$. Here $d\log A$ denotes the abelian differential of 
the third
kind on $R_S(k)$ corresponding to a divisor $A$.

On the other hand,  over a compact Riemann surface  $M$ with punctures 
$P_1,\dots,P_N$
with respect to an effective divisor $D=\sum e_jP_j$, we may introduce the 
group
${\rm Div}_0^{\bf Q}(M;D)$ of degree zero {\bf Q}-divisors along $D$ on $M$ 
by collecting
all degree zero {\bf Q}-divisors of the form $\sum_j{{a_j}\over 
{e_j}}P_j+E$ with
$a_j\in {\bf Z}$ and $E$ a ordinary integral divisor on $M$. Denote the 
induced
(rational equivalence) divisor class group  ${\rm Cl}_0^{\bf Q}(M;D)$. Then,
${\rm Cl}_0^{\bf Q}(M;D)$ is simply the collection of all geo-ari bundles 
of rank 1
introduced in 2.2.2. Hence Theorem 2.4.2 then implies the following
\vskip 0.30cm
\noindent
{\bf Theorem.} {\it With the same notation as above, there is a one-to-one
correspondence between the set ot all isomorphism classes of finite abelian 
coverings
$\pi:X\to M$ branched at most at $D$ and the set of all finite subgroups $S$ of
${\rm Cl}_0^{\bf Q}(M;D)$. Moreover, the correspondence $\pi\mapsto S(\pi)$ 
satisfies
that

\noindent
(i) (Reciprocity Law) $S(\pi)\simeq G_\pi$; and

\noindent
(ii) (Ordering Theorem) $\pi\leq \pi'$ if and only if $S(\pi)\subset S(\pi')$.}
\vskip 0.30cm
Therefore, it seems to be very crucial to
understand the precise relation between Kawada-Tate's results and this 
latest theorem, in
order to develop a (non-abelian) CFT for global fields.
\vskip 0.30cm
\noindent
{\li A.3. Towards Non-Abelian CFT for Global Fields}
\vskip 0.30cm
\noindent
{\bf A.3.1. Weil-Narasimhan-Seshadri Type Correspondence}
\vskip 0.3cm
\noindent
{\bf A.3.1.1. Grometric Representations}
\vskip 0.30cm
 From what discussed above, in order to develop a non-abelian CFT
for local and  global fields, the first step should be the one to establish a
micro reciprocity law, i.e., a Weil-Narasimhan-Seshadri type correspondence.
Therefore, we are supposed to

\noindent
(1) introduce suitable classes of  representations of Galois groups;

\noindent
(2) find corresponding classes for bundles in terms of intersection;

\noindent
(3) establish natural correspondences between classes in (1) and (2).
 
Hence, it is then more practicle to divide the problem into two.
Namely, a general one in the sense of Weil Correspondence,
and a refined one in the sense of Narasimhan-Seshadri correspondence.

We start with the Weil Correspondence. Here, we  are then supposed  to first
introduce a general notion for representation of Galois groups such that, 
naturally,
associated to such representations are special vector bundles together with 
additional
structure.

As an example, let me  explain what I have in mind in the case for number
fields. So, let $F$ be a number field with a finite subset $S$ of places of 
$F$. Denote the
corresponding Galois group by $G_{F,S}$. Naturally, by a representation, it 
should be
first a continuous group homolorphism $\rho:G_{F,S}\to {\rm GL}_n({\bf 
A}_{F,S})$, where
${\bf A}_{F,S}$ denotes the ring of $S$-adeles. Moreover, amony others, we 
should
assume that

\noindent
(a) for all  places $v$ of $F$, the induced representations
$\rho^v:G_{F,S}\to {\rm GL}_n(F_v)$ are unramified almost everywhere;

\noindent
(b) for a fixed places $p$ of {\bf Q}, there should be a compactibility 
condition
for all  places $v$ above $p$.

In addition, it also seems to be very natural to assume that

\noindent
(c) for all induced $\rho_v:G_v\to {\rm GL}_n(F_v)$,  there are invariant 
lattices $M_v$ of
$F_v$ which are induced from a global lattice $M$ over $F$; and

\noindent
(d) at places $v\in S$, there should be naturally a certain weighted filtration
induced by the action of Frobenius.

Thus in particular, associated to such a representation,  is naturally a 
well-defined vector
bundle equipped with a parabolic structure.
We are expecting that geometric bundles are of Arakelov degree zero.

Based on the Weil correspondence, we should be able to enter the level of
Narasimhan-Seshadri type correspondence. Here essential difficulties should 
appear.
Chiefly, what should be a natural analog of being unitary? A suitable 
candinate seems to be
that of Fountaine's semi-stablility at finite places and (unitary) at 
infinite places.
However, we are not very sure about this, as somehow we believe that the 
condition of
semi-stable representation is too restricted.  For this reason, we propose 
the follows;
 
\noindent
(e) for all places $v$, the images of the induced representation 
$\rho_v:G_v\to {\rm
Gl}_n(F_v)$  are contained in maximal compact subgroups; moreover, certain
compactibility conditions are satified by $\rho_v$ for all $v$ over a fixed 
place of {\bf Q};

\noindent
(f)  there should be a natural deformation
theory for these representations such that

\noindent
(i) the size of all equivalence classes of these representations can be 
controlled;

\noindent
(ii) a certain completeness holds.
\vskip 0.30cm
Clearely, here the reader would find the work of Rapoport and Zink [RZ] and 
the lecture
notes of Tilouine [Ti] are of great use. The cases for local fields and 
function fields may be
similarly discussed.  If success, we call such a representation a geometric 
representation.
\vskip 0.30cm
\noindent
{\bf A.3.1.2. Semi-Stability in terms of Intersection}
\vskip 0.30cm
This is the algebraic side of the micro reciprocity law. We here only study
what happens for function fields over finite fields, while leave a detailed 
defintion for
number fields  in Part (B).

With this restriction of fields, the situation becomes much simpler: We may 
use the
existing stability condition of Seshadri for parabolic bundles.

That is to say, what we care here are the so-called parabolic semi-stable 
bundles of
parabolic degree zero defined on algebraic curves over finite fields, 
introduced by Seshadri.
However, in doing so, we are afraid that  our program leads only a 
non-abelian CFT with
tame ramifications. So it seems that, to deal with wild ramifications, 
additional works
are needed.  Therefore, in gereral, for this algebraic aspect, we propose 
the follows:

\noindent
(0) Once we have an  intersection semi-stabilities,  the analog of
the existence and uniqueness of Harder-Narasimhan filtration, the existence of
Jordan-H\"older filtration and the uniqueness of the graded Jordal-H\"older
objects should hold; moreover, morphisms of stable objects should be either 
zero or
isomorphisms;

\noindent
(1) The so-called Langton's  Completeness Principal should hold for such 
intersection
semi-stabilities;

\noindent
(2) The intersection semi-stability should be naturally related with a GIT 
stability.
As a direct consequence, then moduli spaces can be formed naturally, and by
(1), are indeed  compact;

\noindent
(3) With the help of the Frobenius, we should be able to
define a special subset (of points) of  moduli spaces, such that
tensor products of the resulting bundles  are represented again
by  points in this special subset.
\vskip 0.30cm
If success, we will call these objects geometric bundles.
\vskip 0.30cm
\noindent
{\bf A.3.1.3. Weil-Narasimhan-Seshadri Type Correspondence}
\vskip 0.30cm
For a Weil type correspondence, from the classical proof, say, the one given in
Gunning's Princeton lecture notes, we  should develop an analog of
the de Rham, the Dolbeault cohomologies as well as a Hodge type theory. 
While, in general,
it is  out of reach, note that our base is of dimension one,  it is still 
hopeful.
For example, for curves over finite fields, we understand that all this is 
known to
experts.

Next, let us consider a Narasimhan-Seshadri type correspondence. Here, we 
should
establish  a natural correspondence between the equivalence classes of 
geometric
representations and the (Seshadri) equivalence classes of geometric bundles.
Hence, the key points are the follows:

\noindent
(1) By definition, a geometric representation should naturally give a geometric
bundle;

\noindent
(2) A Weil type correspondence holds. In particular, this implies that, by 
definition,
geometric bundles come naturally from representations of Galois groups;

\noindent
(3) Representations resulting from geometric bundles in (2) should be
geometric. To establish this, we should use the deformation theories of 
geometric
representations and geometric bundles as proposed in 3.1.1 and 3.1.2 above.
For example, the final justification should be based on the compactness
of both geometric representations and geometric bundles, via a direct counting.
(In geometry, the counting is possible since we know the structure of the 
fundamental
group. But in arithmetic, this is quite difficult. For this, we find that 
the work of Fried and
V\"olklein [FV] seems to be useful.)
\vskip 0.30cm
\noindent
{\bf A.3.2. Harder-Narasimhan Correspondence}
\vskip 0.30cm
This is specially designed to select special subclasses of geometric 
representations and
geometric bundles so as to get what we call geo-ari representations and 
geo-ari bundles.
There are two main reasons for doing so. The first is based on the facts 
that for function
fields, there are examples of stable bundles whose tensor products  are no 
longer
semi-stable; and more importantly that if the bundles and their associated 
Frobenius twists
are all  semi-stable,  then,   so are their tensor
products.  The second comes from our construction of non-abelian zeta 
functions in
Part (B). To define such non-abelian zeta functions, we  use only
rational points (over constant fields) of the moduli spaces, based on a 
result of Harder and
Narasimhan, which guarantees the coincidence of the rationality of moduli 
points and
the rationality of the corresponding bundles.
 
Thus, for example, for function fields over finite fields, a geo-ari bundle 
is defined to a
geometric bundle which satisfies not only all the conditions for geo-ari 
bundles over
Riemann surfaces, but an additional one, which says that all its Frobenius 
twists are
semi-stable as well.

  So for our more general purpose,  key points here are the follows;

\noindent
(1) to give a proper definition of geo-ari representations so that they 
form a natural
abelian category;

\noindent
(2) to give a suitable definition of geo-ari bundles so that they are 
closed under
tensor product; and

\noindent
(3) to establish a natural correspondence between (1) and (2). If success, 
we call such a
result  a Harder-Narasimhan type corespondence.
\vskip 0.30cm
\noindent
{\bf A.3.3. $KR^2$-Trick}
\vskip 0.30cm
This is specially designed to give an algebraic proof of the following
key statement appeared as 3.2.(2): The tensor products of geo-ari bundles 
are again
geo-ari bundles.
 
We start with geo-ari bundles over function fields. When no parabolic 
structures is
involved, this is solved by Kempf and  Ramanan-Ramanthan.  In their proof,
key points are the follows.

\noindent
(1) Rationality and uniqueness of instability flags, where the Frobenius 
twists are used;

\noindent
(2) Existence of  a GIT stable modification associated to any non 
semi-stable one; and

\noindent
(3) Existence of GIT points for bundles such that GIT stability implies 
intersection
stability by definition;

\noindent
(4) Relation between the GIT modification in (2) and the intersection 
stability, where
the interesction stability for each component of the tensor product plays a 
key rule.
 
More precisely, if the tensor is not intersection stable, by (3), the 
associated GIT point
would not be GIT stable. Thus from (2) there exists a GIT stable 
modification which is well
understood by (1). Therefore, finally, by (4), i.e., the intersection 
stability of each
components,  we may finally conclude the intersection stability for the 
tensor product with
the help of (2).
\vskip 0.30cm
Thus, the problem left here is to see how the work of Ramanan-Rananthan 
could be
developed to deal with parabolic structure. In theory, this may be done
by working on product of varieties instead of just a single
one. As this process is more or less similar to what happens in 
constructing moduli
space of parabolic semi-stable bundles,  an experct should be able
to carry the details out.

However, there is a more elementary approach as well,  essentially given in the
supplementary works of Faltings and Totaro. It goes as follows.  First, as what
Faltings does, write parabolic subbundles for the tensor in terms of 
weighted filtrations
on the fibers (supported over points which are disjoint from punctures);
then introduce a GIT stability with respect to weighted filtrations. So 
there are two
possibilities: if the weighted filtration resulted from the parabolic 
subbundle is GIT
stable, then by (3) above, we get the intersection stability of the tensor 
product.
Otherwise, by (1) and (2) above, we may obtain a well-associated 
modification which is
then GIT stable. This then by applying (4) completes the proof.

Before going to number fields, I would like to draw the  attention of the 
reader
to a relation between the stability of the generic point and the stability 
of special
points, due to Mumford: GIT stability for objects over the  generic point 
implies that
over almost all but finitely many special points, the associated points are 
GIT  semi-stable.
\vskip 0.30cm
Now we come to number fields. For this, first, we  should develop an 
Arakelov style
GIT. This sounds difficult, but due to the following fact we expect that it 
is still workable:
Kempf's functional related to the Hilbert-Mumford criterion
for GIT stability is essentially compactible with Arakelov theory. (See 
e.g., [Bu] and
[Zh].)
 
Based on such a new GIT, then we should find a relation between the 
intersection stability
in the definition of geo-ari bundles  and this new type of GIT stability. 
In particular,
we expect that an analog of $KR^2$-trick works.  Thus the final key point 
should be

\noindent
(5) An Arakelov style GIT  exists such that (1),  (2),  (3) and (4) above 
work equally
well under the framework of this new GIT over number fields.
\vskip 0.30cm
Note also that not only for the purpose of $KR^2$-trick, a new GIT is 
needed, for the
construction of moduli spaces of geometric bundles, such a new GIT is 
supposed to be crucial.
For the related point, see also ??? in Part (B).
\vskip 0.30cm
\noindent
{\bf A.3.4. Tannakian Category Theory over Arbitary Bases}
\vskip 0.30cm
While for a non-abelian class field theory of function fields over complex 
numbers,
the standard theory of Tannakian category  may be applied directly, I do not
expect this in general. For example, for number fields, we need to develop 
a more general
theory, a theory of Tannakian categories over the ring of integers.
The key  points for this new type of Tannakian categories are supposed to 
be the follows.

\noindent
(1) Fiber functors should give us a group. Thus whether fiber functors are
to categories of vector spaces over a field is not really important. In 
fact, for number
fields,  fiber functors should be faithful exact tensor
functor to the categories of finitely generated projective modules of the
corresponding ring of integers, among others;

\noindent
(2) An analog of Tannakian duality holds; and

\noindent
(3) Tannakian category should contain the so-called finitely completed 
Tannakian
subcategories with respect to which an analog of van Kampen
Completeness Theorem holds.
\vskip 0.30cm
If success, as a direct consequence,  by (1) and (2), using 
Narasimhan-Seshadri  and
Harder-Narasimhan correspondences,  we  get naturally  a reciprocity map. 
This then
by (3), leads to a completed non-abelian CFT.
\vskip 0.30cm
\noindent
{\bf A.3.5. Non-Abelian CFT for Local and Global Fields}
\vskip 0.30cm
All in all, from the above discussion, what we then expect is the following 
conjecture, or
better,
\vskip 0.30cm
\noindent
{\bf Working Hypothesis.}  {\it For local and  global fields,

\noindent
(1) there are well-defined geometric
representations and geometric bundles such that a
Weil-Narasimhan-Seshadri type correspondence holds;

\noindent
(2) there are refined well-defined geo-ari representations and geo-ari bundles
such that a Harder-Narasimhan type correspondence holds;

\noindent
(3) there are well-defined GIT type stability  such that the
intersection stability as appeared in the definition of geo-ari bundles
could be understood in terms of this new GIT stability.
Moreover, an analog  of $KR^2$-trick works;

\noindent
(4) there is a well-established Tannakian type category theory, for which a 
Tannaka
type duality and van Kampen type completeness theorem hold; moreover the
category of geo-ari objects forms naturally such Tannakian type categories.

\noindent
In particular, the fundamental results in non-abelian class field
theory, such as  existence theorem, conductor theorem and
reciprocity law, hold.}
\vfill\eject
\vskip 0.30cm
\centerline {\we B. Moduli Spaces, Riemann-Roch, and New Non-Abelian Zeta
functions}
\vskip 0.45cm
\noindent
{\li B.1. New Local and Global Non-Abelian Zeta Functions for Curves}
\vskip 0.30cm
\noindent
{\bf B.1.1. Local Non-Abelian Zeta Functions}
\vskip 0.30cm
\noindent
{\bf B.1.1.1. Artin Zeta Functions for Curves}
\vskip 0.30cm
We start this program by recalling the construction and basic properties of the
classical Artin zeta functions of curves defined over finite fields.

Let $C$ be a projective irreducible reduced regular curve of genus
$g$ defined over a finite field $k:={\bf F}_q$ with $q$ elements. Then the 
arithmetic degree
$d(P)$ of a closed point $P$ on $C$ is defined to be $d(P):=[k(P):k]$, 
where $k(P)$ denotes
the residue class field of $C$ at $P$. So  $q^{d(P)}$ is nothing but the
number $N(P)$ of elements in $k(P)$. Extending this to all divisors, we 
then may define the
Artin zeta function $\zeta_C(s)$ for curve $C$ over $k$ by setting
$$\zeta_C(s):=\prod_P(1-N(P)^{-s})^{-1}=\sum_{D\geq 0}N(D)^{-s}=\sum_{D\geq
0}\big(q^{-s}\big)^{d(D)},\qquad{\rm Re}(s)>1,$$ where the sum is taken 
over all
the effective divisors $D$ on $C$. As usual, set $t=q^{-s}$, and
$Z_C(t)=\sum_{D\geq 0} t^{d(D)}.$

Note that the number of positive
divisors which are rational equivalent to a fixed divisor $D$ is
$(q^{h^0(C,D)}-1)/(q-1)$. So, $$Z_C(t)=\sum_{D\geq 0}t^{d(D)}
=\sum_{d\geq 0}\sum_{\cal D}\sum_{D\in {\cal D}}t^d=h(C)\sum_{d\geq
0}{{q^{h^0(C,D)}-1}\over {q-1}}\cdot t^d,$$ where the sum $\sum_{\cal D}$ 
is taken over the
rational divisor classes of degree
$d$, and $h(C)$ denotes  the cardinality of ${\rm Div}_0(C)/{\rm 
Div}_p(C)$, the group of
degree zero divisors modulo the subgroup of principal divisors.  Thus by 
Riemann-Roch
theorem,
   for a positive
divisor $D$, $h^0(C,D)\leq d(D)+1$. So up to  finitely many convergent 
terms, the
convergence of
$Z_C(t)$ is the same as that for $\sum_{d\geq 0}(d+1)(qt)^d$. Hence $\zeta
_C(s)$ is well-defined.

Being well-defined, $$(q-1)\zeta_C(s)=\sum_{d\geq 0}\sum_{\cal D}(q^{h^0({\cal
D})}-1)q^{-ds} =\sum_{d\geq 0}\sum_{\cal D}q^{h^0({\cal D})-ds}-{{h(C)}\over
{1-q^{-s}}}.$$ Set $F(q^{-s}):=\sum_{d\geq 0}'\sum_{\cal D}q^{h^0({\cal
D})-ds}-{{h(C)}\over {1-q^{-s}}}$ and $G(q^{-s}):=\sum_{d\geq 0}''\sum_{\cal
D}q^{h^0({\cal D})-ds},$ where
$\sum'$ and $\sum''$ are taken over $d\geq 2g-2+1$ and $0\leq d \leq 2g-2$ 
respectively.
Then,

\noindent
(i) $(q-1)Z_C(t)=F(t)+G(t);$

\noindent
(ii) Being a sum of finite many terms, $G(t)$ is rational; and moreover, by
Riemann-Roch,

\noindent
(iii)
$F(t)=h(C){\sum}_d'q^{d-g+1}t^d-{{h(C)}\over
{1-t}}=h(C)q^{1-g}{{(qt)^{2g-2+1}}\over{1-qt}}-{{h(C)}\over
{1-t}}.$

\noindent
So,  $Z_C(t)$ is indeed a rational function of $t$. Further,  from (iii),
$F(t)=q^{g-1}t^{2g-2}F\Big({1\over {tq}}\Big).$ Note also that in $G(t)$, 
the sum
is taken the sum over all divisors whose degrees are between 0 and $2g-2$. 
Thus   by
the duality and Riemann-Roch, for any canonical divisor $K_C$ of $C$,
$$G(t)=\sum_{\cal D}"q^{h^0({\cal D})}t^{d{\cal D}}
=\sum_{\cal D}"q^{d({\cal D})-g+1+h^0(K_C-{\cal D})}t^{d{\cal D}}
=\sum_{\cal D}"(qt)^{d({\cal D})}q^{-g+1}q^{h^0(K_C-{\cal D})}
=q^{-g+1}t^{2g-2}G\Big({1\over{qt}}\Big).$$ That is to say, $\zeta_C(s)$
satisfies the functional equation. So we have the following
\vskip 0.30cm
\noindent
{\bf Theorem.} {\it With the same notation as above, $\zeta_C(s)$ is 
well-defined, admits a
meromorphic to the whole complex $s$-plane. Moreover,
$$\zeta_C(s)=N(K_C)^{1/2-s}\zeta_C(1-s)$$ and there exists a polynomial 
$P_C(t)$ of
degree
$2g$ such that
$$Z_C(t)={{P_C(t)}\over{(1-t)(1-qt)}}.$$}
\vskip 0.30cm
\noindent
{\bf B.1.1.2. Too Different Generalizations: Weil Zeta Functions and A New 
Approach}
\vskip 0.30cm
It is well-known that based on the so-called reciprocity law the above 
Artin zeta function
may be written as
$$Z_C(t)=\exp\Big(\sum_{m\geq 1}{{N_m}\over m}\cdot t^m\Big)$$ where 
$N_m:=\#C({\bf
F}_{q^m})$ denotes the number of ${\bf F}_{q^m}$-rational points of $C$. 
This then leads to
a far reaching generalization to the so-called Weil zeta functions of 
higher dimensional
varieties defined over finite fields, the study of which dominantes what we 
call Arithmetic
Geometry in the second half of 20th centery.

On the other hand, for the purpose of developing a non-abelian zeta 
function theory,
we here do it differently. In terms of Artin zeta functions, the key then 
is the following
observations;

\noindent
(1)  $\sum_{\cal D}$ in 1.1.1 may be viewed as taking summation over the
degree $d$ Picard group of the curve;

\noindent
(2) Picard groups for curves may be viewed as moduli spaces of 
(semi-stable) line
bundles;

\noindent
(3) moduli spaces of semi-stable bundles exist and all the terms appeared 
in the
summation for Artin zeta functions, such as $h^0$ and degree, make sense 
for vector
bundles as well.
\vskip 0.30cm
\noindent
{\bf B.1.1.3. Moduli Spaces of Semi-Stable Bundles}
\vskip 0.30cm
  Let $C$ be a regular, reduced and irreducible projective curve
defined over an algebraically closed field
$\bar k$. Then according to Mumford [M], a vector bundle $V$ on $C$ is 
called semi-stable (resp.
stable) if for any proper subbundle $V'$ of $V$,  $$\mu(V'):={{d(V')}\over 
{r(V')}}\leq {\rm
(resp.}< {\rm )}{{d(V)}\over {r(V)}}=:\mu(V).$$ Here $d$ denotes the degree 
and $r$ denotes the
rank.

\noindent
{\bf Proposition 1.} {\it Let $V$ be a vector bundle over $C$. Then

\noindent
(a) ([HN]) there exists a unique
filtration of subbundles of $V$, the so-called Harder-Narasimhan filtration 
of $V$,
$$\{0\}=V_0\subset V_1\subset V_2\subset\dots\subset V_{s-1}\subset V_s=V$$
such that for $1\leq i\leq s-1$, $V_i/V_{i-1}$ is semi-stable and
$\mu(V_i/V_{i-1})>\mu(V_{i+1}/V_i);$

\noindent
(b) (see e.g. [Se]) if moreover $V$ is semi-stable,
there exists a
filtration of subbundles of $V$, a Jordan-H\"older filtration of $V$,
$$\{0\}=V^{t+1}\subset V^t\subset\dots\subset V^{1}\subset V^0=V$$
such that for all $0\leq i\leq t$, $V^i/V^{i+1}$ is stable and 
$\mu(V^i/V^{i+1})=\mu (V)$.
Moreover, ${\rm Gr}(V):=\oplus_{i=0}^tV^i/V^{i+1}$, the associated
(Jordan-H\"older) graded bundle of
$V$, is determined uniquely by $V$.}
\vskip 0.30cm
  Following Seshadri, two semi-stable vector bundles $V$
and $W$ are called $S$-equivalent, if their associated Jordan-H\"older 
graded bundles are
isomorphic, i.e.,
${\rm Gr}(V)\simeq {\rm Gr}(W)$. Applying Mumford's general result on 
geometric invariant theory
([M]), Seshadri proves the following
\vskip 0.30cm
\noindent
{\bf Theorem 2.} ([Se]) {\it Let $C$ be a regular, reduced, irreducible 
projective curve  defined
over an algebraically closed field. Then over the set ${\cal M}_{C,r}(d)$ 
of $S$-equivalence
classes of rank $r$ and degree $d$ semi-stable vector
bundles over $C$, there is a natural normal, projective
  algebraic variety structure.}

\vskip 0.30cm
  Now  assume that $C$ is defined over a finite field $k$. Naturally
we may talk about
$k$-rational bundles over $C$, i.e., bundles which are defined over $k$. 
Moreover, from geometric
invariant theory, projective varieties
${\cal M}_{C,r}(d)$ are defined over a certain finite extension of $k$. 
Thus it makes sense to
talk about
$k$-rational points of these moduli spaces too. The relation between these 
two types of rationality
is given by  Harder and Narasimhan based on a discussion about Brauer groups:
\vskip 0.30cm
\noindent
{\bf Proposition 3.} ([HN]) {\it There exists a finite field ${\bf F}_q$ 
such that
for any $d$, the subset of ${\bf F}_q$-rational points of
${\cal M}_{C,r}(d)$  consists exactly of all $S$-equivalence classes of
${\bf F}_q$-rational bundles in ${\cal M}_{C,r}(d)$.}
\vskip 0.30cm
 From now on, without loss of generality, we always assume that  finite fields
${\bf F}_q$ (with $q$ elements) satisfy the property in this Proposition. 
Also for simplicity, we
write ${\cal M}_{C,r}(d)$ for ${\cal M}_{C,r}(d)({\bf F}_q)$, the subset of 
${\bf F}_q$-rational
points, and call them moduli spaces by an abuse of notations.
\vskip 0.45cm
\noindent
{\bf B.1.1.4. New Local Non-Abelian Zeta Functions}
\vskip 0.30cm
  Let $C$ be a  regular, reduced, irreducible projective curve defined over 
the finite field ${\bf
F}_q$ with $q$ elements. Define the {\it rank $r$ non-abelian zeta 
function} $\zeta_{C,r,{\bf
F}_q}(s)$ by setting
$$\zeta_{C,r,{\bf F}_q}(s):=\sum_{V\in [V]\in {\cal M}_{C,r}(d), d\geq 
0}{{q^{h^0(C,V)}-1}\over
{\#{\rm Aut}(V)}}\cdot (q^{-s})^{d(V)},\qquad {\rm Re}(s)>1.$$ Clearly, we 
have the
following

\noindent
{\bf Fact.} {\it With the same notation as above, $\zeta_{C,1,{\bf 
F}_q}(s)$ is nothing
but the classical Artin zeta function for curve $C$.}
\vskip 0.30cm
\noindent
{\bf B.1.1.5. Basic Properties for  Non-Abelian Zeta Functions}
\vskip 0.30cm
Many basic properties for classical Artin zeta functions are satisfied by 
our non-abelian
zeta functions as well. More precisely, we have the following
\vskip 0.30cm
\noindent
{\bf Theorem.}  {\it With the same notation as above,

\noindent
(1) The non-abelian zeta function $\zeta_{C,r,{\bf F}_q}(s)$ is 
well-defined for ${\rm
Re}(s)>1$, and admits a meromorphic extension to the whole complex $s$-plane;

\noindent
(2) ({\rm Rationality}) If we  set
$t:=q^{-s}$ and introduce the non-abelian $Z$-function of $C$ by setting
$$\zeta_{C,r,{\bf F}_q}(s)=:Z_{C,r,{\bf F}_q}(t):=\sum_{V\in [V]\in {\cal 
M}_{C,r}(d),d\geq
0}{{q^{h^0(C,V)}-1}\over {\#{\rm Aut}(V)}}\cdot t^{d(V)}, \qquad |t|<1,$$ then
  there exists a polynomial $P_{C,r,{\bf F}_q}(s)\in {\bf Q}[t]$ such that
$$Z_{C,r,{\bf F}_q}(t)={{P_{C,r,{\bf F}_q}(t)}\over {(1-t^r)(1-q^rt^r)}};$$

\noindent
(3) ({\rm Functional Equation}) If we set the rank $r$
non-abelian $\xi$-function
$\xi_{C,r,{\bf F}_q}(s)$ by
$$\xi_{C,r,{\bf F}_q}(s):=\zeta_{C,r,{\bf F}_q}(s)\cdot (q^{s})^{r(g-1)},$$ 
then
$$\xi_{C,r,{\bf F}_q}(s)=\xi_{C,r,{\bf F}_q}(1-s).$$}

One may prove this theorem  by using the vanishing theorem, duality, and the
Riemann-Roch theorem.  See e.g., ???  for details.
\vskip 0.30cm
\noindent
{\bf Corollary.} {\it With the same notation as above,

\noindent
(1) $P_{C,r,{\bf F}_q}(t)\in {\bf Q}[t]$ is a degree $2rg$ polynomial;

\noindent
(2) Denote all reciprocal roots of $P_{C,r,{\bf F}_q}(t)$ by 
$\omega_{C,r,{\bf F}_q}(i),
i=1,\dots, 2rg$. Then after a suitable rearrangement,
$$\omega_{C,r,{\bf F}_q}(i)\cdot \omega_{C,r,{\bf F}_q}(2rg-i)=q,\qquad 
i=1,\dots,rg;$$

\noindent
(3) For each $m\in {\bf Z}_{\geq 1}$, there exists a rational number 
$N_{C,r,{\bf F}_q}(m)$ such
that
$$Z_{r,C,{\bf F}_q}(t)=P_{C,r,{\bf F}_q}(0)\cdot\exp\Big(\sum_{m=1}^\infty 
N_{C,r,{\bf
F}_q}(m){{t^m}\over m}\Big).$$ Moreover, $$N_{C,r,{\bf
F}_q}(m)=\cases{r(1+q^m)-\sum_{i=1}^{2rg}\omega_{C,r,{\bf F}_q}(i)^m,& if 
$r\ |m$;\cr
-\sum_{i=1}^{2rg}\omega_{C,r,{\bf F}_q}(i)^m,& if
$r\not| m$;\cr}$$

\noindent
(4) For any  $a\in {\bf Z}_{>0}$, denote by $\zeta_{a}$ a primitive $a$-th 
root of unity and set
$T=t^a$.  Then
$$\prod_{i=1}^aZ_{C,r}(\zeta_{a}^it)=(P_{C,r,{\bf F}_q}(0))^a\cdot
\exp\Big(\sum_{m=1}^\infty N_{r,C,{\bf F}_q}(ma){{T^m}\over m}\Big).$$}
\vskip 0.45cm
\noindent{\bf  B.1.2. Global Non-Abelian Zeta Functions for Curves}
\vskip 0.30cm
\noindent {\bf B.1.2.1. Preparations}
\vskip 0.30cm
  Let $C$ be a regular, reduced, irreducible
projective curve of genus $g$ defined over the finite field
${\bf F}_q$ with $q$  elements. Then the rationality of $\zeta_{C,r,{\bf 
F}_q}(s)$ says that there
exists a degree $2rg$ polynomial $P_{C,r,{\bf F}_q}(t)\in {\bf Q}[t]$ such that
$$Z_{C,r,{\bf F}_q}(t)={{P_{C,r,{\bf F}_q}(t)}\over {(1-t^r)(1-q^rt^r)}}.$$
Set $$P_{C,r,{\bf F}_q}(t)=\sum_{i=0}^{2rg}a_{C,r,{\bf F}_q}(i)t^i.$$ By 
the functional equation
for $\xi_{C,r,{\bf F}_q}(t)(s)$, we have
$$P_{C,r,{\bf F}_q}(t)=P_{C,r,{\bf F}_q}({1\over {qt}})\cdot q^{rg}\cdot 
t^{2rg}.$$ So, for
$i=0,1,\dots,rg-1$,
$a_{C,r,{\bf F}_q}(2rg-i)=a_{C,r,{\bf F}_q}(i)\cdot q^{rg-i}.$

To further determine these coefficients, following  Harder and Narasimhan 
(see e.g. [HN] and
[DR]), who first consider the
$\beta$-series  invariants below, we introduce the following invariants:
$$\alpha_{C,r,{\bf F}_q}(d):=\sum_{V\in [V]\in {\cal M}_{C,r}(d)({\bf
F}_q)}{{q^{h^0(C,V)}}\over {\#{\rm Aut}(V)}},\qquad
  \beta_{C,r,{\bf F}_q}(d):=\sum_{V\in [V]\in
{\cal M}_{C,r}(d)({\bf F}_q)}{{1}\over {\#{\rm Aut}(V)}},$$
and $\gamma_{C,r,{\bf F}_q}(d):=\alpha_{C,r,{\bf F}_q}(d)-\beta_{C,r,{\bf 
F}_q}(d).$
One checks that all
$\alpha_{C,r,{\bf F}_q}(d),\beta_{C,r,{\bf F}_q}(d)$ and $\gamma_{C,r,{\bf 
F}_q}(d)$'s may be
calculated from $\alpha_{C,r,{\bf F}_q}(i),\beta_{C,r,{\bf F}_q}(j)$ with 
$i=0,\dots, r(g-1)$ and
$j=0,\dots,r-1$.
\vskip 0.30cm
\noindent
{\bf An Ugly Formula} {\it With the same notation as above,
$$\eqalign{~&a_{C,r,{\bf F}_q}(i)\cr
=&\cases{\alpha_{C,r,{\bf
F}_q}(d)-\beta_{C,r,{\bf F}_q}(d),& if\ $0\leq i\leq r-1$;\cr
\alpha_{C,r,{\bf F}_q}(d)
-(q^r+1)\alpha_{C,r,{\bf F}_q}(d-r)+q^r\beta_{C,r,{\bf
F}_q}(d-r),& if\ $r\leq i\leq 2r-1$;\cr
\alpha_{C,r,{\bf F}_q}(d)
-(q^r+1)\alpha_{C,r,{\bf F}_q}(d-r)+q^r\alpha_{C,r,{\bf
F}_q}(d-2r),& if\ $2r\leq i\leq r(g-1)-1$;\cr
-(q^r+1)\alpha_{C,r,{\bf F}_q}(r(g-2))
+q^r\alpha_{C,r,{\bf F}_q}(r(g-3))
+\alpha_{C,r,{\bf F}_q}(r(g-1)),& if\ $i=r(g-1)$;\cr
\alpha_{C,r,{\bf F}_q}(d)
-(q^r+1)\alpha_{C,r,{\bf F}_q}(d-r)+\alpha_{C,r,{\bf
F}_q}(d-2r)q^r,& if\ $r(g-1)+1\leq i\leq rg-1$;\cr
2q^r\alpha_{C,r,{\bf F}_q}(r(g-2))
-(q^r+1)\alpha_{C,r,{\bf F}_q}(r(g-1)),& if\ $i=rg$.\cr}\cr}$$}
\vskip 0.30cm
\noindent
{\bf B.1.2.2. Global Non-Abelian Zeta Functions for Curves}
\vskip 0.30cm
Let ${\cal C}$ be a regular, reduced, irreducible projective curve of genus
$g$ defined over a number field $F$. Let $S_{\rm bad}$ be the collection of 
all infinite places and
these finite places of $F$ at which ${\cal C}$  does not have good 
reductions. As usual, a place
$v$ of
$F$ is called good if $v\not\in S_{\rm bad}$.
 
Thus, in particular, for any good place $v$ of $F$,   the $v$-reduction of 
${\cal C}$, denoted as
${\cal C}_v$, gives a regular, reduced, irreducible projective curve 
defined over the residue field
$F(v)$ of $F$ at $v$. Denote the cardinal number of $F(v)$ by $q_v$.  Then, 
by 1.1, we obtain
the associated rank $r$ non-abelian zeta function
$\zeta_{{\cal C}_v,r,{\bf F}_{q_v}}(s)$. Moreover, from the rationality of 
$\zeta_{{\cal
C}_v,r,{\bf F}_{q_v}}(s)$, there exists a degree $2rg$ polynomial $P_{{\cal 
C}_v,r,{\bf F}_{q_v}}(t)\in {\bf Q}[t]$
such that
$$Z_{{\cal C}_v,r,{\bf F}_{q_v}}(t)={{P_{{\cal C}_v,r,{\bf 
F}_{q_v}}(t)}\over {(1-t^r)(1-q^rt^r)}}.$$
Clearly, $P_{{\cal C}_v,r,{\bf F}_{q_v}}(0)=\gamma_{{\cal C}_v,r,{\bf 
F}_{q_v}}(0)\not=0.$
Thus it makes sense to introduce the polynomial $\tilde P_{{\cal 
C}_v,r,{\bf F}_{q_v}}(t)$ with
constant term 1 by setting
$$\tilde P_{{\cal C}_v,r,F(v)}(t):={{P_{{\cal C}_v,r,F(v)}(t)}\over 
{P_{{\cal C}_v,r,F(v)}(0)}}.$$
Now by definition, {\it the rank $r$ non-abelian
zeta function $\zeta_{{\cal C},r,F}(s)$ of} ${\cal C}$ over $F$ is the 
following Euler product
$$\zeta_{{\cal C},r,F}(s)
:=\prod_{v:{\rm good}}{1\over{
\tilde P_{{\cal C}_v,r,{\bf F}_{q_v}}(q_v^{-s})}},\hskip 2.0cm {\rm 
Re}(s)>>0.$$

Clearly, when $r=1$, $\zeta_{{\cal C},r,F}(s)$ coincides with the classical 
Hasse-Weil zeta function
for $C$ over $F$.
\vskip 0.30cm
\noindent
{\bf Conjecture.} {\it For a regular, reduced, irreducible projective curve 
${\cal C}$ of
genus $g$ defined over a number field $F$,  its associated rank $r$ global 
non-abelian
zeta function
$\zeta_{{\cal C},r,F}(s)$  admits a meromorphic continuation to the whole 
complex $s$-plane.}

Recall that even  when $r=1$, i.e., for the classical Hasse-Weil zeta 
functions, this conjecture
is still open. However, in general, we have the following
\vskip 0.30cm
\noindent
{\bf Theorem 1.}  {\it  Let ${\cal C}$ be a regular, reduced, irreducible 
projective curve
defined over a number  field $F$.   When ${\rm Re}(s)> 1+g+(r^2-r)(g-1)$, 
the associated rank r
global non-abelian zeta function $\zeta_{{\cal C},r,F}(s)$ converges.}

This theorem may be deduced from a result of (Harder-Narasimhan) Siegel on 
Tamagawa numbers of
$SL_r$, the ugly yet very precise formula for local zeta function in 1.2.1, 
Clifford Lemma for
semi-stable bundles, and Weil's theorem on the Riemann hypothesis for Artin 
zeta
functions. In fact we have the following
\vskip 0.30cm
\noindent
{\bf Proposition 2.} {\it With the same notation as above, when $q\to\infty$,

\noindent
(a) For $0\leq d\leq r(g-1)$,
$${{\alpha_{C,r,{\bf F}_q}(d)}\over
{q^{d/2+r+r^2(g-1)}}}=O(1);$$

\noindent
(b) For all $d$, $$\beta_{C,r,{\bf
F}_q}(d)=O\Big(q^{r^2(g-1)}\Big);$$

\noindent
(c) $${{q^{(r-1)(g-1)}}\over{\gamma_{C,r,{\bf
F}_q}(0)}}=O\Big(1\Big).$$}
 
\noindent
{\bf B.1.2.3. Working Hypothesis}
\vskip 0.30cm
Like in the theory for abelian zeta functions, we want to use our
non-abelian zeta functions  to study  non-abelian aspect of
arithmetic of curves. Motivated by the classical analytic class number 
formula for Dedekind zeta
functions and its counterpart BSD conjecture for Hasse-Weil zeta functions 
of elliptic curves, we
expect that our non-abelian zeta function could be used to understand the 
Weil-Petersson volumes of
moduli spaces of stable bundles.

For doing so, we then also need to introduce  local factors for \lq bad' 
places. This may be
done as follows.  For
$\Gamma$-factors, we take these coming from the functional equation for
  $\zeta_F(rs)\cdot\zeta_F(r(s-1))$, where $\zeta_F(s)$ denotes
the standard Dedekind zeta function for $F$; while for finite bad places, 
first,
use the semi-stable reduction for curves to find a semi-stable model for 
${\cal C}$, then use
Seshadri's moduli spaces of parabolic bundles to construct polynomials for 
singular fibers, which
usually have degree lower than $2rg$. With all this being done, we then can 
introduce the so-called
completed rank
$r$ non-abelian zeta function  for
${\cal C}$ over
$F$, or better,  the  completed rank $r$ non-abelian zeta function 
$\xi_{X,r,{\cal O}_F}(s)$ for a
semi-stable model $X\to {\rm Spec}({\cal O}_F)$ of ${\cal C}$. Here ${\cal 
O}_F$ denotes the ring of
integers of $F$. (If necessary, we take a finite extension of $F$.)

\noindent
{\bf Conjecture.} {\it $\xi_{X,r,{\cal O}_F}(s)$ is holomorphic and 
satisfies the functional equation
$$\xi_{X,r,{\cal O}_F}(s)=\pm\, \xi_{X,r,{\cal O}_F}(1+{1\over r}-s).$$}

Moreover, we expect that for certain classes of curves, the inverse Mellin 
transform of our
non-abelian zeta functions are naturally associated to certain modular 
forms of weight $1+{1\over
r}$.
\vskip 0.30cm
\noindent
{\bf Example.}  For elliptic curves ${\cal E}$ defined over {\bf Q}, we
obtain the following  \lq absolute Euler product' for rank 2 zeta functions 
of elliptic
curves
$$\zeta_{{\cal E}, 2,{\bf Q}}(s)=\zeta_2(s)=\prod_{p\ {\rm prime}}{1\over
{1+(p-1)p^{-s}+(2p-4)p^{-2s}+(p^2-p)p^{-3s}+p^2p^{-4s}}}.$$

At this point, it may be better to recall the following result of
Andrianov ([An]). (We thank Kohnen for drawing our attention to this point.)
The so-called genus two spinor
$L$-function stands in the form
$$\prod_p{1\over{1-\lambda(p)p^{-s}+(\lambda(p)^2-\lambda(p^2)-p^{2k-4})
p^{-2s}-\lambda(p)p^{2k-3}p^{-3s}+p^{4k-6}p^{-4s}}}.$$ Clearly, if we set 
$k=2$,
$\lambda(p)=1-p$ and $\lambda(p^2)=p^2-4p+4$, we  see that {\it formally} 
the above two zeta
functions coincide. This suggests that there
might be a close relation between them. The following is a speculation I 
made after the
discussion with Deninger and Kohnen.
\vskip 0.30cm
\noindent
{\it The relation between the above two zeta functions should be  in the 
same style as the Shimura
correspondence for half weight and integral weight modular forms.}
\vskip 0.30cm
To convince the reader, let me point out the following facts:

\noindent
(1) Andrianov's zetas have a Hecke theory,
are coming from certain weight 2 modular forms, and have  the local factors
$${{1-\lambda(p)t+\big(\lambda(p)^2-\lambda(p^2)-p^{2k-4}\big)
t^2-\lambda(p)p^{2k-3}t^3+p^{4k-6}t^4}\over {1-p^{2k-4}t^2}};$$

\noindent
(2) Our working hypothesis concerning weight $3/2$ modular forms
are made mainly from the fact that our local factor takes the form
$${{1+(p-1)t+(2p-4)t^2+(p^2-p)t^3+p^2t^4}\over {(1-t^2)(1-p^2t^2)}},$$
in which an additional factor $1-p^2t^2$ appears in the denominator.
\vskip 0.30cm
\noindent
{\bf B.1.3. Refined Brill-Noether Locus for Elliptic Curves: Towards A 
Reciprocity Law}
\vskip 0.30cm
\noindent
{\bf 1.3.1 Results of Atiyah}
\vskip 0.30cm
Let $E$ be an elliptic curve defined over $\overline{{\bf F}_q}$,
an algebraic closure of the finite field ${\bf F}_p$ with $q$-elements.

Recall that a vector bundle $V$ on $E$ is called indecomposable if $V$ is
not the direct sum of two proper subbundles, and that every
vector bundle on $E$ may be written as a direct sum of indecomposable
bundles, where the summands and their multiplicities are uniquely
determined up to isomorphism. Thus to understand vector bundles, it
suffices to study the indecomposable ones. To this end, we have the
following result of Atiyah [At]. In the sequel, for simplicity, we always
assume that  the characteristic of ${\bf F}_q$ is strictly bigger than the
rank of $V$.

\noindent
{\bf Theorem 1.} (Atiyah) {\it (a) For any $r\geq 1$, there is a unique
indecomposable vector bundle $I_r$ of rank $r$ over $E$, all of whose
Jordan-H\"older constituents are isomorphic to ${\cal O}_E$. Moreover,
the bundle $I_r$ has a canonical filtration $$\{0\}\subset
F^1\subset\dots\subset F^r=I_r$$ with $F^i=I_i$ and $F^{i+1}/F^i={\cal
O}_E$;

\noindent
(b) For any $r\geq 1$ and any integer $a$, relative prime to $r$ and each
line bundle $\lambda$ over $E$ of degree $a$, there exists up to
isomorphism a unique indecomposable bundle $W_r(a;\lambda)$ over $E$ of
rank $r$ with $\lambda$ the determinant;

\noindent
(c) The bundle $I_r(W_{r'}(a;\lambda))=I_r\otimes W_{r'}(a;\lambda)$is
indecomposable and every indecomposable bundle is isomorphic to
$I_r(W_{r'}(a;\lambda))$ for a suitable choice of $r,r',\lambda$. Every
bundle $V$ over $E$ is a direct sum of vector bundles of the form
$I_{r_i}(W_{r_i'}(a_i;\lambda_i))$, for  suitable choices of
$r_i,r_i'$ and $\lambda_i$. Moreover, the triples
$(r_i,r_i',\lambda_i)$ are uniquely specified up to permutation by the
isomorphism type of $V$.}
\vskip 0.30cm
Here note in particular that  $W_r(0,\lambda)\simeq \lambda$, and
that indeed $I_r(W_{r'}(a;\lambda))$ is the unique indecomposable bundle
of rank $rr'$ such that all of whose successive quotients in the
Jordan-H\"older filtration are isomorphic to $W_{r'}(a;\lambda)$.
Now for a vector bundle $V$ over $E$, define its
slop $\mu(V)$ by $\mu(V):={\rm deg}(V)/{\rm rank}(V)$.
Then, $I_r(W_{r'}(a;\lambda))$ is semi-stable with
$\mu(I_r(W_{r'}(a;\lambda)))=a/r'$.
\vskip 0.30cm
\noindent
{\bf Theorem 2.} {\it (a) (Atiyah) Every bundle $V$ over $E$ is isomorphic
to a direct sum $\oplus_iV_i$ of semi-stable bundles, where
$\mu(V_i)>\mu(V_{i+1})$;

\noindent
(b) (Atiyah) Let $V$ be a semi-stable bundle over $E$ with slop 
$\mu(V)=a/r'$ where
$r'$ is a positive integer and $a$ is an integer relatively prime to
$r'$. Then $V$ is a direct sum of bundles of the form
$I_r(W_{r'}(a;\lambda))$, where $\lambda$ is a line bundle of degree
$a$;

\noindent
(c) (Atiyah, Mumford-Seshadri) There exists a natural projective algebraic 
variety
structure on  $${\cal M}_{E,r}(\lambda):=\{V:{\rm semi}-{\rm stable}, {\rm
det}V=\lambda, {\rm rank}(V)=r\}/\sim_S.$$  Moreover, if $\lambda\in {\rm
Pic}^0(E)$, then ${\cal M}_{E,r}(\lambda)$ is simply the projective space ${\bf
P}^{r-1}_{\overline {\bf F}_q}$.}
\vskip 0.30cm
\noindent
{\bf 1.3.2. Refined Brill-Noether Locus}
\vskip 0.30cm
Now  let $E$ be an elliptic curve defined over a
finite field ${\bf F}_q$. Then over $\overline E=E\times_{{\bf
F}_q}\overline{{\bf F}_q}$, from 1.3.1, we have the
moduli spaces  ${\cal M}_{{\bar E},r}(\lambda)$
(resp. ${\cal M}_{{\bar E},r}(d)$) of semi-stable bundles of
rank $r$ with determinant $\lambda$ (resp. degree $d$) over $\bar E$. As
algebraic varieties, we may consider ${\bf F}_q$-rational points of these
moduli spaces. Clearly, by  definition, these
rational points of moduli spaces correspond exactly to these classes
of semi-stable bundles which themselves are defined over
${\bf F}_q$. (In the case for ${\cal M}_{{\bar E},r}(\lambda)$,
$\lambda$ is assumed to be rational over ${\bf F}_q$.) Thus for simplicity,
we simply write ${\cal M}_{E,r}(\lambda)$ or ${\cal M}_{E,r}(d)$ for
the corresponding  subsets of ${\bf F}_q$-rational points. For
example, we then simply write
${\rm Pic}^0(E)$ for ${\rm Pic}^0(E)({\bf F}_q)$.
\vskip 0.30cm
Note that if $V$ is semi-stable with strictly positive degree $d$, then
$h^0(E,V)=d$. Hence the standard Brill-Noether locus is either the whole
space or empty. In this way, we are lead to study the case when $d=0$.

For this, recall that for
$\lambda\in {\rm Pic}^0(E)$,
$${\cal M}_{E,r}(\lambda)=\{V:{\rm semi}-{\rm stable}, {\rm
rank}(V)=r,{\rm det}(V)=\lambda\}/\sim_S$$  is identified with
$$\{V=\oplus_{i=1}^rL_i:\otimes_i L_i=\lambda,L_i\in {\rm Pic}^0(E),
i=1,\dots,r\}/\sim_{\rm iso}\simeq {\bf P}^{r-1}$$ where $/\sim_{\rm iso}$
means modulo isomorphisms.

Now introduce the standard Brill-Noether locus
$$W_{E,r}^a(\lambda):=\{[V]\in {\cal M}_{E,r}(\lambda):h^0(E,{\rm
gr}(V))\geq a\}$$  and its \lq stratification' by
$$W_{E,r}^a(\lambda)^0:=\{[V]\in W_{E,r}(\lambda):h^0(E,{\rm
gr}(V))=a\}=W_{E,r}^a(\lambda)\backslash \cup_{b\geq
a+1}W_{E,r}^b(\lambda).$$
One checks easily that $W_{E,r}^a(\lambda)\simeq {\bf P}^{(r-a)-1}$,
$W_{E,r+1}^{a+1}(\lambda)\simeq
W_{E,r}^a(\lambda),$ and $W_{E,r+1}^{a+1}(\lambda)^0\simeq 
W_{E,r}^a(\lambda)^0.$
\vskip 0.30cm
The Brill-Noether theory is based on the consideration of $h^0$. But in
the case for elliptic curves, for arithmetic consideration, such a theory
is not fine enough: not only $h^0$ plays a crucial role, the automorphism
groups are important as well. Based on this, we introduce,
for a fixed
$(k+1)$-tuple non-negative integers
$(a_0;a_1,\dots,a_k)$, the subvariety of $W_{E,r}^{a_0}$ by
setting
$$W_{E,r}^{a_0;a_1,\dots,a_k}(\lambda):=\{[V]\in W_{E,r}^{a_0}(\lambda):
{\rm gr}(V)={\cal
O}_E^{(a_0)}\oplus\oplus_{i=1}^kL_i^{(a_i)},\ \otimes_iL_i^{\otimes
a_i}=\lambda,\ L_i\in {\rm Pic}^0(E),\ i=1,\dots,k\}.$$ Moreover, we
define the associated \lq stratification' by setting
$$W_{E,r}^{a_0;a_1,\dots,a_k}(\lambda)^0:=\{[V]\in
W_{E,r}^{a_0;a_1,\dots,a_k}(\lambda),\ \#\{{\cal
O}_E,L_1,\dots,L_k\}=k+1\}.$$
Easily we see that
$W_{E,r+1}^{a_0+1;a_1,\dots,a_k}(\lambda)\simeq
W_{E,r}^{a_0;a_1,\dots,a_k}(\lambda)$,
$W_{E,r+1}^{a_0+1;a_1,\dots,a_k}(\lambda)^0\simeq
W_{E,r}^{a_0;a_1,\dots,a_k}(\lambda)^0,$ and ${\cal
M}_{E,r}(\lambda)=\cup_{a_0;a_1,\dots,a_k}
W_{E,r}^{a_0;a_1,\dots,a_k}(\lambda)^0,$ where the union is a disjoint
one.

Now recall that for elliptic curves $E$,

\noindent
(1) The quotient space $E^{(n)}/S_n$ is isomorphic to the ${\bf
P}^{n-1}$-bundle over $E$; and

\noindent
(2) The quotient of $E^{(n-1)}/S_n$ is isomorphic to ${\bf
P}^{(n-1)}$. Here we embed $E^{(n-1)}$ as a subspace of $E^{(n)}$ under
the map: $$(x_1,\dots,x_n)\mapsto (x_1,\dots,x_{n-1},x_n)$$ with
$x_n=\lambda-(x_1+x_2+\dots+x_{n-1}).$  Thus, 
$W_{E,r}^{a_0;a_1,\dots,a_k}(\lambda)$ may
be described explicitly as follows.
\vskip 0.30cm
\noindent
{\bf Proposition.} {\it With the same notation as above, regroup
$(a_0;a_1,\dots,a_k)$ as $(a_0;b_1^{(s_1)},\dots,b_l^{(s_l)})$ with the
condition that
$b_1>b_2>\dots>b_l$ and $s_1,s_2,\dots,s_l\in {\bf Z}_{>0}$, then

\noindent
(1) if $b_l=1$, $$W_{E,r}^{a_0;a_1,\dots,a_k}(\lambda)\simeq
\prod_{i=1}^{l-1}{\bf P}_E^{s_i-1}\times {\bf P}^{s_l};$$

\noindent
(2) if $b_l>1$,
$$W_{E,r}^{a_0;a_1,\dots,a_k}(\lambda)\simeq
\prod_{i=1}^{l}{\bf P}_E^{s_i}.$$}
\vskip 0.30cm
\noindent
{\it Remark.} When $\lambda={\cal O}_E$, the refined Brill-Noether loci
$W_{E,r}^{a_0;a_1,\dots,a_k}({\cal O}_E)$  are isomorphic to products of 
(copies of)
projective bundles over $E$ and (copies of) projective spaces, which are 
special subvarieties
in ${\cal M}_{E,r}({\cal O}_E)={\bf P}^{r-1}$. It appears that the 
intersections among these
refined Brill-Noether loci are quite interesting. So define the 
Brill-Noether tautological ring
${\bf BN}_{E,r}({\cal O}_E)$ to be the subring  generated by
all the associated refined Brill-Noether loci (in the corresponding Chow ring).
  What can we say about   it?  As an example, we consider cases when $r=2,3$.
\vskip 0.30cm
\noindent
(1) If $r=2$, then this
ring consists of only two elements: 1-dimensional one $W_{E,2}^{2;0}({\cal
O}_E)=\{[{\cal O}_E\oplus {\cal O}_E]\}$ and the whole ${\bf P}^1$. So 
everything is simple;

\noindent
(2) If $r=3$, then (generators of) this ring contains five elements: 2 of
0-dimensional objects: $W_{E,3}^{3;0}({\cal O}_E)=\{[{\cal O}_E^{(3)}]\}$
and
$W_{E,3}^{1;2}({\cal O}_E)=\{[{\cal O}_E\oplus T_2^{(2)}]:T_2\in E_2\}$
containing 4 elements; 2 of 1-dimensional objects:
$W_{E,3}^{1;1,1}=\{[{\cal O}_E\oplus L\oplus L^{-1}]:L\in {\rm
Pic}^0(E)\}\simeq {\bf P}^1$, a degree 2 projective line contained in
${\bf P}^2={\cal M}_{E,3}({\cal O}_E);$ and
$W_{E,3}^{0;2,1}=\{[L^{(2)}\oplus L^{-2}]:L\in {\rm Pic}^0(E)\}$ a degree
3 curve which is isomorphic to
$E$; and finally the whole space. Moreover, the intersection of
$W_{E,3}^{1;1,1}={\bf P}^1$ and $W_{E,3}^{0;2,1}=E$ are supported on
0-dimensional locus $W_{E,3}^{1;1,1}$, with the multiplicity 3 on the
single point locus $W_{E,3}^{3;0}({\cal O}_E)$ and 1 on the completement
of the points in $W_{E,3}^{1;1,1}$.
\vskip 0.45cm
\noindent
{\bf 1.3.3. Towards A Reciprocity Law: Measuring
Refined Brill-Noether Locus Arithmetically}
\vskip 0.30cm
To measure the Brill-Noether locus, we introduce the following arithmetic
invariant $\alpha_{E,r}(\lambda)$ by setting
$$\alpha_{E,r}(\lambda):=\sum_{V\in [V]\in {\cal
M}_{E,r}(\lambda)}{{q^{h^0(E,V)}}\over {\#{\rm Aut}(V)}}.$$
Also set
$$\alpha^{a_0+1;a_1,\dots,a_k}_{E,r}(\lambda):=\sum_{V\in [V]\in
W^{a_0+1;a_1,\dots,a_k}_{E,r}(\lambda)^0}{{q^{h^0(E,V)}}\over {\#{\rm
Aut}(V)}}.$$

Before going further, we remark that above, we write $V\in [V]$ in the
summation. This is because in each
$S$-equivalence class
$[V]$, there are usually more than one vector bundles $V$. For example,
$[{\cal O}_E^{(4)}]$ consists of ${\cal O}_E^{(4)}$,  ${\cal
O}_E^{(2)}\oplus I_2$, $I_2\oplus I_2$, ${\cal O}_E\oplus I_3$, and $I_4$
\vskip 0.30cm
Due to the importance of automorphism groups, following
Harder-Narasimhan, and Desale-Ramanan, we introduce the following
$\beta$-series invariants $\beta_{E,r}(d)$, $\beta_{E,r}(\lambda)$
and
$\beta_{E,r}^{a_0;a_1,\dots,a_k}(\lambda)$
by setting
$$\beta_{E,r}(d):=\sum_{V\in [V]\in {\cal
M}_{E,r}(d)}{{1}\over {\#{\rm Aut}(V)}},\qquad 
\beta_{E,r}(\lambda):=\sum_{V\in [V]\in {\cal
M}_{E,r}(\lambda)}{{1}\over {\#{\rm Aut}(V)}},$$
and
$$\beta^{a_0;a_1,\dots,a_k}_{E,r}(\lambda):=\sum_{V\in [V]\in
W^{a_0;a_1,\dots,a_k}_{E,r}(\lambda)^0}{{1}\over {\#{\rm
Aut}(V)}}.$$

In particular, we have the following deep
\vskip 0.30cm
\noindent
{\bf Theorem.} ([HN] \& [DR]) {\it For all
$\lambda,\lambda'\in {\rm Pic}^d(E)$,
$$\beta_{E,r}(\lambda)=\beta_{E,r}(\lambda').$$
Moreover,
$$N_1\cdot \beta_{E,r}(\lambda)={{N_1}\over
{q-1}}\cdot\prod_{i=2}^r\zeta_E(i)-\sum_{\Sigma_1^kr_i=r,\Sigma_i
d_i=d,{{d_1}\over {r_1}}>\dots>{{d_k}\over{r_k}}, k\geq
2}\prod_i\beta_{E,r_i}(d_i){1\over{q^{\Sigma_{i<j}(r_jd_i-r_id_j)}}}.$$
Here $N_1$ denotes $\#E(:=\#E({\bf F}_q))$ and $\zeta_E(s)$ denotes the
Artin zeta function for elliptic curve $E/{\bf F}_q$.}
\vskip 0.30cm
\noindent
{\it Remark.} I would like to thank Ueno here, who many years ago draw my 
attentions
to Atiyah and Bott's comments about their Morse theoretical
approach and Harder-Narasimhan's adelic approach towards Poincar\'e 
polynomials of the
associated moduli spaces.
\vskip 0.30cm
Thus, we are lead to introduce the $\gamma$-series invariants
$\gamma_{E,r}(\lambda)$ and
$\gamma_{E,r}^{a_0+1;a_1,\dots,a_k}(\lambda)$ by setting
$$\gamma:=\alpha-\beta.$$ That is to say,
$$\gamma_{E,r}(\lambda):=\sum_{V\in [V]\in {\cal
M}_{E,r}(\lambda)}{{q^{h^0(E,V)}-1}\over {\#{\rm Aut}(V)}},\qquad
\gamma^{a_0;a_1,\dots,a_k}_{E,r}(\lambda):=\sum_{V\in [V]\in
W^{a_0;a_1,\dots,a_k}_{E,r}(\lambda)^0}{{q^{h^0(E,V)}-1}\over {\#{\rm
Aut}(V)}}.$$
\vskip 0.30cm
Clearly, for $\lambda\in {\rm
Pic}^0(E)$,
$$\alpha_{E,r}(\lambda)=\sum_{a_0;a_1,\dots,a_k}\alpha^{a_0;a_1,\dots,a_k}_{ 
E,r}(\lambda),
\qquad
\beta_{E,r}(\lambda)=\sum_{a_0;a_1,\dots,a_k}\beta^{a_0;a_1,\dots,a_k}_{E,r} 
(\lambda)$$
and hence
$$\gamma_{E,r}(\lambda)=\sum_{a_0;a_1,\dots,a_k}\gamma^{a_0;a_1,\dots,a_k}_{ 
E,r}(\lambda).$$

Motivated by the above theorem, we make  the following
\vskip 0.30cm
\noindent
{\bf Conjecture.} {\it For all $\lambda\in {\rm Pic}^0(E)$,
$\alpha_{E,r}(\lambda)=\alpha_{E,r}({\cal
O}_E).$ Hence also
$\gamma_{E,r}(\lambda)=\gamma_{E,r}({\cal O}_E).$}
\vskip 0.30cm
\noindent
{\bf 1.3.4.  Examples In Ranks Two and Three: A Precise Reciprocity Law}
\vskip 0.30cm
Let $E$ be an elliptic curve defined over the finite field ${\bf F}_q$.
\vskip 0.30cm
\noindent
(I) If rank $r$ is two, we need  only to calculate
$\beta_{E,2}(0),\beta_{E,2}(1)$ and $\gamma_{E,r}(0)$.
\vskip 0.30cm
First consider $\beta_{E,2}(0)$. By our discussion on
Brill-Noether locus, it suffices to calculate
$\beta_{E,2}({\cal O}_E)$.
Now $${\cal M}_{E,2}({\cal O}_E)=W_{E,2}^{2;0}({\cal
O}_E)^0\cup W_{E,2}^{0;2}({\cal O}_E)^0\cup W_{E,2}^{0;1,1}({\cal
O}_E)^0.$$ Clearly, $$W_{E,2}^{2;0}({\cal
O}_E)^0=\{[V]: {\rm gr}(V)={\cal O}_E^{(2)}\}$$ consisting of just 1
element;
$$W_{E,2}^{0;2}({\cal O}_E)^0=\{[V]: {\rm gr}(V)=T_2^{(2)},T_2\in
E_2, T_2\not={\cal O}_E\}$$ consisting of 3 elements coming from
non-trivial $T_2\in E_2$, 2-torsion subgroup of $E$; while
$$W_{E,2}^{0;1,1}({\cal O}_E)^0= \{[V]: {\rm gr}(V)=L\oplus L^{-1},L\in
{\rm Pic}^0(E),L\not=L^{-1}\}$$ is simply the complement
of the above 4 points in ${\bf P}^1$.
With this, one checks that
$$\beta_{E,2}(0)=\Big({1\over {(q^2-1)(q^2-q)}}+{1\over
{(q-1)q}}\Big)+3\cdot \Big({1\over {(q^2-1)(q^2-q)}}+{1\over
{(q-1)q}}\Big)+ \big(q+1-(3+1)\big)\cdot {1\over
{(q-1)^2}}
={{q+3}\over {q^2-1}}.$$ And hence
$$\beta_{E,2}(0)=N_1\cdot {{q+3}\over {q^2-1}}.$$

As for $\beta_{E,2}(1)$, it is very simple, since any degree one rank two
semi-stable bundle is stable. Moreover, by the result of Atiyah cited in
1.3.1,  there is exactly one stable
rank two bundle whose determinant is the fixed line bundle. Thus
$$\beta_{E,2}(1)=N_1\cdot{1\over {q-1}}.$$

Finally, we study $\gamma_{E,2}(0)$.
Clearly if $\lambda\not={\cal O}_E$, then $\gamma_{E,2}(\lambda)$ is
supported on $$W^{1;1}_{E,2}(\lambda)=\{[V]:{\rm gr}(V)={\cal
O}_E\oplus\lambda\}$$ consisting only one element with $V={\rm gr}(V)={\cal
O}_E\oplus\lambda.$ So
$$\gamma_{E,2}(\lambda)={{q-1}\over {(q-1)^2}}={1\over {q-1}}.$$

On the other hand,
$\gamma_{E,2}({\cal O}_E)$ is
supported on $$W^{2;0}_{E,2}({\cal O}_E)=\{[V]:{\rm gr}(V)={\cal
O}_E^{(2)}\}$$ consisting only one element too.
Since in the class $[V]$ with ${\rm gr}(V)={\cal
O}_E^{(2)}$, there are two elements, i.e., ${\cal
O}_E^{(2)}$ and $I_2$,
$$\gamma_{E,2}({\cal O}_E)={{q^2-1}\over
{(q^2-1)(q^2-q)}}
+{{q-1}\over{(q-1)q}}={1\over {q-1}}=\beta_{E,1}({\cal O}_E).$$
Thus we have the following
\vskip 0.30cm
\noindent
{\bf Proposition.} {\it With the same notation as above,
$$Z_{E,2, {\bf F}_q}(t)={{N_1}\over
{q-1}}\cdot{{1+(q-1)t+(2q-4)t^2+(q^2-q)t^3+q^2t^4}\over
{(1-t^2)(1-q^2t^2)}}.$$}

This then gives the beautiful absolute Euler product mentioned in 1.2.3.
\vskip 0.30cm
\noindent
(II) When the rank is three,  first, by the fact that
${\rm Aut}({\cal O}_E\oplus I_2)=(q-1)^2q^3$, we have $$\sum_{V,{\rm 
gr}(V){\cal O}_E^{(2)}}{1\over {\#{\rm
Aut}(V)}}={1\over
{(q^2-1)(q^2-q)}}+{1\over {(q-1)q}}$$ and
$$\sum_{W,{\rm gr}(W)={\cal O}_E^{(3)}}{{q^{h^0(E,V)}-1}\over {\#{\rm
Aut}(V)}}={{q^3-1}\over
{(q^3-1)(q^3-q)(q^3-q^2)}}+{{q^2-1}\over {(q-1)^2q^3}}+{{q-1}\over
{(q-1)q^2}}.$$ Consequently,
$$\gamma_{E,3}(0)=N_1\cdot \gamma_{E,3}({\cal
O}_E)=N_1\cdot\beta_{E,2}({\cal O}_E)=N_1\cdot{{q+3}\over {q^2-1}}.$$

So we are left to study $\beta_{E,3}(d)$, $d=0,1,2$. Easily,
$$\beta_{E,3}(1)=\beta_{E,3}(2)=N_1\cdot {1\over q-1}$$ since here
all semi-stable bundles become stable. Thus we are led to consider only
$\beta_{E,3}(0)$, and hence $\beta_{E,3}(\lambda)$ for any
$\lambda\not={\cal O}_E$. (Despite the fact that
$\beta_{E,r}(\lambda)=\beta_{E,r}({\cal O}_E)$ for any $\lambda\in {\rm
Pic}^0(E)$, in practice, the calculation of
$\beta_{E,r}(\lambda)$ with $\lambda\not={\cal O}_E$ is  easier than
that for $\beta_{E,r}({\cal O}_E)$.)

Now $${\cal M}_{E,3}(\lambda)=\Big((W_{E,3}^{2;1}(\lambda)^0)\cup
W_{E,3}^{1;2}(\lambda)^0\cup W_{E,3}^{1;1,1}(\lambda)^0\Big)\cup
W_{E,3}^{0;3}(\lambda)^0
\cup W_{E,3}^{0;2,1}(\lambda)^0\cup W_{E,3}^{0;1,1,1}(\lambda)^0.$$
Moreover, we have

\noindent
(1) $W_{E,3}^{2;1}(\lambda)^0$ consists a single class [V], i.e., the one
with
${\rm gr}(V)={\cal O}_E^2\oplus\lambda$, which contains two vector bundles,
i.e., ${\cal O}_E^2\oplus\lambda$ and $I_2\oplus\lambda$;

\noindent
(2) $W_{E,3}^{2;1}(\lambda)^0\cup
W_{E,3}^{1;2}(\lambda)^0\cup W_{E,3}^{1;1,1}(\lambda)^0\simeq {\bf P}^1$
with
$W^{1;2}(\lambda)^0$ consists of 4 classes $[V]$, i.e., these such that
${\rm gr}(V)={\cal O}_E\oplus \big(\lambda^{1\over 2}\big)^{(2)}$, where
$\lambda^{1\over 2}$ denotes any of the four square roots of $\lambda$.
Clearly then in each class $[V]$, there are also two vector bundles
${\cal O}_E\oplus \big(\lambda^{1\over 2}\big)^{(2)}$ and ${\cal
O}_E\oplus I_2\otimes \lambda^{1\over 2}$;

\noindent
(3) $W_{E,3}^{0;3}(\lambda)^0
\cup W_{E,3}^{0;2,1}(\lambda)^0\cup W_{E,3}^{0;1,1,1}(\lambda)^0={\bf
P}^2\backslash {\bf P}^1$.

\noindent
(3.a) $W_{E,3}^{0;3}(\lambda)^0$ consists of 9 classes $[V]$,
i.e., these $[V]$ with ${\rm gr}(V)=\big(\lambda^{1\over 3}\big)^{(3)}$
where
$\lambda^{1\over 3}$ denotes any of the 9 triple roots of $\lambda$.
Moreover, in each $[V]$, there are three bundles, i.e.,
$\big(\lambda^{1\over 3}\big)^{(3)}$, $\lambda^{1\over
3}\oplus I_2\otimes \lambda^{1\over 3}$ and $I_3\otimes\lambda^{1\over
3}$.

\noindent
(3.b)  $\Big(W_{E,3}^{2;1}(\lambda)^0\cup
W_{E,3}^{1;2}(\lambda)^0\Big)\cup\Big(W_{E,3}^{0;3}(\lambda)^0\cup
W_{E,3}^{0;2,1}(\lambda)^0\Big)$ is isomorphic to $E$. Moreover, each class
$[V]$ in $W^{0;2,1}(\lambda)^0$ consists of two bundles, i.e.,
$L^{(2)}\oplus \lambda\otimes L^{-2}$ and
$I_2\otimes L\oplus \lambda\otimes L^{-2}$ when ${\rm gr}(V)=L^{(2)}\oplus
\lambda\otimes L^{-2}$.

(One checks that in fact the refined Brill-Noether loci ${\bf P}^1$ and
$E$ appeared above are embedded in ${\bf P}^2$ as degree 2 and 3
regular curves. And hence the intersection should be 6: The
intersection points are at
$[V]$ with ${\rm gr}(V)={\cal O}_E^{(2)}\oplus \lambda$ with multiplicity
2, and ${\cal O}_E\oplus \big(\lambda^{1\over 2}\big)^{(2)}$
corresponding to four square roots of $\lambda$ with multiplicity one.
That is to say, the intersection actually are supported on 
$W_{E,3}^{2;1}(\lambda)^0\cup
W_{E,3}^{1;2}(\lambda)^0$. So it would be very interesting in general
to study the intersections of the refined Brill-Noether loci  as well, as 
stated in 1.3.2.)

 From this analysis, we conclude that
$$\eqalign{\beta_{E,3}(\lambda)=&\Big({1\over {(q^2-1)(q^2-q)(q-1)}}
+{1\over {(q-1)q(q-1)}}\Big)\cr
&+4\Big({1\over
{(q^2-1)(q^2-q)(q-1)}}+{1\over
{(q-1)q(q-1)}}\Big)+(q-4)\cdot\Big({1\over
{(q-1)^3}}\Big)\cr &+9\Big({1\over {(q^3-1)(q^3-q)(q^3-q^2)}}+{1\over
{(q-1)^2q^3}}+{1\over {(q-1)q^2}}\Big)\cr
&+\big(N_1-(9+4+1)\big)\cdot\Big({1\over
{(q-1)(q^2-1)(q^2-q)}}+{1\over {(q-1)(q-1)q}}\Big)\cr
&+\big(q^2-(N_1-4-1)\big)\cdot{1\over {(q-1)^3}}.\cr}$$
Therefore, to finally write down the associated non-abelian zeta function, 
it suffices to
use the ugly formula. We leave this to the reader.
\vskip 0.30cm
\noindent
{\bf 1.3.5. Why  Use only Semi-Stable Bundles}
\vskip 0.30cm
At the first glance, it seems  that in
the definition of non-abelian zeta functions we should consider all vector
bundles, just as what happens in the theory of automorphic $L$-functions.
However, we here use an example with $r=2$ to indicate the opposite.

Thus we first introduce a new  zeta
function $\zeta_{E,r}^{\rm all}(s)$  by
$$\zeta_{E,r}^{\rm all}(s):=\sum_{V: {\rm
rank}(V)=2}{{q^{h^0(V)}-1}\over {\#{\rm Aut}(V)}}\cdot q^{-sd(V)}.$$ Then
by our discussion on the non-abelian zeta functions associated to
semi-stable bundle, we only need to consider the contribution of rank 2
bundles which are not semi-stable.

  Assume that $V$ is
not semi-stable of rank 2. Let $L_2$ be the line subbundle of $V$ with
maximal degree, then
$V$ is obtained from the extension of $L_2:=V/L_1$ by $L_1$ $$0\to L_1\to
V\to L_2\to 0.$$ But $V$ is not semi-stable implies that all such
extensions are trivial.  Thus  $V=L_1\oplus L_2$.
For later use, set $d_i$ to be the degree of
$L_i, i=1,2$. Then $d_1+d_2=d$ the degree of $V$, and
$$\#{\rm Aut}(V)=(q-1)^2\cdot q^{h^0(E,L_1\otimes L_2^{-1})}=(q-1)^2\cdot
q^{d_1-d_2}.$$

Next we study the the contribution of degree 0 vector bundles of rank 2
which are not semi-stable. Note that the support of the summation should
have non-vanishing $h^0$. Thus  $V=L_1\oplus L_2$ where
$L_1\in {\rm Pic}^{d_1}(E)$ with $d_1>0$.
So the contributions of these bundles are  given by
$$\eqalign{\zeta_{E,2}^{=0}(s)=&Z_{E,2}^{=0}(t)\cr
=&\sum_{d=1}^\infty\sum_{L_1\in
{\rm Pic}^d(E),L_2\in {\rm
Pic}^{-d}(E)}{{q^{h^0(L_1)}-1}\over{(q-1)^2q^{h^0(L_1\otimes L_2^{\otimes
-1})}}} ={{N_1^2}\over
{(q-1)^2}}\cdot\sum_{d=1}^\infty{{q^d-1}\over{q^{2d}}}\cr
=&{{qN_1^2}\over{(q^2-1)(q-1)^2}}.\cr}$$

Now we consider all degree strictly positive rank 2 vector bundles which
are not semi-stable. From above we see that $V=L_1\oplus L_2$ with
$d_1>d_2$.
Thus for $h^0(E,V)$, there are three  cases:

\noindent
(i) $d_2>0$, clearly then $h^0(E,V)=d$;

\noindent
(ii) $d_2=0$. Here there are two subcases, namely, (a) if $L_2={\cal
O}_E$, then $h^0(E,V)=d_1+1$; (b) If $L_2\not={\cal O}_E$, then
$h^0(E,V)=d_1$;

\noindent
(iii) $d_2<0$. Then $h^0(E,V)=d_1$.

Therefore, all in all the contribution of strictly positive degree rank 2
bundles which are not semi-stable to the zeta function
$\zeta_{E,r}^{\rm all}(s)$ is given by
$$\zeta_{E,2}^{>0}(s)=Z_{E,2}^{>0}(t)=\Big(\sum_{(i)}+\sum_{(ii.a)}+\sum_{(i 
i.b)}
+\sum_{(iii)}\Big){{q^{h^0(V)}-1}\over {\#{\rm Aut}(V)}}t^d$$ where
$\sum_{(*)}$ means the summation is taken for all vector bundles in case
(*). Hence, we have
$$\eqalign{\sum_{(i)}{{q^{h^0(V)}-1}\over {\#{\rm
Aut}(V)}}=&N_1^2\cdot\sum_{d=1}^\infty\sum_{d_1+d_2=d,
d_1>d_2>0}{{q^d-1}\over{(q-1)^2q^{d_1-d_2}}}t^d,\cr
\sum_{(ii.a)}{{q^{h^0(V)}-1}\over
{\#{\rm
Aut}(V)}}=&N_1\cdot\sum_{d=1}^\infty{{q^{d+1}-1}\over{(q-1)^2q^{d}}}
t^{d},\cr
\sum_{(ii.a)}{{q^{h^0(V)}-1}\over {\#{\rm
Aut}(V)}}=&N_1(N_1-1)\cdot\sum_{d=1}^\infty
{{q^{d}-1}\over{(q-1)^2q^{d}}}t^d,\cr
\sum_{(iii)}{{q^{h^0(V)}-1}\over
{\#{\rm Aut}(V)}}=&N_1^2\cdot\sum_{d=1}^\infty\sum_{d_1+d_2=d,
d_1>0>d_2}{{q^{d_1}-1}\over{(q-1)^2q^{d_1-d_2}}}t^d.\cr}$$

By a direct calculation, we find that
$$\eqalign{\sum_{(i)}{{q^{h^0(V)}-1}\over {\#{\rm
Aut}(V)}}=&{{N_1^2t^3}\over {q-1}}\cdot
{{q^2+q+1+q^2t}\over{(1-t^2)(1-q^2t^2)(q-t)}},\cr
\sum_{(ii.a)}{{q^{h^0(V)}-1}\over
{\#{\rm
Aut}(V)}}=&{{N_1t}\over {q-1}}\cdot{{q+1-t}\over{(q-t)(1-t)}},\cr
\sum_{(ii.a)}{{q^{h^0(V)}-1}\over {\#{\rm
Aut}(V)}}=&{{N_1(N_1-1)}\over {q-1}}\cdot {t\over {(q-t)(1-t)}},\cr
\sum_{(iii)}{{q^{h^0(V)}-1}\over
{\#{\rm
Aut}(V)}}=&{{N_1^2t}\over{(q-1)^2(q^2-1)}}\cdot{{q^2+q-1-qt}\over
{(1-t)(q-t)}}.\cr}$$

Finally we consider the contribution of bundles with strictly negative
degree. First we have the following classification according to
$h^0(E,V)$.

\noindent
(i) $d_1>0>d_2$. Then $h^0(E,V)=d_1$;

\noindent
(ii) $d_1=0>d_2$. Here two subcases. (a) $L_1={\cal O}_E$,
then $h^0(E,V)=1$; (b) $L_1\not={\cal O}_E$, then $h^0(V)=0$;

\noindent
(iii) $0>d_1>d_2$. Here $h^0(V)=0$.

Thus note that the support of $h^0(E,V)$ is only on the cases (i) and
(ii.a), we see that similarly as before,
the contribution of strictly negative degree rank 2 bundles
which are not semi-stable to the zeta function is given by
$$\zeta_{E,2}^{<0}(s)=Z_{E,2}^{<0}(t)=\Big(\sum_{(i)}+\sum_{(ii.a)}\Big){{q^ 
{h^0(V)}-1}\over
{\#{\rm Aut}(V)}}t^d,$$
which may be checked to be
$$\eqalign{\zeta_{E,2}^{<0}(s)=Z_{E,2}^{<0}(t)
=&
N_1^2\sum_{d=-1}^{-\infty}\sum_{d_1>0>d_2,
d_1+d_2=d}{{q^{d_1}-1}\over{(q-1)^2q^{d_1-d_2}}}t^d+N_1\cdot\sum_{d=-1}^{-\i 
nfty}
{{q-1}\over {(q-1)^2q^{-d}}}t^d\cr
=&{{N_1^2}\over {(q-1)^2}}\cdot {q\over {(qt-1)(q^2-1)}}
+{{N_1}\over{q-1}}\cdot {1\over {qt-1}}.\cr}$$

I hope now the reader is fully convinced that our definition of
non-abelian zeta function by using moduli space of semi-stable bundles is
much better: Not only our semi-stable zeta functions
have much neat structure, we also have well-behavior geometric and hence
arithmetic spaces ready to use. In a certain sense, we think the picture
of our non-abelian zeta function is quite similar to that the so-called
new forms: Only after  removing these not-semi-stable contributions, we
can see the intrinsic beautiful structures.
\vskip 0.30cm
\noindent
{\it Remark.} Another  way to introduce non-abelian zeta functions using
semi-stable bundles is that when taking the summation, do not take all 
elements in a single
Seshadri equivalence class; instead, choose only one single representative, 
say the one
with maximal automorphism group. We leave the details to the reader.
\vskip 0.5cm
\noindent
{\li Appendix to B.1:  Weierstrass  Groups}
\vskip 0.30cm
Motivated by Kato's construction of Euler systems for elliptic curves in 
terms of elements
in $K_2$ using torsion points, we here introduce what I call Weierstrass 
groups using
Weierstrass divisors for curves.
\vskip 0.3cm
\noindent
{\bf 1. Weierstrass Divisors}
\vskip 0.30cm
\noindent
(1.1) Let $M$ be a compact Riemann surface of
genus $g\geq 2$. Denote its degree $d$ Picard
variety by ${\rm Pic}^d(M)$. Fix a Poincar\'e
line bundle ${\cal P}_d$ on $M\times
{\rm Pic}^d(M)$. (One checks easily that  our
constructions  do not depend on a
particular choice of  the Poincar\'e line bundle.)
Let $\Theta$ be the theta divisor of
${\rm Pic}^{g-1}(M)$, i.e., the image of the
natural map $M^{g-1}\to {\rm Pic}^{g-1}(M)$
defined by $(P_1,\dots,P_{g-1})\mapsto
[{\cal O}_M(P_1+\dots+P_{g-1})]$. Here $[\cdot]$
denotes the class defined by $\cdot$. We will
view the theta divisor as a pair $({\cal
O}_{{\rm Pic}^{g-1}(M)}(\Theta),{\bf 1}_\Theta)$
with ${\bf 1}_\Theta$ the defining section of $\Theta$ via
the structure exact sequence $0\to
{\cal O}_{{\rm Pic}^{g-1}(M)}\to {\cal
O}_{{\rm Pic}^{g-1}(M)}(\Theta)$.

Denote by $p_i$ the $i$-th projection of $M\times
M$ to $M$, $i=1,2$. Then for any degree $d=g-1+n$
line bundle on
$M$, we get a line bundle $p_1^*L(-n\Delta)$
on $M\times M$ which has relative $p_2$-degree
$g-1$. Here, $\Delta$ denotes the diagonal
divisor on $M\times M$. Hence, we get a
classifying map
$\phi_{L}: M\to {\rm Pic}^{g-1}(M)$ which makes the
following diagram commute:
$$\matrix{ M\times M&\to &M\times
{\rm Pic}^{g-1}(M)\cr
p_2\downarrow&&\downarrow \pi\cr
M&\buildrel\phi_{L}\over\to&
{\rm Pic}^{g-1}(M).}$$
One checks that there are  canonical
isomorphisms
$$\lambda_\pi({\cal P}^{g-1})\simeq {\cal
O}_{{\rm Pic}^{g-1}(M)}(-\Theta)$$ and
$$\lambda_{p_2}(p_1^*L(-n\Delta))\simeq
\phi_{L}^*{\cal O}_{{\rm
Pic}^{g-1}(M)}(-\Theta).$$ Here, $\lambda_\pi$
(resp. $\lambda_{p_2}$) denotes the
Grothendieck-Mumford cohomology determinant with
respect to $\pi$ (resp. $p_2$). (See e.g.,
[L].)

Thus, $\phi_{L}^*{\bf 1}_{\Theta}$ gives a
canonical holomorphic section of the dual of the
line bundle
$\lambda_{p_2}(p_1^*L(-n\Delta))$, which in turn
gives an effective divisor $W_L(M)$ on $M$, the
so-called {\it Weierstrass divisor associated to}
$L$.
\vskip 0.30cm
\noindent
{\it Example.} With the same notation as
above, take $L=K_M^{\otimes m}$ with $K_M$ the
canonical  line bundle of $M$ and $m\in
{\bf Z}$. Then we get an effective
divisor
$W_{K_M^{\otimes m}}(M)$ on $M$, which will be
called the $m$-{\it th  Weierstrass divisor}
associated to $M$. For simplicity, denote
$W_{K_M^{\otimes m}}(M)$ (resp.
$\phi_{K_M^{\otimes m}}$)
   by
$W_m(M)$ (resp. $\phi_m$).
\vskip 0.30cm
One checks easily that the degree of $W_m(M)$
is $g(g-1)^2(2m-1)^2$ and we have an
isomorphism ${\cal O}_M(W_m(M))\simeq
K_M^{\otimes g(g-1)(2m-1)^2/2}$.
Thus, in particular,
$$f_{m,n}:={{(\phi_m^*{\bf 1}_\Theta)^{\otimes
(2n-1)^2}}\over {(\phi_{n\ }^*{\bf
1}_\Theta)^{\otimes (2m-1)^2}}}$$ gives a
canonical meromorphic function on $M$ for all
  $m,n\in {\bf Z}$.
\vskip 0.30cm
\noindent
{\it Remark.} We may also assume that $m\in
{1\over 2}{\bf Z}$. Furthermore, this
construction has a relative
version as well, for which we assume that
$f:{\cal X}\to B$ is a semi-stable family of
curves of genus $g\geq 2$. In that case, we get
an effective divisor $({\cal O}_{\cal
X}(W_m(f)),{\bf 1}_{W_m(f)})$ and canonical
isomorphism
$$\eqalign{({\cal O}_{\cal
X}(W_m(f)),&{\bf 1}_{W_m(f)})\cr
\simeq&
({\cal O}_{\cal
X}(W_1(f)),{\bf 1}_{W_1(f)})^{\otimes
(2m-1)^2}\otimes ({\cal O}_{\cal
X}(W_{1\over 2}(f)),{\bf 1}_{W_{1\over
2}(f)})^{\otimes 4m(1-m)}.\cr}$$ The proof may be
given by using Deligne-Riemann-Roch theorem, which
in general, implies that we have the following
canonical isomorphism:
$$({\cal O}_{\cal X}(W_L(f)), {\bf 1}_{W_L(f)})\otimes
f^*\lambda_f(L)\simeq
L^{\otimes n}\otimes K_f^{\otimes n(n-1)/2}.$$
(See e.g. [Bu].) To allow   $m$ be a half
integer, we then should assume
that $f$ has a spin structure. Certainly,
without using spin structure,  a  modified
canonical isomorphism, valid  for integers, can
be given.
\vskip 0.30cm
\noindent
{\bf 2. K-Groups}
\vskip 0.30cm
\noindent
(2.1) Let $M$ be a compact
Riemann surface of genus $g\geq 2$. Then by the
localization theorem, we get the following exact
sequence for
$K$-groups
$$K_2(M)\buildrel\lambda\over\to K_2({\bf C}(M))\buildrel \coprod_{p\in 
M}\partial_p\over\to
\coprod_{p\in M}{\bf C}_p^*.$$

Note that the middle term may also be written as
  $K_2({\bf C}(M\backslash S))$ for any finite
subset $S$ of $M$, we see that naturally by a theorem of Matsumoto,
the Steinberg symbol $\{f_{m,n}, f_{m',n'}\}$ gives
a well-defined element in $K_2({\bf C}(M))$.
Denote the subgroup generated by all $\{f_{m,n},
f_{m',n'}\}$ with $m,n,m',n'\in {\bf Z}_{>0}$ in
$K_2({\bf C}(M))$ as $\Sigma(M)$.
\vskip 0.30cm
\noindent
{\bf Definition.} With the same notation as above,
the {\it first Weierstrass group} $W_I(M)$ of
$M$ is defined to be the $\lambda$-pull-back of
$\Sigma(M)$, i.e., the subgroup
$\lambda^{-1}(\Sigma(M))$ of
$K_2(M)$.
\vskip 0.45cm
\noindent
(2.2) For simplicity, now let $C$ be a regular
projective irreducible curve of genus $g\geq 2$
defined over
${\bf Q}$. Assume that $C$ has a semi-stable
regular module
$X$ over ${\bf Z}$ as well.
Then we have a natural
morphism $K_2(X)\buildrel\phi\over\to
K_2(M)$. Here $M:=C({\bf C})$.
\vskip 0.30cm
\noindent
{\bf Conjecture I.} {\it With the same notation as
above,  $\phi\big(K_2(X)\big)_{\bf Q}
=W_I(M)_{\bf Q}.$}
\vskip 0.45cm
\noindent
{\bf 3. Generalized Jacobians}
\vskip 0.30cm
\noindent
(3.1) Let $C$ be a projective, regular,
irreducible curve. Then for any effective divisor
$D$, one may canonically construct the so-called
generalized Jacobian $J_D(C)$ together with a
rational map
$f_D:C\to J_D(C)$.

More precisely, let $C_D$ be the group of classes
of divisors prime to $D$ modulo these which can
be written as ${\rm div}(f)$. Let $C_D^0$ be the
subgroup of $C_D$ which consists of all
elements of degree zero. For each $p_i$ in the
support of
$D$, the invertible elements modulo those
congruent to 1 (mod $D$) form an algebraic group
$R_{D,p_i}$ of dimension $n_i$, where $n_i$ is the
multiplicity of $p_i$ in $D$. Let $R_D$ be the
product of these $R_{D,p_i}$. One checks easily
that ${\bf G}_m$, the multiplicative group of
constants naturally embeds into $R_D$. It is a
classical result that we then have the short
exact sequence
$$0\to R_D/{\bf G}_m\to C_D^0\to J\to 0$$ where
$J$ denotes the standard Jacobian of $C$. (See
e.g., [S].) Denote $R_D/{\bf G}_m$ simply by
${\bf R}_D$.

Now the map $f_D$ extends naturally to a
bijection from $C_D^0$ to $J_D$. In this way
the commutative algebraic group $J_D$ becomes an
extension as algebraic groups of the standard
Jacobian by the group ${\bf R}_D$.
\vskip 0.30cm
\noindent
{\it Example.} Take the field of constants as
${\bf C}$ and $D=W_m(M)$, the $m$-th Weierstrass
divisor of a compact Riemann surface $M$ of
genus $g\geq 2$. By (1.1), $W_m(M)$ is
effective. So we get the associated generalized
Jacobian $J_{W_m(M)}$. Denote it  by
$WJ_m(M)$ and call it the $m$-th
{\it Weierstras-Jacobian} of $M$. For example, if
$m=0$, then $WJ_0(M)=J(M)$ is the standard
Jacobian of $M$. Moreover, one knows that
  the dimension of $R_{W_m(M),p}$ is at most
$g(g+1)/2$. For later use denote ${\bf
R}_{W_m(M)}$ simply by ${\bf R}_m$.
\vskip 0.45cm
\noindent
(3.2) The above construction works on any base
field as well. We leave the detail to the reader
while point out that if the curve is defined over
a field $F$, then its associated $m$-th
Weierstrass divisor is rational over the same
field as well. (Obviously, this is not true for
the so-called Weierstrass points, which behavior
in a rather random way.) As a consequence, by the
construction of the generalized Jacobian, we see
that the
$m$-th Weierstrass-Jacobians are also defined
over $F$. (See e.g., [S].)
\vskip 0.45cm
\noindent
{\bf 4. Galois Cohomology Groups}
\vskip 0.30cm
\noindent
(4.1) Let $K$ be a perfect field, $\overline K$
be an algebraic closure of $K$ and $G_{\overline
K/K}$ be the Galois group of $\overline K$ over
$K$. Then for any $G_{\overline K/K}$-module $M$,
we have the Galois cohomology groups
$H^0(G_{\overline K/K},M)$ and $H^1(G_{\overline
K/K},M)$ such that if $$0\to M_1\to M_2\to M_3\to
0$$ is an exact sequence of $G_{\overline
K/K}$-modules, then we get  a natural long exact
sequence
$$\eqalign{0\to &H^0(G_{\overline
K/K},M_1)\to
H^0(G_{\overline
K/K},M_2)\to H^0(G_{\overline
K/K},M_3)\cr
&\to H^1(G_{\overline
K/K},M_1)\to H^1(G_{\overline
K/K},M_2)\to H^1(G_{\overline K/K},M_3).\cr
}$$
Moreover, if $G$ is a subgroup
of $G_{\overline K/K}$ of finite index or a
finite subgroup, then
$M$ is naturally a $G$-module. This leads a
restriction map on cohomology
${\rm res}: H^1(G_{\overline K/K},M)\to
H^1(G,M)$.
\vskip 0.30cm
\noindent
(4.2) Now let $C$ be a projective,
regular irreducible curve defined over a number
field $K$. Then for each place $p$ of $K$, fix an
extension of $p$ to $\overline K$, which then
gives
  an embedding
$\overline K\subset \overline K_p$ for the
$p$-adic completion $K_p$ of $K$ and a
decomposition group $G_p\subset G_{\overline
K/K}$.

Now apply the  construction in (3.1) to the short
exact sequence
$$0\to {\bf R}_m\to WJ_m(C)\to J(C)\to 0$$ over
$K$. Then we have the following
long exact sequence
$$\eqalign{0\to& {\bf R}_m(K)\to WJ_m(K)\to
J(K)\cr
&\to H^1(G_{\overline K/K}, {\bf R}_m(K))\to
H^1(G_{\overline K/K}, WJ_m(K))\buildrel
\psi\over\to H^1(G_{\overline K/K},J(K)).\cr}$$

Similarly, for each place $p$ of $K$, we have  the
following exact sequence
$$\eqalign{0\to& {\bf R}_m(K_p)\to WJ_m(K_p)\to
J(K_p)\cr
&\to H^1(G_p, {\bf R}_m(K_p))\to
H^1(G_p, WJ_m(K_p))\buildrel
\psi_p\over\to
H^1(G_p,J(K_p)).\cr}$$

Now the natural inclusion $G_p\subset
G_{\overline K/K}$ and $\overline
K\subset\overline K_p$ give restriction maps on
cohomology, so we  arrive at a natural
morphism
$$\Phi_m:\psi\Big(H^1(G_{\overline K/K},
WJ_m(K))\Big)\to
\prod_{p\in M_K}\psi_p\Big(H^1(G_p,
WJ_m(K_p))\Big).$$
Here $M_K$ denotes the set of all places over $K$.
\vskip 0.30cm
\noindent
{\bf Definition.} With the same notation as above,
{\it the second Weierstrass group} $W_{II}(C)$ of
$C$ is defined to be the subgroup of
$H^1(G_{\overline K/K}, J(C)(K))$
generated by all ${\rm Ker}\,\Phi_m$, the kernel
of $\Phi_m$, i.e.,
$W_{II}(C):=\langle {\rm Ker}\,\Phi_m:m\in {\bf
Z}_{>0}\rangle_{\bf Z}.$
\vskip 0.30cm
\noindent
{\bf Conjecture II.} {\it With the same notation
as above, the second Weierstrass group $W_{II}(C)$
is   finite.}
\vskip 0.45cm
\noindent
{\bf 5.  Deligne-Beilinson Cohomology}
\vskip 0.30cm
\noindent
(5.1) Let $C$ be a projective regular curve of
genus $g$. Let $P$ be a finite set of
$C$. For simplicity, assume that all of them
are defined over {\bf R}. Then
we have the associated Deligne-Belinsion
cohomology group
$H_{\cal D}^1(C\backslash P,{\bf R}(1))$ which
leads to the following short exact sequence:
$$0\to {\bf R}\to H_{\cal D}^1(C\backslash P,{\bf
R}(1))\buildrel{\rm div}\over\to {\bf R}[P]^0\to
0$$ where ${\bf R}[P]^0$ denotes  $({\rm the\
group\ of\ degree\ zero\ divisors\ with\ support\
on\ } P)_{\bf R}$.

The standard cup product  on
Deligne-Beilinson cohomology leads to a
well-defined map:
$$\cup:H_{\cal D}^1(C\backslash P,{\bf
R}(1))\times H_{\cal D}^1(C\backslash P,{\bf
R}(1))\to H_{\cal D}^2(C\backslash P,{\bf
R}(2)).$$

Furthermore, by Hodge theory, there is a canonical
short exact sequence
$$0\to H^1(C\backslash P,{\bf R}(1))\cap
F^1(C\backslash P)\to H_{\cal D}^2(C\backslash
P,{\bf R}(2))\buildrel {p_{\cal D}}\over \to
  H_{\cal D}^2(C,{\bf R}(2))\to 0$$ where
$F^1$ denotes the $F^1$-term of the Hodge
filtration on $H^1(C\backslash P,{\bf C})$.

All this then leads to a well-defined morphism
$$[\cdot,\cdot]_{\cal D}:\wedge^2 {\bf R}[P]^0\to
H_{\cal D}^2(C,{\bf R}(2))=H^1(C,{\bf R}(1))$$
which make the associated diagram coming from the
above two short exact sequences commute. (See e.g.
[Bei].)

\noindent
(5.2) Now applying the above construction with
$P$ being the union of the  supports of $W_1$,
$W_m$ and $W_n$ for $m,n>0$. Thus for fixed
$m,\,n$, in ${\bf R}[P]^0$,  we get two
elements
${\rm div}(f_{1,m})$ and ${\rm div}(f_{1,n})$.
This then gives $[{\rm div}(f_{1,m}), {\rm
div}(f_{1,n})]_{\cal D}\in H^1(X,{\bf R}(1)).$ Thus, by a
simple argument using the Stokes
formula, we obtain the following

\noindent
{\bf Lemma.} {\it For any holomorphic
differential 1-form $\omega$ on $C$, we have
$$\eqalign{~&\langle [{\rm div}(f_{1,m}), {\rm
div}(f_{1,n})]_{\cal D},\omega\rangle
:=-{1\over{2\pi {\sqrt -1}}}\int [{\rm
div}(f_{1,m}), {\rm div}(f_{1,n})]_{\cal
D}\wedge\bar\omega\cr =&-{1\over{2\pi {\sqrt -1}}}
\int g({\rm div}(f_{1,m}),z)d g({\rm
div}(f_{1,n}),z)\wedge\bar\omega,\cr}$$ Here
$g(D,z)$ denotes the Green's function of $D$ with
respect to any fixed normalized (possibly
singular) volume form of quasi-hyperbolic type. (See e.g., [We])}

\noindent
(5.3) With exactly the same notation as in (4.2), then in
$H^1(X,{\bf R}(1))$ we get a collection of
elements $[{\rm div}(f_{1,m}), {\rm
div}(f_{1,n})]_{\cal D}$ for $m,n\in {\bf
Z}_{>0}$.

\noindent
{\bf Definition.}  With the same notation as
above, assume that $C$ is defined over {\bf Z}.
Define the {\it --first quasi-Weierstrass group}
$W_{-I}'(C)$ of $C$ to be the subgroup of
$H^1(X,{\bf R}(1))$ generated by $[{\rm div}(f_{1,m}), {\rm
div}(f_{1,n})]_{\cal D}$ for all $m,n\in {\bf
Z}_{>0}$ and call  $W_{-I}'(C)_{\bf Q}$
the {\it --first Weierstrass group} $W_{-I}(C)$ of
$C$. That is to say,
$W_{-I}(C):=\langle [{\rm div}(f_{1,m}), {\rm
div}(f_{1,n})]_{\cal D}: m,n\in {\bf
Z}_{>0}\rangle_{\bf Q}.$
\vskip 0.30cm
\noindent
{\bf Conjecture III.} {\it With the same notation
as above, $W_{-1}(C)_{\bf R}$ is the full space,
i.e. equals to
$H^1(X,{\bf R}(1))$.}
\vskip 0.30cm
\noindent
That is to say, Weierstrass divisors should
give a new rational structure for $H^1(X,{\bf
R}(1))$, and hence the
corresponding regulator should give the leading
coefficient of the $L$-function of $C$ at $s=0$,
up to rationals.
\vskip 0.5cm
\noindent
{\li  B.2. New Non-Abelian Zeta Functions for Number Fields}
\vskip 0.30cm
\noindent
{\bf B.2.1. Iwasawa's ICM Talk on Dedekind Zeta Functions}
\vskip 0.30cm
As  for function fields, here we start with a discussion
on  abelian zeta functions for number fields, i.e.,
Dedekind zeta functions. However,
we will not adapt the classical approach, rather we would like to recall
Iwasawa's interpretation. (Based on the fact
that Iwasawa's original choice of certain auxiliary functions do not 
naturally lead to any
meaningful cohomology, some subtle changes are made.)

Let $F$ be a number field. Denote by $S$ the collection of all
(unequivalent) normalized places of $F$. Set $S_\infty$ to be the
collection of all Archimedean places of $F$ and $S_{\rm fin}:=S\backslash
S_\infty$.

Denote by {\bf I} the idele group of $F$, $N:{\bf I}\to {\bf R}_{\geq 0}$
and ${\rm deg}:{\bf I}\to {\bf R}$ the norm map and the degree map on
ideles respectively. Also introduce the following subgroups of
{\bf I}:

$$\eqalign{{\bf I}^0:=&\{a=(a_v)\in {\bf I}:{\rm deg}(a)=0\},\cr
F^*:=&\{a=(a_v)\in {\bf
I}^0:a_v=\alpha\in F\backslash \{0\},\forall v\in S\},\cr
U:=&\{a=(a_v)\in {\bf I}^0:|a_v|_v=1\forall v\in S\},\cr
{\bf I}_{\rm fin}:=&\{a=(a_v)\in {\bf I}:a_v=1\forall v\in
S_\infty\},\cr
{\bf I}_\infty:=&\{a=(a_v)\in {\bf I}:a_v=1\forall v\in
S_{\rm fin}\}.\cr}$$
  Set $U_{\rm fin}:=U\cap {\bf I}_{\rm fin}$. Then, with
respect to the natural topology on {\bf I}, we have

\noindent
(1) $F\hookrightarrow {\bf I}$ is discrete and ${\bf I}^0/F^*$ is compact.
Write ${\rm Pic}(F)={\bf I}/F^*$;

\noindent
(2) $U\hookrightarrow {\bf I}$ is compact;

\noindent
(3) $U_{\rm fin}\hookrightarrow {\bf I}_{\rm fin}$ is both open and compact.
Moreover, the morphism $I:[a=(a_v)]\mapsto I(a):=\prod_{v\in S_{\rm
fin}}P_v^{{\rm ord}_v(a_v)}$ induces an isomorphism between ${\bf I}_{\rm
fin}/U_{\rm fin}$ and the ideal group of $F$, where $P_v$ denotes the
maximal ideal of the ring of integers ${\cal O}_F$ corresponding to the
place $v$; and, $N(a)=\prod_{v\in S}|a_v|_v^{N_v:=[F_v:{\bf
Q}_p]}=N(I(a))^{-1}$ with $N(I(a))$ the norm of the ideal $I(a)$;

\noindent
(4) ${\bf I}={\bf I}_{\rm fin}\times {\bf I}_\infty$. Hence we may
write an idele $a$ as
$a=a_{\rm fin}\cdot a_\infty$ with $a_{\rm fin}
\in {\bf I}_{\rm fin}$ and $a_\infty\in {\bf I}_\infty$ respectively.
In particular, if $d\mu(a)$ denotes the normalized Haar measure on {\bf
I} as say in Weil's Basic Number Theory, we have $$d\mu(a)=d\mu(a_{\rm
fin})\cdot d\mu(a_\infty)$$ corresponding to the decomposition
${\bf I}={\bf I}_{\rm fin}\times {\bf I}_\infty$.

Set $$e(a_{\rm fin}):=\cases{1,&if $I(a_{\rm fin})\subset {\cal O}_F$,\cr
0,&otherwise,\cr}$$  $$e(a_\infty)=\exp(-\pi\sum_{v:{\bf
R}}a_v^2-2\pi\sum_{v:{\bf C}}|a_v|^2)$$ and
$$e(a_{\rm fin}\cdot a_\infty)=e(a_{\rm fin})\cdot e(a_\infty),$$
for  $a_{\rm fin}\in {\bf I}_{\rm fin}$ and $a_\infty\in {\bf I}_\infty$,
Denote by $\Delta_F$ (the absolute values of) the
discriminant of $F$, $r_1$ and $r_2$ the number of real and complex
places in $S_\infty$ as usual.
\vskip 0.30cm
Now we are ready to write down Iwasawa's
interpretation of the Dedekind zeta function for $F$ in the form suitable
for our later study. This goes as follows.

For $s\in {\bf C}, {\rm Re}(s)>1$,

$$\eqalign{\xi_F(s):=&\Delta_F^{s\over 2}(2\pi^{-{s\over 2}}\Gamma({s\over
2}))^{r_1}((2\pi)^{-s}\Gamma(s))^{r_2}\sum_{0\not=I\subset {\cal
O}_F}N(I)^{-s}\cr
=&\Delta_F^{s\over 2}\cdot \sum_{0\not=I\subset {\cal
O}_F}N(I)^{-s}\int_{t_v\in F_v,v\in
S_\infty}\prod_v|t_v|_v^s\exp(-\pi\sum_{v:{\bf R}}t_v^2-2 \pi\sum_{v:{\bf
C}}|t_v|^2)\prod_vd^*t_v\cr
=&\Delta_F^{s\over 2}\cdot \sum_{0\not=I\subset {\cal
O}_F}N(I)^{-s}\int_{{\bf
I}_\infty}N(a_\infty)^se(a_\infty)d\mu(a_\infty)\cr
  =&\Delta_F^{s\over
2}\cdot {1\over{{\rm vol}(U_{\rm fin})}}\cdot
\Big(\int_{{\bf
I}_{\rm fin}}N(a_{\rm fin})^se(a_{\rm fin})d\mu(a_{\rm fin})
\cdot
\int_{{\bf
I}_\infty}N(a_\infty)^se(a_\infty)d\mu(a_\infty)\Big)\cr
=&{1\over{{\rm vol}(U_{\rm fin})}}\cdot\Delta_F^{s\over 2}\cdot
\int_{{\bf
I}_{\rm fin}\times {\bf
I}_\infty}\Big(N(a_{\rm fin})N(a_\infty)\Big)^s\Big(e(a_{\rm
fin})e(a_\infty)\Big)\Big(d\mu(a_{\rm fin})d\mu(a_\infty)\Big)\cr
=&{1\over{{\rm vol}(U_{\rm fin})}}\cdot\Delta_F^{s\over 2}\cdot
\int_{\bf I}N(a)^se(a)d\mu(a).\cr}$$

Now denote by $d\mu([a])$ the induced Haar measure on the Picard group
${\rm Pic}(F):={\bf I}/F^*$. Note that
$$e(a_{\rm fin}a_\infty):=\cases{e(a_\infty),&if $I(a_{\rm
fin})\subset {\cal O}_F$,\cr 0,&otherwise,\cr}$$
and that by (1) above, $F^*\hookrightarrow {\bf I}$ is discrete, (hence
taking integration over $F^*$ means taking summation), we get

$$\eqalign{\xi_F(s)
=&{1\over{{\rm vol}(U_{\rm fin})}}\cdot\Delta_F^{s\over 2}\cdot
\int_{{\bf I}/F^*}\Big(\int_{F^*}
N(\alpha a)^se(\alpha a)d\alpha\Big)d\mu([a])\cr
=&{1\over{{\rm vol}(U_{\rm fin})}}\cdot\Delta_F^{s\over 2}\cdot
\int_{{\bf I}/F^*}\Big(\sum_{\alpha\in F^*}
e(\alpha a)\Big)\cdot N([a])^sd\mu([a])\cr
&\qquad({\rm by\ the\ product\ formula})\cr
=&{1\over{{\rm vol}(U_{\rm fin})}}\cdot\Delta_F^{s\over 2}\cdot
\int_{{\bf I}/F^*} N([a])^sd\mu([a])
\cdot \Big(\sum_{\alpha\in F^*}
e(\alpha a_{\rm fin})e(\alpha a_\infty)\Big)\cr
=&{1\over{{\rm vol}(U_{\rm fin})}}\cdot\Delta_F^{s\over 2}\cdot
\int_{{\bf I}/F^*} N([a])^sd\mu([a])
\cdot \Big(\sum_{\alpha\in F^*,I(\alpha a_{\rm fin})\subset{\cal O}_F}
e(\alpha a_\infty)\Big)\cr
=&{1\over{{\rm vol}(U_{\rm fin})}}\cdot\Delta_F^{s\over 2}\cdot
\int_{{\bf I}/F^*} N([a])^sd\mu([a])
\cdot \Big(\sum_{\alpha\in F^*,\alpha a_v\subset{\cal O}_v,\forall v\in
S_{\rm fin}} e(\alpha a_\infty)\Big).\cr}$$
Now  for an idele  $L=(a_v)$, define the $0$-th algebraic cohomology
group of the idele
$L^{-1}$ by $$H^0(F,L^{-1}):=\{\alpha\in F^*,\alpha
a_v\subset{\cal O}_v,\forall v\in S_{\rm fin}\}=I(L^{-1});$$ moreover
for the associated idele class, introduce its associated geometric 0-th 
geo-ari cohomology
via
$$h^0(F,L^{-1}):=\log\Big(\sum_{\alpha\in
H^0(F,L^{-1})\backslash\{0\}}\exp
\Big(-\pi\sum_{v:{\bf R}}|g_v\alpha|^2-2\pi\sum_{v:{\bf 
C}}|g_v\alpha|^2\Big).$$
With this, then easily, we have $$\eqalign{\xi_F(s) =&{1\over{{\rm vol}(U_{\rm
fin})}}\cdot\Delta_F^{s\over 2}\cdot
\int_{{\rm Pic}(F)}{1\over{w_F}}\sum_{\alpha\in
H^0(F,L^{-1})\backslash\{0\}}\exp
\Big(-\pi\sum_{v:{\bf R}}|g_v\alpha|^2-2\pi\sum_{v:{\bf C}}|g_v\alpha|^2
\Big) N(L)^sd\mu(L)\cr
=&{1\over{{\rm vol}(U_{\rm fin})}}\cdot\Delta_F^{s\over 2}\cdot
\int_{{\rm
Pic}(F)}{1\over{w_F}}\Big(e^{h^0(F,L^{-1})}-1\Big)N(L)^sd\mu(L)\cr
=&{1\over{{\rm vol}(U_{\rm fin})}}\cdot {1\over{w_F}}\cdot\Delta_F^{s\over
2}\cdot
\int_{{\rm Pic}(F)}\Big(e^{h^0(F,L)}-1\Big)N(L)^{-s}d\mu(L)\cr
=&{1\over{{\rm vol}(U_{\rm fin})}}\cdot\Delta_F^{s\over 2}\cdot
\int_{{\rm Pic}(F)}{{e^{h^0(F,L)}-1}\over{{\rm
Aut}(L)}}\cdot\Big(e^{{\rm deg}(L)}\Big)^{-s}d\mu(L).\cr}$$ Here $w_F$
denotes the number of units of
$F$, and ${\rm Aut}(L)(=w_F)$ denotes the number of automouphisms of $L$.
\vskip 0.30cm
It is very clear that, formally, this version of (the completed) Dedekind 
zeta functions for
number fields stands exactly the same as our interpretation of Artin zeta 
functions for
function fields. So to introduce an non-abelian zeta function for number 
fields, key points
are the follows:

\noindent
(1) A suitable stability in terms of intersection for bundles over number 
fields should be
introduced;

\noindent
(2) Stable bundles over number fields should form  moduli spaces, over which
there exist  natural measures; and

\noindent
(3) There should be a geo-ari cohomology such that duality and Riemann-Roch
type results hold.
\vskip 0.30cm
\noindent
{\bf B.2.2. Intersection Stability}
\vskip 0.30cm
\noindent
{\bf B.2.2.1.  Classification of Unimodular Lattices: A Global Approach}
\vskip 0.30cm
Even though it  has not yet been very popular, intersection stability in 
arithmetic
is indeed quite fundamental,  as what we are going to see.

Recall that a full rank lattice
$\Lambda\subset {\bf R}^r$ is said to be {\it integral} if for any
$x\in \Lambda$, $(x,x)$ is an integer,  and that an integral lattice 
$\Lambda$ is called {\it
unimodular}, if the volume of its fundamental domain is  one. It is a 
classical yet still very
challenging problem to classify all unimodular lattices.

Roughly speaking,  classifications of unimodular lattices consist of two 
different aspects, i.e.,
the local and the global one. For the local study, we are mainly interested in
enumerating all unimodular lattices, which in recent years proves to be 
very fruitful.
However for the global study, besides the pioneer works done by Minkowski 
and Siegel,
such as the mass formula and asymptotic upper bounds for the numbers of 
unimodular
lattices in terms of  volumes of Siegel domains, less progress has been 
recorded.

One of the main difficulties in the global study is that Siegel domains
are hard to be understood. Thus, it seems to be very essential to find a
natural method to divide these domains into certain well-behavior blocks. 
Motivated by
what happens for bundles over function fields, in particular, the so-called 
Mumford stability
and the associated Harder-Narasimhan filtration, we introduce the following
\vskip 0.30cm
\noindent
{\bf Definition.} A lattice $\Lambda$ is called {\it stable} (resp. {\it
semi-stable}) if for any proper sublattice $\Lambda'$,
$${\rm Vol}(\Lambda')^{{\rm rank}(\Lambda)}\qquad>\qquad({\rm
resp.}\geq)\qquad {\rm Vol}(\Lambda)^{{\rm rank}(\Lambda')}.$$

Standard properties about Harder-Narasimhan filtrations
and Jordan-H\"older filtrations hold here as well. That is to say, we
have the following:
\vskip 0.30cm
\noindent
{\bf Propposition.} {\it Let $\Lambda$ be a lattice. Then

\noindent
(1) There exists a unique filtration of proper sublattices,
$$0=\Lambda_0\subset\Lambda_1\subset\dots\dots\subset
\Lambda_s=\Lambda$$ such that $\Lambda_i/\Lambda_{i-1}$ is semi-stable and
$${\rm Vol}(\Lambda_{i+1}/\Lambda_{i})^{{\rm
rank}(\Lambda_i/\Lambda_{i-1})}>
{\rm Vol}(\Lambda_{i}/\Lambda_{i-1})^{{\rm
rank}(\Lambda_{i+1}/\Lambda_{i})}.$$

\noindent
(2) If moreover $\Lambda$ is semi-stable, then there exists a
filtration of proper sublattices,
$$0=\Lambda^{t+1}\subset\Lambda^{t}\subset\dots\dots\subset
\Lambda^0=\Lambda$$ such that $\Lambda^j/\Lambda^{j+1}$ is stable and
$${\rm Vol}(\Lambda^{j}/\Lambda^{j+1})^{{\rm
rank}(\Lambda^{j-1}/\Lambda^{j})}=
{\rm Vol}(\Lambda^{j-1}/\Lambda^{j})^{{\rm
rank}(\Lambda^{j}/\Lambda^{j+1})}.$$ Furthermore, the graded lattice
${\rm Gr}(\Lambda):=\oplus \Lambda^{j-1}/\Lambda^{j}$, the so-called 
Jordan-H\"older graded
lattice of $\Lambda$, is uniquely determined by $\Lambda$.}
\vskip 0.30cm
Thus, in particular,  for unimodular lattices, we have the
following
\vskip 0.30cm
\noindent
{\bf Corollary.} {\it Unimodular lattices are semi-stable. Moreover, a 
unimodular lattice
is stable if and only if it contains no proper  unimodular sublattice.}
\vskip 0.30cm
In this sense, to classify all unimodular lattices, it suffices to classify 
all stable unimodular
lattices. This then leads to the
following consideration.

Denote by ${\cal M}_{{\bf Q},r}(1)$ the collection, or better, the moduli 
space, of all rank $r$
semi-stable lattices of  volume one.
Then one checks that ${\cal M}_{{\bf Q},r}(1)$ admits
a natural metric and is indeed compact.
\vskip 0.30cm
\noindent
{\bf Example.} With the help of the  reduction theory
from geometry of numbers, $${\cal M}_{{\bf 
Q},2}(1)\simeq\Big\{\Big(\matrix{a&0\cr
b&{1\over a}\cr}\Big):1\leq a\leq\sqrt{2\over{\sqrt 
3}},\sqrt{a^2-a^{-2}}\leq b\leq
a-\sqrt{a^2-a^{-2}}\Big\}\cdot {\rm SO}(2).$$ Moreover, 
$$\Big\{\Big(\matrix{a&0\cr b&{1\over
a}\cr}\Big):1\leq a\leq\sqrt{2\over{\sqrt 3}},\sqrt{a^2-a^{-2}}\leq b\leq
a-\sqrt{a^2-a^{-2}}\Big\}$$ may be viewed as a closed bounded domain in the 
upper half
plane. As a direct consequence, ${\cal M}_{{\bf Q},2}(1)$ admits a natural 
metric as well,
induced from the Poincar\'e metric.
\vskip 0.30cm
To go back to unimodular lattices, we may now view them
naturally as certain special points in our geometric moduli spaces. (Recall 
that
unimodular lattices are integral.) In this way,
the problem of classifying unimodular lattices looks very much
similar to that of finding rational points in algebraic varieties. Thus, 
along with the line of
Minkowski's geometry of numbers, the first thing we have to do is to evaluate
volumes of these moduli spaces with respect to the associated natural 
metrics. It is for this
purpose that
we introduce our non-abelian zeta functions for number fields: Theory of
Dedekind zeta functions tells us that,  regulators of number fields, or better,
volumes of the lattices generated by  fundamental units, may be read from 
the residues of
Dedekind zeta functions at the simple poles $s=1$.
\vskip 0.30cm
\noindent
{\bf B.2.2.2. Semi-Stable  Bundles over Number Fields}
\vskip 0.30cm
Over general number fields, we may also introduce the 
intersection-stability for
(parabolic $G$) bundles according to the following observations:

\noindent
(1) there exists a well-developed Arakelov theory, from which in particular 
we have the
concept like hermitian vector sheaves, rank and degree;

\noindent
(2) over each local fields, the relation between 1-PS and weighted 
filtration is
well-understood for reductive groups.
\vskip 0.30cm
For example, in terms of Arakelov theory, let $({\cal E},\rho)$ be a 
hermitian vector
sheaf over a number field $F$, or, better over, the spectrum of the ring of 
integers.
Then the rank and the Arakelov degree makes sense. Introduce the
(Arakelov) $\mu$-invariant by
$\mu({\cal E},\rho):={{{\rm deg}({\cal E},\rho)}\over {{\rm rank}({\cal E})}}$.
Then by definition, $({\cal E},\rho)$ is called semi-stable (resp. stable) 
if for any
metrized vector subsheaf $({\cal E}_1,\rho_1=\rho|_{{\cal E}_{1,\infty}})$,
$$ \mu({\cal E}_1,\rho_1)\ \ \leq\quad({\rm resp.}\ <)\ \ \mu({\cal F},\rho).$$
 
Moreover, just as in 2.2.1 over {\bf Q}, in general, standard properties about
Harder-Narasimhan filtration and Jordan-H\"older filtrations holds here as 
well, based on
the fact that, as a subset of {\bf R},  $$\{\mu({\cal E}_1,\rho_1): {\cal 
E}_1\ {\rm is\ a\ vector\
subsheaf\ of}\ {\cal E},\ {\rm and}\ \rho_1=\rho|_{{\cal E}_{1,\infty}}\}$$ 
is  discrete and
bounded from above.

We will leave the corresponding generalization to parabolic $G$-bundles to 
the reader.
Instead, we want to introduce an adelic version of the stability with the 
aim to construct
the corresponding moduli spaces and the associated Tamagawa measures.
\vskip 0.30cm
\noindent
{\bf B.2.2.3. Adelic Moduli and Its Associated Tamagawa Measure}
\vskip 0.30cm
As above, let $F$ be a number field, i.e.,   a finite extension of {\bf Q}. 
Denote
by $S=S_F$ the collection of all (unequivalent)
normalized places of $F$. Set $S_\infty$  be the
collection of all Archimedean places in
$S$, and $S_{\rm fin}:=S\backslash S_\infty$.

For each $v\in S$, denote by $F_v$ the $v$-completion of $F$.
If $v\in S_\infty$,  $F_v$ is {\bf R} or {\bf C}.
We will then call  $v$ (resp. denote $v$) a real or a complex place (resp.
{\bf R} or {\bf C}) accordingly. If
$v\in S_{\rm fin}$, denote by ${\cal O}_v$ the ring of $v$-adic integers of
$F_v$, and ${\cal M}_v$ its maximal ideal. Fix also a generator $\pi_v$
of ${\cal M}_v$.

Fix a positive integer $r$. For all $v\in S$, let
$G_v:={\rm GL}_r(F_v)$. If
$g_v=(g_{v,ij})\in G_v$, denote its inverse by $g_v^{-1}=(g_v^{ij})$.
With this, for
$v\in S_{\rm fin}$, introduce a subgroup $U_v$ of $G_v$  by setting
$$U_v:=\{g_v\in G_v:g_v=(g_{v,ij})\ {\rm with}\ g_{v,ij}\in {\cal O}_v\ {\rm
and}\ g_v^{ij}\in {\cal O}_v\}.$$ Now, following Weil,  define the
associated adelic group $G({\bf A}_F)$ via
$$G({\bf A}_F):={\rm GL}_r({\rm A}_F):=\Big\{g=(g_v)_{v\in S}: g_v\in G_v\ {\rm
s.t.\ for\ almost\ all\ but\ finitely\ many}\ v\in S_{(\rm fin)}, g_v\in
U_v\Big\}.$$

Note that  $G(F):={\rm GL}_r(F)$ may be naturally
embedded into ${\rm GL}_r({\bf A}_F)$ via the diagonal map $\alpha\mapsto
(\alpha,\dots,\alpha,\dots)$. One checks that with respect to the
natural topology on $G({\bf A}_F)$,
$G_r(F)$ is a discrete subgroup of $G_r({\bf A}_F)$. So
we may form the quotient group
$G_r(F)\backslash G_r({\bf A}_F)$. By
definition, a rank $r$ {\it pre vector bundle} on a number field $F$ is an 
element
$[g]\in {\rm GL}_r(F)\backslash {\rm GL}_r({\bf A}_F)$. Also for our own
convenience, we call an element $g\in {\rm GL}_r({\bf A}_F)$ a rank $r$
{\it matrix divisor} on $F$. Two rank $r$ matrix divisor $g=(g_v)$ and
$g'=(g_v')$ are said to be ({\it rationally}) {\it equivalent}
if there is an element
$\alpha\in G(F)$ such that $g_v=\alpha\cdot g_v'$ for all $v\in S$.

For example, if $r=1$, then $G({\bf A}_F)$ is simply the
collection of invertible elements in ${\bf A}_F$, the ring of adeles of $F$,
i.e., ${\rm GL}_1({\bf A}_F)={\bf I}_F$, the group of ideles of
$F$. Hence a pre line bundle on  $F$ is
indeed an element in $F^*\backslash {\bf I}_F$. Moreover, we may view 
elements in
${\bf I}_F$, the rank 1 matrix divisors on $F$ in our language, as
divisors on $F$.

Associated to a rank $r$ matrix divisor $g=(g_v)$ is naturally a hermitian 
vector sheaf
$({\cal E}(g),\rho(g))$ over $F$. Indeed, we may set ${\cal E}(g)$ to be 
$\{\alpha\in
F^r:g_v^{-1}\cdot \alpha\in {\cal O}_v^r,\forall v\in S_{\rm fin}\}$ which 
is a rank $r$ vector
sheaf on ${\rm Spec}({\cal O}_F)$ and may be naturally embedded into $({\bf 
R}^{r_1}\times
{\bf C}^{r_2})^r$, where as usual, $r_1$ and $r_2$ denote the real and 
complex embeddings
of $F$ respectively. View $(g_\sigma)_{\sigma\in S_\infty}$ as an ismorphism of
$({\bf R}^{r_1}\times {\bf C}^{r_2})^r$, and then define $\rho(g)$ as the 
natural metric on
  ${\cal E}(g)_\infty$  induced from the Euclidean metric on this latest 
$({\bf R}^{r_1}\times
{\bf C}^{r_2})^r$. Clearly if $g$ and $g'$ are rational equivalent, then 
the associated
hermitian vector sheaves $({\cal E}(g),\rho(g))$ and $({\cal 
E}(g'),\rho(g'))$ are isometric to
each other. Hence it makes sense to talk the associated hermitian vector 
sheaf for a pre
vector bundle. If the Arakelov  degree of $({\cal E}(g),\rho(g))$ is $d$, 
the $g$ and $[g]$
are said to be of degree $d$.

By definition, a pre vector bundle $[g]$ of a number field $F$ is 
semi-stable (resp. stable),
if its associated hermitian vector sheaf $({\cal E}(g),\rho(g))$ is 
semi-stable (resp. stable);
and  the adelic moduli space ${\cal M}_{{\bf A}_F,r}(d)$ of semi-stable pre 
vector bundles of
rank $r$ and degree $d$ is a collection of all  semi-stable pre vector 
bundles of
rank $r$ and degree $d$.
 
Clearly, semi-stability is a closed condition. Moreover, for a fixed 
degree, semi-stability
is also a bounded condition. Thus, we conclude that the moduli space ${\cal 
M}_{{\bf
A}_F,r}(d)$ is indeed compact.

The advantage of using the adelic moduli space ${\cal M}_{{\bf A}_F,r}(d)$ 
is that then we
may obtain a natural measure, the Tamagawa one. In fact, as a compact 
subset in ${\rm
Gl}_r(F)\backslash {\rm GL}_r({\bf A}_F)$, ${\cal M}_{{\bf A}_F,r}(d)$ 
inherits a natural
measure from the Tamagawa measure on the total space. For simplicity, we 
simply call this
induced measure on ${\cal M}_{{\bf A}_F,r}(d)$ as the (associated) Tamagawa 
measure.

On the other hand, using Seshadri type equivalence, we may also introduce 
the so-called
moduli space ${\cal M}_{F,r}(d)$ of semi-stable vector sheaves of rank $r$ 
and degree $d$ on
$F$. By the uniqueness of the Jordan-H\"older graded hermitian vector 
sheaves, we obtain
a well-defined continuous map $$\pi_F:{\cal M}_{{\bf A}_F,r}(d)\to {\cal 
M}_{F,r}(d).$$
Now by the so-called finiteness result of Borel on adelic groups over 
number fields, we see
that the fiber of $\Pi_F$ is indeed compact. Thus, naturally,  we get a 
natural finite  measure
on ${\cal M}_{F,r}(d)$ as well. (In fact, for vector bundles, only a weak 
version of Borel's
result is needed here. But if we want to study parabolic $G$-bundles, then 
a full version of
Borel's result has to be used.)
\vskip 0.30cm
\noindent
{\bf Remark.} As suggested in Part (A), we should develop a GIT in terms of 
Arakelov
theory. If so, then we may have a new construction of the moduli space of 
semi-stable
vector bundles, for which, the above map $\Pi_F$ may be viewed as an analog 
of the
moment map.
\vskip 0.30cm
Clearly the relation between the Tamagawa volume of ${\cal M}_{{\bf A}_F,r}(d)$
and the volume of ${\cal M}_{F,r}(d)$ (with respect to the measure induced 
from that of
Tamagawa measure via $\pi_F$) deserves a thoughtful study. For example, 
when $r=1$, this
is carried by Tate in his thesis. In this sense, what we just ask is a 
non-abelian
generalization of what Tate does. So it would be quite interesting to see 
how our problem is
related Bloch's work on the so-called Tamagawa numbers as well.
\vskip 0.30cm
\noindent
{\bf B.2.3. Geo-Ari Duality and Riemann-Roch: A Practical Geo-Ari 
Cohomology following
Tate}
\vskip 0.30cm
\noindent
{\bf B.2.3.1. An Example}
\vskip 0.30cm
To introduce a more general zeta function, from Artin's
definition of (abelian) zeta function for curves defined over finite fields 
([A]),
and Iwasawa's interpretation of Dedekind zeta function ([Iw]), we see that 
it is
better to have a cohomology theory in arithmetic such that the duality and the
Riemann-Roch  are satisfied. We claim that this can  be
rigorously developed following Tate's Thesis at least in geo-ari dimension one.
To explain the basic idea, we offer the following simplest example.

Consider only rank 1 lattices over {\bf Q}: They are parametrized by
${\bf R}_{>0}$, say the lattices with the forms $\Lambda_t:={\bf Z}\cdot \sqrt
t$, $t\in {\bf R}_{>0}$. For $\Lambda_t$,  the Poisson summation formula 
says that
$$\sum_{n\in {\bf Z}}e^{-\pi tn^2}={1\over {\sqrt t}}
\sum_{n\in {\bf Z}}e^{-\pi n^2/t}.$$ Namely,
$$\sum_{\alpha\in \Lambda_t}e^{-\pi|\alpha|^2}={1\over {{\rm
Vol}(\Lambda_t)}}\sum_{\beta\in \Lambda_t^\vee}e^{-\pi|\beta|^2}.$$
Here $\Lambda_t^\vee$ denotes the dual lattice of
$\Lambda_t$. Thus, if we set
$h^0({\bf Q},\Lambda_t):=\log \big(\sum_{\alpha\in
\Lambda_t}e^{-\pi|\alpha|^2}\big)$, then
$$h^0({\bf Q},\Lambda_t)-h^0({\bf Q}, \Lambda_t^\vee)
={\rm deg}(\Lambda_t).\eqno(*)$$ This is simply the analogue of the
Riemann-Roch  in geometry. Indeed, as for {\bf Q},
the metrized dualizing sheaf is simply the standard lattice ${\bf
Z}\subset {\bf R}={\bf Q}_\infty$. (See e.g. [La2].) So (*) becomes
$$h^0({\bf Q},\Lambda_t)-h^0({\bf Q}, K_{\bf Q}\otimes \Lambda_t^\vee)
={\rm deg}(\Lambda_t)-{1\over 2}{\rm deg}K_{\bf Q},$$
where ${\rm deg}K_{\bf Q}=\log|\Delta_{\bf Q}|=\log 1=0$ with $\Delta_{\bf 
Q}$  the discriminant
of {\bf Q}.
\vskip 0.30cm
\noindent
{\bf B.2.3.2. Canonical Divisors and  Space of Different Forms}
\vskip 0.30cm
  Let $F$ be a number field, i.e.,   a finite extension of {\bf Q}. Denote
by $S=S_F$ the collection of all (unequivalent)
normalized places of $F$. Set $S_\infty$  be the
collection of all Archimedean places in
$S$, and
$S_{\rm fin}:=S\backslash S_\infty$.

For any $v\in S_{\rm fin}$, denote by
$\lambda_0$ the composition of natural morphisms
$${\bf Q}_p\to {\bf Q}_p/{\bf Z}_p\hookrightarrow {\bf Q}/{\bf 
Z}\hookrightarrow
{\bf R}/{\bf Z}.$$  Then we get a natural map
$\lambda_v:F_v\to {\bf R}/{\bf Z}$ defined by $\lambda_v:=\lambda_0\circ {\rm
Tr}_{F_v/{\bf Q}_p}$. Here $p$ is the place of {\bf Q} under $v$, and ${\rm
Tr}_{F_v/{\bf Q}_p}:F_v\to {\bf Q}_p$ denotes the local trace. With this, the
local different $\partial_v$ of
$F$ at
$v$ is characterized by $$\partial_v^{-1}:=\Big\{\alpha_v\in
F_v:\lambda_v(\alpha\cdot {\cal O}_v)=0\Big\}.$$
Moreover, being an ideal of the
discrete valuation ring ${\cal O}_v$ with a parameter $\pi_v$,
$\partial_v=\pi_v^{{\rm ord}_v{(\partial_v)}}\cdot {\cal O}_v$.

By definition, two ideles $a=(a_v)$ and $b=(b_v)$ are called {\it strictly
equivalent}, written as $a\sim_{\rm se}b$, if,
for all $v\in S_\infty$, $a_v=b_v$;
while for all $v\in S_{\rm fin}$, there exists $v$-adic units $u_v$ such that
$a_v=u_vb_v$. By an abuse of notation, denote the associated equivalence class
of $a$ by $[a]_{\rm st},\,[a]$, or even $a$, and denote ${\bf 
I}_F/\sim_{\rm se}$ also
by ${\bf I}_F$.

Now define an {\it idelic canonical element} $\omega_F$ of the number
field $F$ as the (strictly) equivalence class associated to the idele
$(\omega_v)$ of $F$. Here
for each $v\in S_\infty$, $\omega_v:=1$; while for each $v\in S_{\rm fin}$,
$\omega_v:=\pi_v^{-{\rm ord}_v(\partial_v)}$. One checks easily that
$[(\omega_v)]_{\rm st}$ is well-defined. We often call $[(\omega_v)]_{\rm 
st}$ a
{\it canonical divisor} of $F$ as well. For our own convenience,  set
$\pi_v^{\pm {\rm ord}_v(\partial_v)}:=1$ for all $v\in S_\infty$,
despite that we do not have $\pi_v$ when $v\in S_\infty$.

Motivated by the study for function fields, we then define the {\it space
of rational differentials} of $F$ by
$$\Omega_F^1:=\Big\{[(\alpha\cdot \pi_v^{{\rm
ord}_v(\partial_v)})^t]_{\rm st}:
\alpha\in F\Big\}.$$
Here $\cdot^t$ denotes the transpose of $\cdot$.
\vskip 0.30cm
\noindent
{\bf B.2.3.3. Algebraic Cohomology for Matrix Divisors}
\vskip 0.30cm
Let $F$ be  a number field. For every rank $r$ matrix divisor
$g=(g_v)$ of $F$, define its {\it 0-th (cohomology) group}
$$H^0(\overline{{\rm Spec}({\cal
O}_F)},g):=\Big\{\alpha=(\alpha_1,\dots,\alpha_r)^t\in
F^r: g_v^{-1}\cdot \alpha\in ({\cal O}_v)^r,\ \forall v\in S_{\rm
fin}\Big\}.$$ (In this part, for any  set $A$,  let $A_r$ to be
the collection of vectors $(a_1,\dots,a_r)$ with $a_i\in A$ and $A^r$ to
be the collection of vectors $(a_1,\dots,a_r)^t$ with
$a_i\in A$.)

In particular, one sees that, if $r=1$,
$$H^0(\overline{{\rm Spec}({\cal
O}_F)},g):=\Big\{\alpha\in F:g_v^{-1}\cdot \alpha\in {\cal O}_v,\ \forall v\in
S_{\rm fin}\Big\}$$ has the following interpretation.

Let ${\rm Cl}({\cal O}_F)$ denote the ideal class group of $F$. Then there
exists a natural morphism
$$\matrix{\psi:&{\bf I}_F/\sim_{\rm st}&\to& {\rm Cl}({\cal O}_F)\cr
&&&\cr
&g=(g_v)&\mapsto&\prod_{v\in S_{\rm fin}}{\cal P}_v^{{\rm ord}_v(g_v)}.\cr}$$
Here ${\cal P}_v$ denotes the prime ideal of ${\cal O}_F$ corresponding to the
place $v$. One checks easily that $$H^0(\overline{{\rm Spec}({\cal
O}_F)},g)=\psi(g)$$ which is nothing but the global section of the
line bundle
${\cal O}\Big(\sum_v-{\rm ord}_v(g_v)[v]\Big)$ on ${\rm Spec}({\cal O}_F)$.

Next let us define the {\rm 1-st (cohomology) group} of a rank $r$ matrix 
divisor $g$ on $F$. To make the picture more clear, we start with $r=1$. In
this case, $H^1(\overline{{\rm Spec}({\cal O}_F)},g)$ should be a
collection of rational differentials on $F$. Thus, as over function
fields, for a  pre-line bundle $g$ over $F$, naturally we define its
first cohomology group by setting $$H^1(\overline{{\rm Spec}({\cal
O}_F)},g):=\Big\{\beta\in \Omega_F^1: \beta_v\cdot g_v\in {\cal O}_v,\
\forall v\in S_{\rm fin}\Big\}.$$
 
 From this definition, we note that there is a natural isomorphism between
$H^1(\overline{{\rm Spec}({\cal O}_F)},g)$ and
$$\Big\{\alpha\in F:g_v\cdot\alpha\cdot \pi_v^{{\rm ord}_v(\partial_v)}\in 
{\cal
O}_v,\ \forall v\in S_{\rm fin}\Big\}$$ which is simply
$$\Big\{\alpha\in F:\big(\pi_v^{-{\rm ord}_v(\partial_v)}\cdot
g_v^{-1}\big)^{-1}\alpha
\in {\cal O}_v,\ \forall v\in S_{\rm fin}\Big\},$$ i.e.,
$H^0(\overline{{\rm Spec}({\cal O}_F)},\omega_F\otimes g^{-1})$. (Here
and later, the tensor product is defined as usual for matrices.)

We now study how $H^i$'s  depend on  rational equivalence classes.
Assume $g=(g_v)\sim_{\rm ra}g'=(g_v')$. So there exists an $\alpha\in F^*$ such
that $g_v=\alpha\cdot g_v'$ for all $v\in S$. Hence
$$H^0(\overline{{\rm Spec}({\cal O}_F)},g)=\Big\{x\in F:g_v^{-1}\cdot x\in 
{\cal
O}_v,\ \forall v\in S_{\rm fin}\Big\}
=\Big\{x\in F:(g_v')^{-1}\cdot (\alpha^{-1}\cdot x)\in {\cal
O}_v,\ \forall v\in S_{\rm fin}\Big\}.$$
That is to say, for every element $x\in H^0(\overline{{\rm Spec}({\cal
O}_F)},g)$,
$\alpha^{-1}\cdot x$ is an element in $H^0(\overline{{\rm Spec}({\cal
O}_F)},g')$. This then gives an effective isomorphism between
these two 0-th cohomology groups
$$H^0(\overline{{\rm Spec}({\cal O}_F)},g)\buildrel\alpha^{-1}\cdot\over\simeq
H^0(\overline{{\rm Spec}({\cal O}_F)},g').$$
Similarly, we have the canonical isomorphism
$$H^1(\overline{{\rm Spec}({\cal O}_F)},g)\buildrel\alpha\cdot\over\simeq
H^1(\overline{{\rm Spec}({\cal O}_F)},g').$$

With this, we are ready to come back to the general situation, i.e., that
for matrix divisors. By definition, for a rank $r$ matrix divisor $g$ over 
$F$, define its first cohomology groups by setting
$$H^1(\overline{{\rm Spec}({\cal O}_F)},g):=
\Big\{\beta=(\beta_1,\dots,\beta_r)\in
(\Omega_F^1)_r:\beta\cdot g_v^t\in ({\cal O}_v)_r,\ \forall v\in S_{\rm
fin}\Big\}.$$

One chacks easily that the above discussion for rank 1 vector divisors also 
holds for
matrix divisors.
\vskip 0.30cm
\noindent
{\bf Proposition.} {\it With the same notation as above,

\noindent
(1) If $g=\alpha\cdot g'$ for $\alpha\in G(F)$,
$$H^0(\overline{{\rm Spec}({\cal O}_F)},g)\buildrel
\alpha^{-1}\cdot\over\simeq H^0(\overline{{\rm Spec}({\cal O}_F)},g').$$

\noindent
(2) If $g=(g_v)$ denote $g^{-1}:=(g_v^{-1})$ the inverse
of $g$, then
$$H^0(\overline{{\rm Spec}({\cal O}_F)},g)\simeq
\Big(H^1(\overline{{\rm Spec}({\cal O}_F)},g^{-1}\otimes
\omega_F)\Big)^t.$$

\noindent
In particular,
$H^1(\overline{{\rm Spec}({\cal O}_F)},g)$ is canonically
isomorphic to
$$\Big\{\alpha\in F^r:\big(g_v(\alpha)\big)\in \big(\partial_v^{-1}\big)^r,\
\forall v\in S_{\rm fin}\Big\}.$$}
\vskip 0.30cm
\noindent
{\bf B.2.3.4. Geo-Arit Cohomology and Its Associated
Riemann-Roch}
\vskip 0.30cm
It is well-known that in geometry, once we have cohomology groups, naturally,
we use (their ranks or) their dimensions over the base field to define
the corresponding $h^0$ and $h^1$. Yet, for arithmetic setting, we must
do it very differently. (By saying this, we do not mean that the
original geometric counting has no implication in arithmetic
setting: recently Deninger ([D]) proposes a formalism of Betti type
(co)homology theory, where he essentially uses the original geometric
way to count (infinite dimensional spaces). It would be quite
interesting to understand the relation between Deninger's geometric
way of counting and the one used here, which we call an arithmetic
counting over finitely generated cohomology groups $H^0$ and $H^1$.)

Let $F$ be a number field and $g$ be a rank $r$ matrix divisor on $F$,
i.e.,
$g\in {\rm GL}_r({\bf A}_F)$. Then the 0-th cohomology
group and the 1-st cohomology group of $g$ are well-defined. Note that
in particular, if $\beta=(\beta_v)\in H^1(\overline{{\rm Spec}({\cal 
O}_F)},g)$,
then for all $v\in S_\infty$, $\beta_v\in (F_v)_r$ are simply real or
complex
$r$-vectors. With this in mind, we define the geometric arithmetic 
cohomology of
$g$ as follows.

First, define the {\it  0-th geometric arithmetic cohomology}
$h^0(F,g)$ via
$$h^0(F,g):=\log\Big(\sum_{\alpha\in
H^0(\overline{{\rm Spec}({\cal O}_F)},g)}
\exp\Big(-\pi\big(\sum_{v:\ {\rm real}} |
g_v^{-1}\cdot \alpha_v|_v^2+2\sum_{v:\ {\rm complex}}
|g_v^{-1}\cdot \alpha_v|_v^2\big)\Big)\Big),$$ and  define  the {\it
1-st geometric arithmetic cohomology}
$h^1(F,g)$ via
$$h^1(F,g):=\log\Big(\sum_{\beta\in
H^1(\overline{{\rm Spec}({\cal O}_F)},g)}
\exp\Big(-\pi\big(\sum_{v:\ {\rm real}} |
\beta_v\cdot g_v|_v^2+2\sum_{v:\ {\rm complex}} |
\beta_v\cdot g_v|_v^2\big)\Big)\Big).$$
With this, then we obtain the following
\vskip 0.30cm
\noindent
{\bf Proposition.} {\it With the same notation as above,

\noindent
(1) $h^0(F,g)$ and $h^1(F,g)$ are well-defined, i.e., the summations on
the right hand sides are convergent;

\noindent
(2) If $g\sim_{\rm ra} g'$, then
$h^i(F,g)=h^i(F,g')$ for $i=0,1$;

\noindent
(3) $h^0(F,g)=h^1(F,g^{-t}\otimes \omega_F).$}
\vskip 0.30cm
Hence, for a rank $r$ pre vector bundle $[g]$ over $F$, we define
$h^0(F,[g])$ and
$h^1(F,[g])$, {\it the 0-th  and the 1-st
arithmetic cohomology of} $[g]$, to be
$h^0(F,g)$ and
$h^1(F,g)$ respectively, for any representative
$g\in {\rm GL}({\bf A}_F)$ (of $[g]$). By Proposition (2) above,
they are well-defined.

With all this, surely, to state the Riemann-Roch theorem, we still need to
define the degree for a vector bundle. This may be done as over function
fields.  That is to say, if
$g=(g_v)\in {\rm GL}_r({\bf A}_F)$ is a matrix divisor, denote
its determinant by ${\rm det}(g)=({\rm det}(g_v))$,
which is simply an idele of $F$.
Moreover, choose a Haar measure $da$ on ${\bf A}_F$. Then for any
idele $b$ of $F$,
set $N(b)$ or $\|b\|$ to be the unique positive number such that
$d(b\cdot a)=:N(b)\cdot da.$ (This is a global way to understand $N(b)$.
Locally, $N(b)=\|b\|$ may be defined as follows: If $b=(b_v)$, then
$$\|b\|:=\prod_{v\in S}\|b_v\|_v=:\prod_{v\in S}|b_v|_v^{N_v}.$$ Here
$N_v:=[F_v:{\bf Q}_p]$ denotes the local degree of the place of $v$ and $p$
is the place of {\bf Q} under $v$.) Finally, define the {\it degree} of
$g$,  denoted by ${\rm deg}(g)$,
by $${\rm deg}(g):=\log\big(N({\rm det}(g))\big).$$ By using the product
formula, one checks that for any rank $r$ vector bundle $[g]$ of $F$,
$${\rm deg}([g]):={\rm deg}(g)$$ is well-defined. We will call this real number
the {\it degree} of the rank $r$ vector bundle $[g]$.
For example, one checks easily that the
degree of the  canonical divisor $\omega_F$ is simply
$\log|\Delta_F|$, where $\Delta_F$ denotes the discriminant of $F$.
\vskip 0.30cm
\noindent
{\bf Geo-Ari Riemann-Roch Theorem over Number Fields.}
{\it Let $F$ be a number field. Then
for any vector bundle $E$ over $F$, we have
$$\chi_{\rm ga}(F,E):=h^0(F,E)-h^1(F,E)={\rm deg}(E)-{\rm rank}(E)\cdot
{1\over 2}{\rm deg}(\omega_F).$$}

\noindent
One may prove this result as follows following Tate. First, recall  the 
standard
Poisson summation formula to the pair $(F^n, {\bf A}^n)$.

\noindent
{\bf Poisson summation formula.} {\it Let $f$ be continuous and in
$L^1({\bf A}^n)$. Assume that $\sum_{\alpha\in F^n}|f(x+\alpha)|$ is
uniformly convergent for $x$ in a compact subset of $ {\bf A}^n$, and
that $\sum_{\alpha\in F^n}\hat f(\alpha)$ is convergent, where $\hat f$
denotes the Fourier transform of $f$.  Then
$$\sum_{\alpha\in F^n}\hat f(\alpha)=\sum_{\alpha\in F^n}
f(\alpha).$$}

But, for any element $g\in {\rm GL}_r({\bf A}_F)$, set $h(x):=f(gx)$,
then one checks that $\hat h(x)={1\over {\|{\rm det}(g)\|}}\hat
f(g^{-t} x).$
Thus, from the Poisson summation formula above, we get the
following more suitable version for our application:
$${1\over {\|{\rm det}(g)\|}}\sum_{\alpha\in F^r}\hat f(g^{-t}\cdot
\alpha) =\sum_{\alpha\in F^r} f(g\cdot \alpha).$$

Now set $f(x):=\prod_{v\in S}f_v(x_v)$ with $f_v$ as follows:

\noindent
(i) If $v$ is real, then $f_v(x_v):=e^{-\pi |x_v|^2}$, where
for $x_v=(x_v^{(1)},\dots, x_v^{(r)})$, $|x_v|^2=\sum_i(x_v^{(i)})^2$.
Obviously, we have $f_v=\hat f_v$;

\noindent
(ii) If $v$ is complex, then $f_v(x_v):=e^{-2\pi |x_v|^2}$, where
for $x_v=(x_v^{(1)},\dots, x_v^{(r)})$, $|x_v|^2=\sum_ix_v^{(i)}\cdot
{\overline {x_v^{(i)}}}$. Similarly, $f_v=\hat f_v$;

\noindent
(iii) If $v$ is finite, then $f_v$ is defined to be the characteristic
function of $(\partial_v^{-1})^r$. One checks that $\hat f_v$
equals to $(N(\partial_v))^{r/2}$ times the characteristic function of
${\cal O}_v^r$.
 
With this, what we have is the following formula:
$$\eqalign{{1\over {\|{\rm det}(g)\|}}&\cdot (N(\omega_F))^{r/2}\cdot
\sum_{g_v^{-t}(\alpha)\in {\cal O}_v^r,\forall
v\in S_{\rm fin}}\Big(\prod_{v:\
{\rm real}}  e^{-\pi |g_v^{-t}(
\alpha)|^2}\prod_{v:\ {\rm complex}}
e^{-2\pi |g_v^{-t}(
\alpha)|^2}\Big)\cr =&\sum_{g_v(\alpha)\in (\partial_v^{-1})^r ,\forall
v\in S_{\rm fin}}\Big(\prod_{v:\
{\rm real}}  e^{-\pi |g_v(\alpha)|^2}\prod_{v:\ {\rm complex}}
e^{-2\pi |g_v( \alpha)|^2}\Big).\cr}$$
This then gives a proof of the theorem by definition and  the Proposition 
2.3.3.
\vskip 0.30cm
\noindent
{\it Remarks.} (1) The practical geo-ari cohomology and Riemann-Roch 
theorem here have
indeed a theoretical treatment along the line in C.3. We are going to 
include this in
the yet to be released second version of [We2].

\noindent
(2) We may understand van der Geer and Schoof's
work as follows:

\noindent
(i) For any idele $a$, define the associated Arakelov
divisor ${\rm div}_{\rm Ar}(a)$ as follows:
$${\rm div}_{\rm Ar}(a):=-\sum_{v\in S_{\rm fin}}{\rm ord}_v(a_v)\cdot [v]
+\sum_{v\in S_\infty}N_v\log|a_v|\cdot [v].$$ Obviously,
all Arakelov divisor may be
constructed in this way.

\noindent
(ii) Define $H^0(\overline{{\rm Spec}{\cal O}_F},{\rm div}_{\rm
Ar}(a)):=H^0(\overline{{\rm Spec}{\cal O}_F},a)$ and
setting $h^0(\overline{{\rm Spec}{\cal O}_F},{\rm div}_{\rm
Ar}(a)):=h^0(\overline{{\rm Spec}{\cal O}_F},a)$. One checks then that this
definition coincides with that of van der Geer and Schoof. Hence, we have
also;

\noindent
{\bf van der Geer-Schoof's Riemann-Roch theorem:}
{\it For any Arakelov divisor $D$ over $\overline{{\rm Spec}({\cal
O}_F)}$, $$h^0(\overline{{\rm
Spec}({\cal O}_F)},D)-h^0(\overline{{\rm Spec}({\cal
O}_F)},K_{{\rm Spec}({\cal O}_F)}-D)={\rm deg}(D)-
{1\over 2}{\rm deg}(\omega_F).$$}
\vskip 0.30cm
\noindent
{\bf B.2.4. Non-Abelian Zeta Function For Number Fields}
\vskip 0.30cm
\noindent
{\bf B.2.4.1. The Construction}
\vskip 0.30cm
Let $F$ be a number fields with discriminant $\Delta_F$. Denote by ${\cal 
M}_{{\bf
A}_F,r}(d)$ the (adelic) moduli space of rank $r$ and degree $d$ 
semi-stable pre vector
bundles over $F$, and $d\mu$ its associated Tamagawa measure. Then we define
the rank $r$  non-abelian (completed)  zeta function $\xi_{F,r}(s)$ of $F$ by
$$\xi_{F,r}(s):=\Big(|\Delta_F|^{{r\over
2}}\Big)^s\int_{g\in {\cal M}_{{\bf A}_F,r}(t),t\in {\bf R}_{>0}}
\Big(e^{h^0({\bf A}_F,g)}-1\Big)\big(e^{-s}\big)^{{\rm
deg}(g)}\cdot d\mu(g),\qquad {\rm Re}(s)>1.$$

Following 2.1, i.e., Iwasawa's interpretation of Dedekind zeta functions, 
we have
$\xi_{F,1}(s)=w_F\cdot \xi_F(s)$ where $w_F$ denotes the number
of roots of unity in $F$ and $\xi_F(s)$ denotes the completed Dedekind zeta 
function for $F$.
Note also that  the terms apeared in our construction above, such as the 
degree,
the geo-ari cohomology $h^0$, the moduli space, and the Tamagawa measure 
are all
canonically and naturally associated with number fields. Hence, our 
non-abelian zeta
should be genuinely related with non-abelian arithmetic properties of 
number fields.
\vskip 0.30cm
\noindent
{\it Remark.} Recall that there is an algebraic moment map $\pi_F:{\cal 
M}_{{\bf
A}_F,r}(d)\to {\cal M}_{F,r}(d)$. Thus our above construction of the zetas 
may be understood
as a kind of algebraic version of Feynman type integral. On the other hand, 
in Part (A), we
propose a Weil-Narasimhan-Seshadri type correspondence, namely, a micro 
reciprocity law.
So it is not unreasonable to expect  that our zetas may be written in terms 
of analytic
Feynman type integrals, and that the global (non-abelian) reciprocity law 
may be obtained
from our zeta functions. I would like to thank Nitta and Okada for their 
discussion here.
\vskip 0.30cm
\noindent
{\bf B.2.4.2.  Basic Properties}
\vskip 0.30cm
Just as for Dedekind zeta functions, our non-abelian zeta functions are 
well-defined, satisfy
functional equation as well. Moreover,  the residues of these zeta 
functions may be
calculated in terms of the volumes of the moduli space. More precisely, we 
have the following
\vskip 0.30cm
\noindent
{\bf Theorem.} {\it With the same notation as above, we have

\noindent
(1) $$\xi_{F,r}(s)=\Big(|\Delta_F|^{{r\over
2}}\Big)^s\int_{g\in {\cal M}_{{\bf A}_F,r}(t),t\in {\bf R}_{>0}}
\Big(e^{h^0({\bf A}_F,g)}-1\Big)\big(e^{-s}\big)^{{\rm
deg}(g)}\cdot d\mu(g)$$ converges
absolutely and uniformly when ${\rm Re}(s)\geq 1+\delta$ for any
$\delta>0$;

\noindent
(2) $\xi_{F,r}(s)$ admits a unique meromorphic continuation to the
whole complex $s$-plane with only two simple poles at $s=0,1$ whose
residues are  ${\rm Vol}({\cal M}_{{\bf A}_F,r}(t))$
for one and hence for all $t$;

\noindent
(3) ({\rm Functional Equation})  $\xi_{F,r}(s)=\xi_{F,r}(1-s)$.}
\vskip 0.30cm
\noindent
{\it Remark.} Most suitable definition for non-abelian zeta functions of 
number fields should be
$$\xi_{F,r}(s):=\Big(|\Delta_F|^{{r\over
2}}\Big)^s\int_{\Lambda\in {\cal M}_{{\bf A}_F,r}(t),t\in {\bf R}_+}
{{e^{h^0({\bf A}_F,\Lambda)}-1}\over{\#{\rm
Aut}(\Lambda)}}\big(e^{-s}\big)^{{\rm deg}(\Lambda)}
d\mu(\Lambda),\qquad {\rm Re}(s)>1$$ where ${\rm Aut}$ denotes the
  automorphism group. Moreover,  instead of using the adelic moduli spaces,
we may introduce the \lq standard' version of non-abelian zeta functions 
using integrations
over moduli spaces of semi-stable lattices. In this way, we may then also 
see what are the
\lq Gamma'-factors and a non-abelian version of Tate's calculation on the 
so-called
analytic class number formula.
\vfill
\eject
\vskip 0.45cm
\centerline{\we C. Explicit Formula, Functional Equation and Geo-Ari 
Intersection}
\vskip 0.45cm
\noindent
{\li C.1. The Riemann Hypothesis for Curves}
\vskip 0.45cm
\noindent
{\bf C.1.1. Weil's Explicit Formula: the Reciprocity Law}
\vskip 0.30cm
Let $C$ be a projective irreducible reduced regular curve of genus
$g$ defined over a finite field $k:={\bf F}_q$. Denote by $\zeta_C(s)$ the 
associated Artin
zeta function. Set $t=q^{-s}$ and $Z_C(t):=\zeta_C(s)$. Then by the 
rationality,
there exists a polynomial
$P_C(t)$ of degree $2g$ such that  $$Z_C(t)={{P_C(t)}\over{(1-t)(1-qt)}}.$$.

Now let $C_n$ be the curve obtained
from $C$ by extending the field of constants from ${\bf F}_q$ to ${\bf 
F}_{q^n}$.
Then by a discussion on covering of curves, we obtain the following
\vskip 0.30cm
\noindent
{\bf Reciprocity Law.} {\it With the same notation as above, 
$Z_{C^n}(t^n)=\prod_\zeta
Z_C(\zeta t)$ where the product is taken over all the $n$-th roots of 1.}
\vskip 0.30cm
Moreover, by the Euler product,
$$t{{Z_C'}\over Z_C}(t)=\sum_{n=1}^\infty\sum_Pd(P)t^{nd(P)}
=\sum_{m=1}^\infty\Big(\sum_Pd(P)\Big) t^m$$ where $\sum_P$ is taken over those
closed points rational over ${\bf F}_q$ whose degree divides $m$. Hence,
$$Z_C(t)=\exp \Big\{\sum_{m=1}^\infty N_m{{t^m}\over m}\Big\}$$ where
$N_m=\sum_{P,d(P)|m}d(P)$. Clearly when $m=1$, the sum $\sum_{P,d(P)|1}1$ 
simply
counts the number of closed points of $C$ rational over ${\bf F}_q$.
Thus, by the reciprocity law above, $$Z_{C^n}(t^n)=\exp\Big\{\sum_{m=1}^\infty
M_m{{t^{nm}}\over m}\Big\}=
\exp \Big\{\sum_{m=1}^\infty N_m{{t^m}\over 
m}\Big(\sum_\zeta\zeta^m\Big)\Big\}.$$
This implies that  $N_n=N_1(C_n)$, i.e., $N_n$ is the number of closed 
points on $C$ which
are rational over ${\bf F}_{q^n}$.  Therefore, we have the following
\vskip 0.30cm
\noindent
{\bf Weil's Explicit Formula.} {\it With the same notation as above,
$$N_n=q^n+1-\sum_{\zeta_C(\rho)=0}\rho^n,$$
where $\rho_1,\dots,\rho_{2g}$ denotes the reciprocals of the roots of 
$P_C(t)$.}
\vskip 0.30cm
\noindent
{\bf C.1.2. Geometric Version of Explicit Formula}
\vskip 0.30cm
As above, let $C$ be an algebraic curve defined over ${\bf F}_p$, the 
finite field with
$p$ elements. Over $C\times C$, for $n\in {\bf Z}$, introduce  (micro) divisors
$A_n$ (via algebraic correspondence) as follows:
$$A_n:=\cases{\Big\{(x,x^{p^n}):x\in C\Big\}, &if\ $n\geq 0$;\cr
\Big\{(x^{p^{-n}},x):x\in C\Big\}, &if\ $n\leq 0$.\cr}$$
Clearly, we have the following relations for the intersections among $A_n$'s.

\noindent
(i) If $n\geq m\geq 0$, $$\langle A_n,A_m\rangle =p^m\langle A_{n-m},
A_0\rangle;$$

\noindent
(ii) If $m\geq n\geq 0$, $$\langle A_n,A_m\rangle =p^n\langle A_{m-n},
A_0\rangle;$$

\noindent
(iii) If $m\geq 0\geq n$, $$\langle A_n,A_m\rangle =\langle A_{n-m},
A_0\rangle;$$

\noindent
(iv)  If $n\leq m\leq 0$, $$p^m\langle A_n,A_m\rangle =\langle A_{m-n},
A_0\rangle;$$

\noindent
(v) If $m\leq n\leq 0$, $$p^n\langle A_n,A_m\rangle =\langle A_{n-m},
A_0\rangle;$$

\noindent
(vi) If $n\geq 0 \geq m$, $$\langle A_n,A_m\rangle =\langle A_{n-m},
A_0\rangle.$$

Hence, we have the following
\vskip 0.30cm
\noindent
{\bf Lemma 1.} {\it With the same notation as above, we have

\noindent
(1) For all $m,n\in {\bf Z}$, $$\langle A_n,A_m\rangle =\langle
A_m,A_n\rangle;$$

\noindent
(2) For all $m,n\in {\bf Z}$, $$\langle A_{-n},A_{-m}\rangle =\langle
A_m,A_n\rangle;$$

\noindent
(3) For all $n\geq m\geq 0$ in {\bf Z}, $$\langle A_n,A_m\rangle =p^m\langle
A_{n-m},A_0\rangle;$$

\noindent
(4) For all $n\geq 0\geq m$ in {\bf Z},
$$\langle A_n,A_m\rangle =\langle
A_{n-m},A_0\rangle.$$}

Obviously, (1) $\sim$ (4) are equivalent to (i) $\sim$ (vi) above. Therefore,
in order to understand $\langle A_n,A_m\rangle$ for all $n,m\in {\bf Z}$, we
only need to know $\langle A_n,A_0\rangle$ for $n\geq 0$.

For this latest purpose, first, note that $A_0$ is simply the diagonal. 
Hence, by definition,
for $n\,>\,0$, $\langle A_n,A_0\rangle$ is  $N_n(C)$,
the number of closed points of $C$ which are rational over ${\bf F}_{p^n}$;

Secondly, by the functional equation of the zeta function $\zeta_C(s)$ of $C$,
which itself is a direct consequence of the Riemann-Roch theorem for the curve
$C$, (and the rationality of $\zeta_C(s)$), we have the following
\vskip 0.30cm
\noindent
{\bf Lemma 2.} ({\bf Explicit Formula of Weil}) {\it With the same notation as
above, if $n\in {\bf Z}_{\geq 0}$,
$$\langle A_n,A_0\rangle=\langle A_n,F_1\rangle+\langle A_n,F_2\rangle
-\sum_{\zeta_C(s)=0}s^n.$$ Here $F_1$ and $F_2$ denotes the fibers in two
directions of $C\times C$ respectively, and  the sum is taken over all zeros of
$\zeta_C(s)$.}
\vskip 0.30cm
\noindent
{\it Remark.} From now on, we may from place to place have some sign 
problems, say
$s^n$ may well mean $s^{-n}$.
\vskip 0.30cm
\noindent
{\bf C.1.3. Riemann Hypothesis for Function Fields}
\vskip 0.30cm
Now, following Weil again, we prove the (Artin-)Riemann hypothesis. (See 
e.g., [Ha].)

Let $f:p^{\bf Z}\to {\bf Z}$ be a function with finite supports. Define its 
Mellin
transform via $\hat f(s):=\sum_nf(p^n)p^{ns}.$ Clearly, if
$f^*(p^n):=f(p^{-n})p^{-n}$,
$\widehat{f^*}(s)=\hat f(1-s)$; moreover, if
$(f*g)(p^n):=\sum_mf(p^m)g(p^{n-m})$,  $\widehat{f*g}(s)=\hat f(s)\cdot \hat
g(s).$ In particular, $\widehat{f*g^*}(s)=\hat f(s)\cdot \hat
g(1-s).$

With $f$, we may  use the divisors $A_n$'s above to define a (global) {\bf 
Q}-divisor
$D_{\hat f}$ on $C\times C$ as follows:
$$D_{\hat f}:=\sum_{n>0}f(p^n)A_n+\sum_{n\geq
0}f(p^{-n})p^{-n}A_n.$$
\vskip 0.30cm
\noindent
{\bf Theorem.} (Weil) {\it With the same notation as above,

\noindent
(1) (Relative Degrees)
$$\langle D_{\hat
f},F_1\rangle=\hat f(1);\qquad
\langle D_{\hat
f},F_2\rangle=\hat f(0);$$

\noindent
(2) (Fixed Points Formula)
$$\langle D_{\hat f},D_{\hat g}\rangle=\langle
D_{\widehat {f*g^*}},{\rm Diag}\rangle,$$
where Diag denotes the diagonal of $C\times C$;

\noindent
(3) (Explicit Formula)
$$\hat f(0)+\hat f(1)-\sum_{\zeta(s)=0}\hat f(s)
=\langle
D_{\hat {f}},{\rm Diag}\rangle.$$}

\noindent
Indeed, by definition,
$$\eqalign{\langle D_{\hat
f},F_1\rangle=&\sum_{n>0}f(p^n)\cdot p^n+\sum_{n\geq 0}
f(p^{-n})p^{-n}\cdot 1=\hat f(1);\cr
\langle D_{\hat
f},F_2\rangle=&\sum_{n>0}f(p^n)+\sum_{n\geq 0}
f(p^{-n})p^{-n}\cdot p^n=\hat f(0).\cr}$$ This gives
(1). (2) is a direct consequence of the relations (i) $\sim$ (vi) for
the intersections of $A_n$'s, and hence comes from Lemma 12.1.
Finally, (3) is simply Lemma 1.2.2 by definition.
\vskip 0.30cm
Next we apply the two dimensional intersection theory, in particular, the Hodge
Index Theorem. Since
$$\langle F_1+F_2,D_{\hat f}-\hat f(1)F_2-\hat f(0)F_1\rangle=0,$$
$$\langle D_{\hat f}-\hat f(1)F_2-\hat f(0)F_1,D_{\hat
f}-\hat f(1)F_2-\hat f(0)F_1\rangle\leq 0.$$
That is,
$$\hat f(0)\cdot\hat f(1)\geq {1\over
2}\langle D_{\hat f},D_{\hat f}\rangle.$$ Thus by Theorem above,
this last equality is equivalent to
$$\sum_{\zeta(s)=0}\hat f(s)\cdot\hat f(1-s)\geq
0.$$
 From this, easily, we get the following  Riemann Hypothesis for Artin zeta 
functions of
curves over finite fields.
\vskip 0.30cm
\noindent
{\bf Theorem.} (Hasse-Weil)
{\it Let $\zeta(s)$ be the zeta function for a curve
defined over a finite field. If $\zeta(s)=0$,
then ${\rm Re}(s)={1\over 2}.$}
\vskip 0.45cm
\noindent
{\li C.2. Geo-Ari Intersection in Dimension Two: A Mathematics Model}
\vskip 0.45cm
\noindent
{\bf C.2.1. Motivation from Cram\'er's Formula}
\vskip 0.30cm
In the above discussion, the summation $\sum_{\zeta_C(s)=0}s^n$ plays a key
role in understanding Artin-Riemann Hypothesis, via the so-called micro 
explicit formula of
Weil. So  naturally, we want to know whether this approach works for number 
fields, and
are led to study the formal summation
$$\sum_{\xi_{\bf Q}(\rho)=0} x^\rho,\qquad {\rm for}\ x\in 
[1,\infty).\eqno(*)$$
Here $\xi_{\bf Q}(s)$ denotes the completed Riemann zeta function, i.e.,
$\xi_{\bf Q}(s):=\pi^{-{s\over 2}}\Gamma({s\over 2})\zeta_{\bf Q}(s)$ with
$\Gamma(s)$ the standard Gamma function and $\zeta_{\bf Q}(s)$ the Riemann zeta
function.

The reader at this point certainly would reject $(*)$,
since the summation does not make any sense. How could I write such a 
monster down?!
Well, do not be so panic!!! After all,
in the study of the prime distributions, the conditional convergent summation
$\sum_{\xi_{\bf Q}(\rho)=0}{{x^\rho}\over\rho}$ does appear. Recall that we 
have the
following Riemann-von Mangoldt formula
$$\sum_{p^n\leq x}'\log p=x-\sum_{\xi_{\bf Q}(\rho)=0}{{x^\rho}\over
\rho}-{{\zeta'}\over \zeta}(0)-{1\over 2}\log(1-x^2),$$ which itself motivats
and hence stands as a special form of the more general form of explicit
formulas. See e.g., Jorgenson and Lang's lecture notes on Explicit Formulas.
More generally, in various discussions about prime distributions, we do use
the summations such as
$\sum_{\xi_{\bf Q}(\rho)=0}{{x^{\rho\pm 1}}\over {\rho(\rho\pm 1)}}$. In 
this sense, the
problem is not whether we should introduce $(*)$, rather,
it should be how to justify it.

Anyway, let me make a change of variables $x:=e^t$ with $t\geq
0$. Then $(*)$ becomes $$V(t):=\sum_{\xi_{\bf Q}(\rho)=0}e^{t\rho}.\eqno(**)$$
So, following Riemann, we may further view $V(t)$ as  a function of complex 
variable $z$, i.e.,
$$V(z):=\sum_{\xi_{\bf Q}(\rho)=0}e^{z\rho}.\eqno(**')$$ Now, we claim that
there is a nice way to regularize $(*)$.

To explain this in a simpler form, which in fact would not really make our 
life any
easier, we instead consider the partial sum $V_+(z)$ defined by
$$V_+(z):=\sum_{\xi_{\bf Q}(\rho)=0,{\rm Re}(\rho)>0}e^{z\rho}.$$
Then we have the following result, dated in 1919.
\vskip 0.30cm
\noindent
{\bf Theorem.} (Cram\'er) {\it The function
$V_+(z)$ converges absolutely for ${\rm Im}(z)>0$. Moreover,
$$2\pi i V_+(z)-{{\log z}\over {1-e^z}}$$ has a meromorphic continuation to
{\bf C}, with simple poles at the points $\pm\pi \pi n$ for all integers 
$n$, and
at the points $\pm \log p^m$ for all powers of primes.}
\vskip 0.30cm
We claim that this theorem acturally offers us a natural analytic way to 
normalize
the formal summation $\sum_{\xi_{\bf Q}(\rho)=0} x^\rho$ for $x\in 
[1,\infty)$, by using
[C], [JL1,2,3], and [DS]. In a sense, this is in a similar way as  what we 
do when normalizing
$\infty!$: By the Stirling formula
$$n!=\sqrt{2\pi}\cdot\sqrt {n}\cdot\Big({n\over
e}\Big)^n\cdot\exp\Big({{\theta_n}\over {12}}\Big)$$ for $n$ sufficient large
with $|\theta_n|<1$, we set $\infty!=\sqrt{2\pi}$.

However,  in this article, we  do it very differently -- We are going to
construct a mathematics model  to  normalize this formal
summation geometrically.
\vskip 0.30cm
\noindent
{\bf C.2.2. Micro Divisors}.
\vskip 0.30cm
We will not use Grothendieck's scheme language.
Instead, formally, we  call (the set theoretical product)
$S:=\overline {{\rm Spec}({\bf Z})}\times \overline {{\rm Spec}({\bf Z})}$ a
geometric arithmetic base (surface).

Similarly as in 1.2, associated to all $x\in [0,\infty]$ are symbols $D_x$ 
which will be
called  {\it micro divisors}. Associated to any two micro divisors
$D_x,\, D_y$ is  the intersection
number $\langle D_x,D_y\rangle\in {\bf R}$. Assume that the
following fundamental relations
are satisfied by our micro intersections:
\vskip 0.30cm
\noindent
(1) ({\it Symmetry}) For any $x,y\in [0,\infty]$,
$$\langle D_x,D_y\rangle=\langle D_y,D_x\rangle;$$

\noindent
(2) ({\it Mirrow Image}) For any $x,y\in [0,\infty]$,
$$\langle D_x,D_y\rangle=\langle D_{1\over x},D_{1\over y}\rangle;$$

\noindent
(3) ({\it Fixed Points 1}) If $0\leq x\leq y\leq 1$, then
$$\langle D_x,D_y\rangle=y\langle D_{x\over y},D_1\rangle;$$

\noindent
(4) ({\it Fixed Points 2}) If $0\leq x\leq 1\leq y\leq \infty$, then
$$\langle D_x,D_y\rangle=\langle D_{x\over y},D_1\rangle.$$

Note that
$[0,\infty]\times [0,\infty]$ is simply the union of
$\{0\leq x\leq y\leq 1\}$, $\{0\leq y\leq x\leq 1\}$,
$\{1\leq x\leq y\leq \infty\}$, $\{1\leq y\leq x\leq \infty\}$,
$\{0\leq x\leq 1, 1\leq y\leq \infty\}$
and  $\{0\leq y\leq 1, 1\leq x\leq \infty\}$.
Thus by the above relations, if we define the pricise intersection
$\langle D_x,D_1\rangle$ for all $x\in [0,1]$, then we have all the
intersections $\langle D_x,D_y\rangle$ for all $x,y\in [0,\infty]$.
\vskip 0.30cm
\noindent
(5) ({\it Explicit Formula 1}) Denote the completed Riemann
zeta function by $\xi_{\bf Q}(s)$. Then, for all $x\in [0,1]$,
$$\langle D_x,D_1\rangle=\langle D_0,D_x\rangle+\langle D_\infty,D_x\rangle
-\sum_{\xi_{\bf Q}(s)=0}x^s.$$

Here, as said in 2.1, we certainly encount with the convergence
problem of $\sum_{\xi_{\bf Q}(s)=0}x^s$. Instead of solving it, in our model,
we simply assume that our micro
intersection offers a natural normalization of
$\sum_{\xi_{\bf Q}(s)=0}x^s$  via the (5). It is in this sense  we  say that
our model gives a geometric way to normalize
$\sum_{\xi_{\bf Q}(s)=0}x^s$.
\vskip 0.30cm
\noindent
{\it Remark.} The compactibility among
($i$)'s, $i=1,2,3,4,5$, is guaranteed by the functional equation for Riemann
zeta function. Moreover, if $x\geq 1$, then by the Mirrow principal, we see 
that
$\langle D_x,D_1\rangle=\langle D_{1\over x},D_1\rangle$. So, by using the 
Explicit
Formula 1, i.e., (5) above, together with the Relations II below, we have
$$\langle D_{1\over x},D_1\rangle=1+{1\over x}-\sum_{\xi_{\bf 
Q}(s)=0}x^{-s}.$$ Multiplying
both sides by $x$, we  get
$$x\langle D_{1\over x},D_1\rangle=1+x-\sum_{\xi_{\bf Q}(s)=0}x^{1-s}=
1+x-\sum_{\xi_{\bf Q}(s)=0}x^s,$$  by the functional equation. On the other 
hand,
by the Fixed Point 1, i.e., (3) above, we get
$$\langle D_x,D_1\rangle=x\langle D_{1\over x},D_1\rangle.$$
That is to say, formally, for all $x$,
$$\langle D_x,D_1\rangle=1+x-\sum_{\xi_{\bf Q}(s)=0}x^s.$$

 From (1), (2), (3) and (4), we may formally get the following
relations.
\vskip 0.30cm
\noindent
{\bf Relations I.} {\it (i) $\langle D_0,D_0\rangle =\langle
D_\infty,D_\infty\rangle$;

\noindent
(ii) $\langle D_0,D_1\rangle =\langle D_\infty,D_1\rangle=\langle D_0,D_\infty
\rangle$.}
\vskip .30cm
Note that $D_0$ and $D_\infty$ are supposed to be the fibers in two 
directions of $\overline
{{\rm Spec}({\bf Z})}\times \overline {{\rm Spec}({\bf Z})}$
  over $\infty\in \overline {{\rm Spec}({\bf Z})}$, so
  we normalize our intersection further by the following
\vskip 0.30cm
\noindent
(6) ({\it Normalization 1}) $\langle D_0,D_0\rangle=0,\quad \langle 
D_0,D_1\rangle=1$.
\vskip 0.30cm
With this, formally we may further get the following
\vskip 0.30cm
\noindent
{\bf Relations II.} {\it (i) For $\langle D_0,D_x\rangle$,
$$\langle D_0,D_x\rangle =\cases{x,&if $x\in [0,1]$;\cr
1,&if $x\in[1,\infty]$;\cr}$$

\noindent
(ii) For $\langle D_\infty,D_x\rangle$,
$$\langle D_\infty,D_x\rangle =\cases{1,&if $x\in [0,1]$;\cr
{1\over x},&if $x\in[1,\infty]$;\cr}$$

\noindent
In particular,  for $x\in [0,1]$,
$$\langle D_x,D_1\rangle=1+x-\sum_{\xi_{\bf Q}(s)=0}x^s.$$}
 
\noindent
{\it Remark 2.} In (iii) above, taking $x=1$, we  get
$$\langle D_1,D_1\rangle =2 + \sum_{\xi_{\bf Q}(s)=0} 1.$$ Thus,
via the so-called Adjunction Formula, which should be one of the fundamental
results for our intersection, $\sum_{\xi_{\bf Q}(s)=0} 1$
is supposed to be related to the
{\it canonical divisor} for our aritmetic dimension one base.
\vskip 0.30cm
\noindent
{\bf C.2.3. Global Divisors and Their Intersections: Geometric Reciprocity Law}
\vskip 0.30cm
Motivated by C.1, i.e., the discussion about Artin-Riemann Hypothesis, we start
with a standard construction in function theory. Let $f:{\bf R}^+\to  {\bf R}$
be a smooth, compactly supported function. Define its {\it Mellin 
transform} via
$$\hat f(s):=\int_0^\infty f(x) x^s {{dx}\over x}.$$ Then,  if 
$f^*(x):=f({1\over
x})\cdot {1\over x}$,
$$\widehat {f^*}(s)=\hat f(1-s);$$ and
if $(f*g)(x):=\int_0^\infty f(y)g({x\over y})
{{dy}\over y}$ denotes the standard multiplicative convolution,
$$\widehat{f*g}(s)=\hat f(s)\cdot\hat g(s).$$ In particular,
$$\widehat{f*g^*}(s)=\hat f(s)\cdot\hat g(1-s).$$

\vskip 0.30cm
Next, we give a parallel construction for our divisors.
Standard wishdom says that we should use linear combinations of generalized
divsors $D_x, x\in [0,\infty]$ to form new type of divisors. But as we
clearly see that such a conventional way does not result sufficiently many
divisors, we do it very differently.

Starting from divisors $D_x$, for any might-be-interesting function $f$,
formally define the associated {\it global divisor} $D_{\hat f}$  by setting
$$D_{\hat f}:=\int_0^1 f(x)\cdot D_x\cdot {{dx}\over x}
+\int_1^\infty f(x)\cdot xD_x\cdot {{dx}\over x}.$$

\noindent
{\it Remark.} In all formal discussions here, we pay no attention to the 
convergence
problem. See however 2.5 below.  In other words, we assume also in our 
model that global
divisors do exist and their relations with  micro divisors are given as above.

With this definition, we may extend the intersection in 2.2 to $D_{\hat 
f}$'s by
linearity. For example, the intersection $\langle D_{\hat f},D_1\rangle$
is by definition given by
$$\langle D_{\hat f},D_1\rangle=\int_0^1 f(x)\cdot \langle D_x,D_1\rangle\cdot
{{dx}\over x} +\int_1^\infty f(x)\cdot x\langle D_x,D_1\rangle\cdot {{dx}\over
x}.$$ Then formally, we have the following
\vskip 0.30cm
\noindent
{\bf Key Relations.} {\it With the same notation as above,

\noindent
(i)  (Relative Degrees in Two Fiber Directions)
$${\rm deg}_1D_{\hat f}:=\langle D_0, D_{\hat f}\rangle =\hat f(1)$$
and
$${\rm deg}_2D_{\hat f}:=\langle D_\infty, D_{\hat f}\rangle =\hat f(0).$$

\noindent
(ii) (Fixed Point Formula)
$$\langle D_{\hat f},D_{\hat g}\rangle=\langle 
D_{\widehat{f*g^*}},D_1\rangle;$$

\noindent
(iii) (Explicit Formula 2)
$$\langle D_{\hat f},D_1\rangle
=\hat f(0)+\hat f(1)-\sum_{\xi_{\bf Q}(s)=0}\hat f(s).$$}

\noindent
In fact, this may be formally checked as follows using axioms. First 
consider the Fixed Point
Formula. Set $D_x':=D_{1\over x}$. Then,
by definition, $$\eqalign{~&\langle
D_{\widehat{f*g^*}},D_1\rangle\cr
=&\langle\int_0^1\widehat{f*g^*}(x)D_x{{dx}\over x}+
\int_0^1\widehat{f*g^*}({1\over x}){1\over x}D_x'{{dx}\over x},D_1\rangle\cr
=&\langle\int_0^1\int_0^\infty f(y)g^*({x\over
y}){{dy}\over y} D_x{{dx}\over x}+
\int_0^1\int_0^\infty f(y)g^*({1\over {xy}}) xy{{dy}\over y}
{1\over x}D_x'{{dx}\over x},D_1\rangle\cr}$$ By changing variables 
$x':={y\over x}$,
the latest quantity is simply
$$=\langle\int_0^\infty f(y){{dy}\over y}\int_y^\infty g(x)x D_{y\over 
x}{{dx}\over x}
+\int_0^\infty f(y){{dy}\over y}\int_0^yg(x) yD_{x\over y} {{dx}\over x}.$$
But for $x\geq y$, we may split the region into three parts, i.e.,

\noindent
(1.1) $1\geq x\geq 0, 1\geq y\geq 0$ and $x\geq y$;

\noindent
(1.2) $\infty\geq x\geq 1, \infty\geq y\geq 1$ and $x\geq y$;

\noindent
(1.3) $\infty\geq x\geq 1$ and $1\geq y\geq 0$.

Similarly, the region of $x<y$ may be split into three parts, i.e.,

\noindent
(2.1) $1\geq x\geq 0, 1\geq y\geq 0$ and $x< y$;

\noindent
(2.2) $\infty\geq x\geq 1, \infty\geq y\geq 1$ and $x< y$;

\noindent
(2.3) $\infty\geq y\geq 1$ and $1\geq x\geq 0$.

Hence the latest quantity is simply, by writing according to (2.1),(2.2) 
and (2.3)
(resp. (1.1), (1.2) and (1.3))
for the first term (resp. second term),
$$\eqalign{~&\langle\Big( \int_0^1f(y){{dy}\over y}\int_0^yg(x)yD_{x\over 
y}{{dx}\over
x} +\int_1^\infty f(y){{dy}\over y}\int_1^yg(x)yD_{x\over y}{{dx}\over x}
+\int_1^\infty f(y){{dy}\over y}\int_0^1g(x)yD_{x\over y}{{dx}\over x}\Big)\cr
&+\Big(\langle \int_0^1f(y){{dy}\over y}\int_y^1g(x)xD_{y\over x}{{dx}\over x}
+\int_1^\infty f(y){{dy}\over y}\int_y^\infty g(x)xD_{y\over x}{{dx}\over x}
+\int_0^1 f(y){{dy}\over y}\int_1^\infty g(x)xD_{y\over x}{{dx}\over
x}\Big),D_1\rangle\cr
=&\langle\Big( \int_0^1f(y){{dy}\over y}\Big[\int_0^yg(x)yD_{x\over 
y}{{dx}\over x}
+\int_y^1g(x)xD_{y\over x}{{dx}\over x}\Big]+\int_1^\infty f(y){{dy}\over
y}\Big[\int_1^yg(x)yD_{x\over y}{{dx}\over x} +\int_y^\infty g(x)xD_{y\over
x}{{dx}\over x}\Big]\cr
&+\int_1^\infty f(y){{dy}\over y}\int_0^1g(x)yD_{x\over
y}{{dx}\over x}\Big) +\int_0^1 f(y){{dy}\over y}\int_1^\infty g(x)xD_{y\over
x}{{dx}\over x}\Big),D_1\rangle\cr
=&\Big( \int_0^1f(y){{dy}\over y}\Big[\int_0^yg(x)y\langle D_{x\over
y},D_1\rangle{{dx}\over x} +\int_y^1g(x)x\langle D_{y\over 
x},D_1\rangle{{dx}\over
x}\Big]\cr
&+\int_1^\infty f(y){{dy}\over y}\Big[\int_1^yg(x)y\langle D_{x\over
y},D_1\rangle{{dx}\over x} +\int_y^\infty g(x)x\langle D_{y\over
x},D_1\rangle{{dx}\over x}\Big]\cr
&+\int_1^\infty f(y){{dy}\over
y}\int_0^1g(x)y\langle D_{x\over y},D_1\rangle{{dx}\over x}\Big) +\int_0^1
f(y){{dy}\over y}\int_1^\infty g(x)x\langle D_{y\over x},D_1\rangle{{dx}\over
x}\Big).\cr}$$ Now accordingly call each of the terms from the beginning as 
(3.1),
(3.2), (4.1), (4.2), (5.1) and (5.2), we then see that for (3.1), (resp.
(3.2), (4.1), (4.2), (5.1) and (5.2)) we have $0\leq x\leq y\leq
1$, (resp. $0\leq y\leq x\leq 1$, $1\leq x\leq y\leq \infty$ and 
$\infty\geq x\geq
y\geq 1$,) hence by our axioms, for $0\leq x\leq y\leq 1$
$$\langle y\langle D_{x\over y},D_1\rangle=\langle D_y,D_x\rangle,$$
(resp. for $0\leq y\leq x\leq 1$ $$x\langle D_{y\over x},D_1\rangle=\langle
D_x,D_y\rangle,$$
for $1\leq x\leq y\leq\infty$  $$\langle D_{x\over y},D_1\rangle =x\langle
D_x,D_y\rangle,$$
for $\infty\geq x\geq y\geq 1$ $$\langle D_{y\over x},D_1\rangle=y\langle
D_x,D_y\rangle,$$
for $0\leq x\leq 1\leq y\leq \infty$
$$y\langle D_{x\over y},D_1\rangle=y\langle D_y,D_x\rangle,$$
and for $0\leq y\leq 1\leq x\leq\infty$
$$x\langle D_{y\over x},D_x\rangle =x\langle D_y,D_x\rangle.)$$
We see that the latest combination is simply
$$\eqalign{~&\Big( \int_0^1f(y){{dy}\over y}\Big[\int_0^yg(x)\langle
D_{y},D_x\rangle{{dx}\over x} +\int_y^1g(x)\langle D_{x},D_y\rangle{{dx}\over
x}\Big]\cr &+\int_1^\infty f(y){{dy}\over y}\Big[\int_1^yg(x)(xy)\langle
D_{x},D_y\rangle{{dx}\over x} +\int_y^\infty g(x)(xy)\langle
D_{x},D_y\rangle{{dx}\over x}\Big]\cr
&+\int_1^\infty f(y){{dy}\over y}\int_0^1g(x)
y\langle D_{y},D_x\rangle{{dx}\over x}\Big)
+\int_0^1 f(y){{dy}\over
y}\int_1^\infty g(x)x\langle D_{y},D_x\rangle{{dx}\over x}\Big)\cr}$$ which 
certainly
is nothing but
$$\eqalign{~&\int_0^1f(y){{dy}\over y}\int_0^1g(x){{dx}\over x}\langle 
D_y,D_x\rangle
+\int_0^1f(y){{dy}\over y}\int_1^\infty g(x){{dx}\over x}\cdot x\langle
D_y,D_x\rangle\cr &+\int_1^\infty f(y){{dy}\over y}\int_0^1g(x){{dx}\over 
x}\cdot
y\langle D_y,D_x\rangle +\int_1^\infty f(y){{dy}\over y}\int_1^\infty 
g(x){{dx}\over
x}\cdot xy\langle D_y,D_x\rangle.\cr}$$ Now by defintion, this latest 
combination is
simply
$$\langle D_{\hat f},D_{\hat g}\rangle.$$ This then completes the proof of 
the Fixed
Point Formula.

To see the relative degree relation, we have the following formal arguments.
$$\eqalign{~&\langle D_{\hat f},D_0\rangle\cr =&\int_0^1 f(x)\langle
D_x,D_0\rangle {{dx}\over x}+\int_1^\infty f(x) x\langle D_x,D_0\rangle 
{{dx}\over
x}\cr =&\int_0^1 f(x)x{{dx}\over x}+\int_1^\infty f(x)x{{dx}\over x}\cr
=&\int_0^\infty f(x)x{{dx}\over x}=\hat f(1)\cr}$$ and
$$\eqalign{~&\langle D_{\hat f},D_\infty\rangle\cr
=&\int_0^1 f(x)\langle
D_x,D_\infty\rangle {{dx}\over x}+\int_1^\infty f(x) x\langle 
D_x,D_\infty\rangle
{{dx}\over x} =\int_0^1 f(x){{dx}\over x}+\int_1^\infty f(x){{dx}\over x}\cr
=&\int_0^\infty f(x){{dx}\over x}=\hat f(0).\cr}$$ So here a standard 
regularization
is needed. For details, see 2.5 below on Not So Serious Convergence Problems.

Finally, let see how the explicit formula is established. Here the functional
equation plays a key role as in function fields case.

Indeed, by definition,
$$\langle D_{\hat f},D_1\rangle
=\int_0^1 f(x)\langle D_x,D_1\rangle {{dx}\over
x}+\int_1^\infty f(x) x\langle D_x,D_1\rangle {{dx}\over x}.$$ Now by the micro
explicit formula, we have
$$\eqalign{~&\langle D_{\hat f},D_1\rangle\cr
=&\int_0^1 f(x)\Big(\langle D_x,D_0\rangle+\langle D_x,D_\infty\rangle
-\sum_{\xi_{\bf Q}(s)=0}x^s\Big)
  {{dx}\over
x}+\int_1^\infty f(x) x\Big(\langle D_x,D_0\rangle+\langle D_x,D_\infty\rangle
-\sum_{\xi_{\bf Q}(s)=0}x^{s-1}\Big)\rangle {{dx}\over x}.\cr}$$
This is because by the local explicit formula, we have for $x\in [0,1]$,
$$\langle D_x,D_1\rangle=\langle D_x,D_0\rangle+\langle D_x,D_\infty\rangle
-\sum_{\xi_{\bf Q}(s)=0}x^s.$$
Hence, if $x\geq 1$, we have
$$x\langle D_x,D_1\rangle=x\langle D_{1\over x},D_1\rangle
=x\langle D_{1\over x},D_0\rangle+x\langle D_{1\over x},D_\infty\rangle
-\sum_{\xi_{\bf Q}(s)=0}x^{1-s}
=x\langle D_{1\over x},D_0\rangle+x\langle D_{1\over x},D_\infty\rangle
-\sum_{\xi_{\bf Q}(s)=0}x^s$$ where in the last step, we use the functional 
equation
for the Riemann zeta function. Thus in particular, we see that
for $x\geq 1$, $$\langle D_x,D_1\rangle
=\langle D_x,D_0\rangle+\langle D_x,D_\infty\rangle-\sum_{\xi_{\bf
Q}(s)=0}x^s$$ holds as well. Certainly, then
$$\langle D_{\hat f},D_1\rangle=
\langle D_{\hat f},D_0\rangle +\langle D_{\hat f},D_\infty\rangle
+\int_0^\infty f(x) \sum_{\xi_{\bf Q}(s)=0}x^s {{dx}\over x}
=\hat f(0)+\hat f(1)-\sum_{\xi_{\bf Q}(s)=0}\hat f(s).$$
\vskip 0.30cm
\noindent
{\bf C.2.4.  The Riemann Hypothesis}
\vskip 0.30cm
To finally relate our intersection
with the Riemann Hypothesis, we should have a certain positivity.
That is to say, we need an analog of the so-called Hodge Index Theorem.
\vskip 0.30cm
\noindent
{\bf A Weak Version of Hodge Index Theorem.} {\it With the same notation as
above,  the self-intersection of the global divisor
$\displaystyle{L_{\hat f}:=D_{\hat f}-\hat f(1) D_\infty
-{\hat f}(0)D_0}$ is non-positive, i.e.,
$$\langle L_{\hat f},L_{\hat f}\rangle\leq 0.\eqno(*)$$}

\noindent
{\it Remark.} We call $(*)$ a weak version of Hodge Index Theorem, since
$$\langle L_{\hat f},D_0+D_\infty\rangle=0.$$

 From now on, let us assume that $(*)$ holds.
Then by the Key Relations in the previous section, with a direct
calculation, we certainly will arrive at
$$\sum_{\xi_{\bf Q}(s)=0}\hat f(s)\cdot\hat f(1-s)\geq
0.$$ This is very nice, since then, following Weil, we may get the Riemann 
Hypothesis from
this latest inequality.
(See e.g., page 342 of the second
edition of Lang's Algebraic Number Theorm for more details.)
\vskip 0.30cm
\noindent
{\bf C.2.5. Not so serious Convergence Problem}
\vskip 0.30cm
We add some remarks on the formal calculation appeared in 2.3. Roughly 
speaking, to
justify them, what we meet is a certain regularized process. This may be 
done as what
Jorgenson and Lang do in their lecture notes on Basic Analysis of 
Regularized Series and
Product. More precisely, motivated by the definitions of  hyperbolic 
Green's functions (of
Selberg, Hejhal, Groos and Zagier,)  Ray-Singer's analytic torsions,
we may first introduce {\it imaginary divisors}
$D_{\hat f,s}$ by setting
$$D_{\hat f,s}={1\over{\Gamma(s)}}
\Big(\int_0^1 f(x) D_x x^s{{dx}\over x}+\int_1^\infty f(x) xD_x 
x^s{{dx}\over x}
\Big)$$ for certain type of suitable funcions $f$ for $s$ whose real parts 
are sufficiently
large; then assume that in our model $D_{\hat f,s}$
has a meromorphic continuation  to the half space ${\rm Re}(s)\geq
-\varepsilon$ with $\varepsilon>0$, from which we could finally get a 
well-defined
$D_{\hat f}$ after removing the singularity at $s=0$.
\vskip 0.30cm
\noindent
{\bf C.2.6. Weil's Explicit Formula and Two Dimensional Geometric Arithmetic
Intersections}
\vskip 0.30cm
I still have not explain why we say the above intersection is indeed an
intersection over the geometric arithmetic surface $\overline {{\rm Spec}({\bf
Z})}\times \overline {{\rm Spec}({\bf Z})}$. To understand this, we have to
use yet another fundamental result of Weil, the Weil Explicit Formula.

Note that by the Key Relations, we obtain the following crucial formula:

$$\langle D_{\hat f}, D_1\rangle = \hat f(0)+\hat f(1)-\sum_{\xi_{\bf
Q}(s)=0}\hat f(s).\eqno(*)$$

On the other hand, as an intersection over $\overline {{\rm Spec}({\bf
Z})}\times \overline {{\rm Spec}({\bf Z})}$,
$\langle D_{\hat f}, D_1\rangle$ should be counted locally over each point
$(p,q)\in \overline {{\rm Spec}({\bf
Z})}\times \overline {{\rm Spec}({\bf Z})}$, i.e., we should have the
decomposition
$$\langle D_{\hat f}, D_1\rangle=\sum_{p,q\leq \infty} \langle D_{\hat f},
D_1\rangle_{(p,q)}.$$ Now note that $D_1$ is the diagonal, so
besides the points $(p,p)$, $p\leq\infty$ on the diagonal,
$\langle D_{\hat f}, D_1\rangle_{(p,q)}$ is naturally zero. With this, then we
would have
$$\langle D_{\hat f}, D_1\rangle=\sum_{p\leq \infty} \langle D_{\hat f},
D_1\rangle_{(p,p)}.$$

Clearly, this latest expression suggests that via $(*)$ above
$$\hat f(0)+\hat f(1)-\sum_{\xi_{\bf
Q}(s)=0}\hat f(s)=\sum_{p\leq \infty} W_p(f),\eqno(**)$$
where for each places $p$ of {\bf Q}, i.e., the primes $p$ and the Archimedean
place $\infty$. Without any mistake, it is then nothing but Weil's explicit 
formula.

On the other hand, this interpretation then naturally leads to a question about
the explicit formula for the micro intersection. Recall that one of the key
assumption for our micro intersection is that,
for $x\in [0,1]$,
$$\langle D_x,D_1\rangle=\langle D_x,D_0\rangle+\langle
D_x,D_\infty\rangle-\sum_{\xi_{\bf Q}(s)=0}x^s.$$
Therefore, if we believe that $\langle\cdot,\cdot\rangle$ is indeed a two
dimensional intersection on $\overline {{\rm Spec}({\bf
Z})}\times \overline {{\rm Spec}({\bf Z})}$, then  we similarly should have
$$1+x-\sum_{\xi_{\bf Q}(s)=0}x^s=\sum_{p\leq \infty}w_p(x),\qquad {\rm for}\
x\in [0,1].\eqno(***)$$ We will call $(***)$  the {\it micro explicit 
formula} of Crem\'er,
motivated by Jongenson and Lang's
\lq ladder principle\rq.
\vskip 0.30cm
\noindent
{\it Remark.} There are also fundamental works of
Deninger and Quillen on  the Riemann Hypothesis, based on certain cohomological
consideration. While  these approaches appear quite different, they share 
one common
part,  the Weil Explicit Formula.
\vskip 0.45cm
\noindent
{\li C.3. Towards A  Geo-Ari Cohomology in Lower Dimensions}
\vskip 0.30cm
\noindent
{\bf C.3.1. Classical Approach in Diemnsion One}
\vskip 0.30cm
It is clear that we should go beyond a geo-ari intersection: To complete 
the picture, for
example, we need to establish an analog of Hodge index theorem in Geometric 
Arithmetic,
since it is where the positivity comes. Thus, from our experiences with
Algebraic Geometry and Arakelov Theory, we  are led to develop a corresponding
geo-ari cohomology theory.
 
It is our belief that a general yet well-behavior cohomology theory is at the
present time beyond our reach. However this does not mean that we cannot do 
anything
about it. After all, what we need is a practical yet uniform cohomology 
(and intersection)
in dimensions one and two, such that  duality,  adjunction formula and
Riemann-Roch are satisfied.

With this in mind, we recall what happens in geometry for cohomology in 
dimension one.
Classical approaches (to cohomology), such as the one cited in Serre's GTM 
on Algebraic
Groups and Class Fields, consist of  two aspects, namely, the algebraic one 
and the analytic
one. Moreover, with an algebraic approach,  we may develop a general sheaf 
cohomology
theory,  thanks to the work of Grothendieck. Relatively speaking, for 
analytic apsect,
we have achieved very little. Thus, we want to  explore it, since we 
understand that a
gro-ari cohomology should be based on the  alnalytic discussion.

Let $C$ be a regular irreducible reduced projective curve of genus $g$
defined over a field $k$ with $D$ a divisor on $C$. Denote by $F$ its 
associated function
fields.  Then the keys to an algebraic cohomology theory  may be summarized 
as follows.

\noindent
(1) by definition, the 0-th cohomology group of (the divisor class associated
to) $D$  is given by $H^0(C,D):=\{f\in F:{\rm div}(f)+D\geq 0\}$;

\noindent
(2) From the short exact sequence of sheaves $$0\to {\cal O}_C(D)\to
{\cal F}\to {\cal F}/{\cal O}_C(D)\to 0,$$ where ${\cal F}$ denotes the 
constant sheaf on $C$
associated to the function field $F$,
we get a long exact sequence $$0\to H^0(C,{\cal O}_C(D))\to
F\to H^0(C,{\cal F}/{\cal O}_C(D))\to H^1(C,{\cal O}_C(D))\to 0.$$ Thus, by 
definition,
there should be  canonical isomorphisms $$H^1(X,{\cal O}_C(D))\simeq {\bf 
A}/({\bf
A}(D)+F).$$ Here {\bf A} denotes the associated
adelic ring, and $${\bf A}(D):=\{(r_p)\in {\bf A}:{\rm ord}_p(r_p)+{\rm 
ord}_p(D)\geq 0\};$$

\noindent
(3) For any point $p$, from the structural exact sequence of sheaves $0\to 
{\cal O}_C(D)\to
{\cal O}_C(D+p)\to {\cal O}_C(D)|_p\to 0$,  we get a long exact sequence of 
cohomology
$$0\to H^0(C,{\cal O}_C(D))\to  H^0(C,{\cal O}_C(D+p))\to  H^0(C,{\cal 
O}_C(D)|_p)\to
H^1(C,{\cal O}_C(D))\to  H^1(C,{\cal O}_C(D+p))\to 0;$$

\noindent
(4) By studying the residue pairing, we get a canonical isomorphism
$${\bf A}/({\bf
A}(D)+F) \simeq (H^0(C,{\cal
O}_C(K_C-D)))^\vee$$ which in particular implies that  there exists a natural
duality between $H^0(C,D)$ and $H^1(C,K_C-D)$;

\noindent
(5) Cohomology groups $H^0$ and $H^1$ are all finite dimensional vector spaces.
Thus in particular, $$h^0(C,D)-h^1(C,D)=h^0(C,D+p)-h^1(C,D+p)-1.$$ This 
then implies
the duality and the Riemann-Roch $$h^0(C,D)-h^1(C,D)=d(D)-(g-1).$$
\vskip 0.30cm
Next, we describe the analytic aspect of the above cohomology, which is 
based on a
study about certain quotient and sub spaces associated to the adelic ring 
{\bf A}.
\vskip 0.30cm
\noindent
(1$'$) For any divisor $D$ on $C$, define its associated cohomology groups 
by $H^0(X,{\cal
O}_C(D))={\bf A}(D)\cap F$ and
$H^1(C,D):={\bf A}/{\bf A}(D)+ F$;

\noindent
(2$'$) There is the following commutative 9-diagram
$\Sigma(D)$, whose rows and colums  are all exact:
$$\matrix{&&0&&0&&0&&\cr
&&\downarrow&&\downarrow&&\downarrow&&\cr
0&\to&{\bf A}(D)\cap F &\to&{\bf A}(D)&\to& {\bf A}(D)/{\bf A}(D)\cap 
F\simeq {\bf
A}(D)+F/F&\to&0\cr
&&\downarrow&&\downarrow&&\downarrow&&\cr
0&\to& F&\to&{\bf A}&\to& {\bf A}/F&\to&0\cr
&&\downarrow&&\downarrow&&\downarrow&&\cr
0&\to&F/{\bf A}(D)\cap F&\to&{\bf A}/{\bf A}(D)&\to& {\bf A}/{\bf A}(D)+ 
F&\to&0\cr
&&\downarrow&&\downarrow&&\downarrow&&\cr
&&0&&0&&0.&&\cr}$$
So, we get the exact sequences
$$0\to H^0(C,{\cal O}_C(D))\to F\to F/{\bf A}(D)\cap F\to 0, \qquad
0\to F/{\bf A}(D)\cap F\to {\bf A}/{\bf A}(D)\to H^1(C,{\cal O}_C(D))\to 
0$$ which clearly is
equivalent to the exact sequence (2) above;

\noindent
(3$'$)  For any $p\in C$,  there exists
a natural morphism from $\Sigma(D)$ to $\Sigma(D+p)$. Thus, by the five lemma,
based on the fact that $F$ and ${\bf A}$ are the same in these two 
9-diagrams, we
obtain the exact sequences
$$0\to {\bf A}(D)\cap F\to {\bf A}(D+p)\cap F\to {\bf A}(D+p)\cap F/{\bf 
A}(D)\cap
F\to 0$$
$$0\to {\bf A}(D+p)\cap F/{\bf A}(D)\cap
F\to  {\bf A}(D+p)/{\bf A}(D)\to {\bf A}(D+p)+ F/{\bf A}(D)+
F\to 0$$ $$0\to {\bf A}(D+p)+ F/{\bf A}(D)+ F\to
{\bf A}/{\bf A}(D)+F\to {\bf A}/{\bf A}(D+p)+F\to 0.$$ Clearly, these are 
equivalent to the
exact sequence (3) above;

\noindent
(4$'$) Residue pairing works at the level of adelic language as well by the 
self-dual property
of {\bf A}.
\vskip 0.30cm
Therefore, provided that  we know how to count the terms
involved, namely, that we have an analog of (5) above,  we can develop
a cohomology theory using only adelic language (which satisfies
duality and the Riemann-Roch theorem).
\vskip 0.30cm
However, generally speaking, the counting is a very difficult one. In algebraic
approach, this is based on the fact that all coherent sheaves are locally 
finitely generated.
In analytic approach, the counting will be based on the spacial analytic 
properties of {\bf A}.
To explain it, as an example, we now consider the simplest case, namely, 
when the constant
field $k$ is
${\bf F}_q$, the finite field with $q$ elements. (Over number fields, we
count the geo-ari cohomology by Tate's  Fourier analysis over {\bf A}.)

Recall that with respect to the natural topology on {\bf A}, $F$ is 
discrete and
${\bf A}(D)$ is compact. Thus in particular, $H^0(C,{\cal O}_C(D))={\bf 
A}(D)\cap F$ is finite.
Similarly, since ${\bf A}(D)$ is compact and ${\bf A}/F$ is compact so 
${\bf A}(D)+F/F$ is
again compact. But $k$ is a finite fields, so compactness implies 
finiteness, and the number of
elements in a finite dimensional space is simply $q$ to the power of the
corresponding dimension.
   Thus, with repsect to the natural Haar measures on  the
associated groups induced from that on {\bf A},
$${\rm Vol}({\bf A}(D))=q^{h^0(C,{\cal O}_C(D))}\cdot {\rm Vol}({\bf 
A}(D)+F/F)$$ and
$${\rm Vol}({\bf A}/F)=q^{h^1(C,{\cal O}_C(D))}\cdot {\rm Vol}({\bf 
A}(D)+F/F).$$
Therefore, $$q^{h^0(C,{\cal O}_C(D))-h^1(C,{\cal O}_C(D))}={{{\rm Vol}{\bf 
A}(D)}
\over {{\rm Vol}({\bf A}/F)}}.$$  Easily by definition, $${\rm Vol}{\bf
A}(D)=q^{d(D)},\quad {\rm Vol}({\bf A}/F)=q^{g-1},$$ so we obtain an 
analytic proof of
the Riemann-Roch theorem $$h^0(C,{\cal O}_C(D))-h^1(C,{\cal 
O}_C(D))=d(D)-(g-1).$$
\vskip 0.30cm
\noindent
{\bf C.3.2. Chevalley's Linear Compacity}
\vskip 0.30cm
The above discussion for curves over finite fields is not valid for general 
curves, since even
the space is finite-dimensional, it is not finite itself. To overcome this 
difficult,
Chevalley introduced his linear compacity. Next, we indicate how
Chevalley's method works. Instead of recalling all the details, I  indicate 
what are the
essential points involved. (For details, we recommend the reader to consult 
Iwasawa's
Princeton lecture notes.)
\vskip 0.30cm
\noindent
(0) There existes the 9-diagram as above;

\noindent
(1) (Additive Structure) Among subquotient groups of the adelic ring are
topological spaces called discrete objects and linearly compact objects. 
Moreover, all groups
used in our 9-diagram are supposed to be  locally linearly compact, i.e.,
they are either extensions of discrete objects by linearly compact objects, 
or extensions of
linearly compact objects by discrete objects, or simply generated by 
finitely many discrete
objects and linearly compact objects. In particular, {\bf A} is selfdual 
and locally linearly
compact.

\noindent
(2) Discrete objects and linearly compact objects are dual to each other.
Thus, if an object is both locally linearly compact and discrete, it is 
then isomorphic to a
finite dimensional $k$ vector space, and hence the dimension may be counted;

\noindent
(3) (Multiplicative Structure) Under the multiplication, ${\bf A}$ is 
self-dual. As a direct
consequence, we get the duality.
\vskip 0.30cm
As such, the counting may be proceeded as follows to offer the Riemann-Roch:

\noindent
(i) By definition, ${\rm dim}({\bf A}(0)\cap F)=1$ and ${\rm dim}({\bf A}/{\bf
A}(0)+F)=g$;

\noindent
(ii) By counting local contribution, $[{\bf A}(D+p):{\bf 
A}(D)]=d(D+p)-d(D)=1$, despite the
fact that
${\bf A}(D)$ and ${\bf A}(D+p)$ cannot be counted;

\noindent
(iii) By definition, $$[{\bf A}(D+p):{\bf A}(D)]={\rm dim} \big({\bf 
A}(D+p)/{\bf A}(D)\big)$$

\noindent
(iv) The fundamental theorem of isomorphisms for groups implies that
$$\eqalign{{\rm dim} \big({\bf A}(D+p)/{\bf A}(D)\big) =&{\rm dim} \big({\bf
A}(D+p)+F/{\bf A}(D)+F\big)+{\rm dim}
\big({\bf A}(D+p)\cap F/{\bf A}(D)\cap F\big)\cr
=&{\rm dim} {\bf A}/{\bf
A}(D)+F-{\rm dim}{\bf A}/{\bf A}(D+p)+F+{\rm dim}
{\bf A}(D+p)\cap F-{\rm dim}{\bf A}(D)\cap F\cr
=&\big(h^0(C,D+p)-h^1(C,D+p)\big)
-\big(h^0(C,D)-h^1(C,D)\big).\cr}$$
\vskip 0.30cm
At this point, I would like to point out that it is quite essential to 
combine the approach here
with our approach to geo-ari cohomology for number fields via Fourier analysis.
\vskip 0.30cm
\noindent
{\bf C.3.3. Adelic Approach in Geometric Dimension Two}
\vskip 0.30cm
Now we  consider two dimensional case. Besides the definition of cohomology 
groups and
the counting of these cohomology groups,
from algebraic geometry and Arakelov theory, we know that the (weak)
Riemann-Roch theorem and duality may be obtained via the adjunction 
formula, a one
dimensional Riemann-Roch and a long exact sequence of cohomology groups 
resulting from
a short exact sequence of sheaves
$$0\to {\cal O}(D)\to{\cal O}(D+C)\to {\cal O}(D)|_C\to 0$$ for a divisor 
$D$ and a regular
curve $C$ on the surface, and a residue discussion.

On the other hand, just as for curves, we want to develop a two dimensional 
cohomology
theory using only adelic language. So two parts are involved:

\noindent
(I) (Algebraic Structure) Definition of cohomology groups and some 
associated structural
exact sequences;

\noindent
(II) (Analytic Structure) Counting of the cohomology groups in (1).
\vskip 0.30cm
We start with the algebraic structure. In algebraic geometry, this is 
essentially given by
Parshin as a by-product of his discussion on residues. As above, we here 
only indicate the
main points. (Interesting reader may consult Prshin's original paper for 
the details.)

So for a surface $S$ with function field $F$, we  introduce its associated 
ring of
adeles {\bf A}. Inside {\bf A} are two subrings, which we denote by ${\bf 
A}_0$  and
${\bf A}_1$, respectively. Similarly, for a divisor $D$ on $S$, introduce
its associated subgroup
${\bf A}(D)$ as in curve case. In particular, then we have the following 
three 9-diagrams:
$$\matrix{&&0&&0&&0&&\cr
&&\downarrow&&\downarrow&&\downarrow&&\cr
0&\to&{\bf A}(D)\cap \big({\bf A}_0\cap {\bf A}_1\big) &\to&{\bf A}(D)&\to& 
{\bf
A}(D)/{\bf A}(D)\cap\big({\bf A}_0\cap {\bf A}_1\big) &\to&0\cr
&&\downarrow&&\downarrow&&\downarrow&&\cr 0&\to& {\bf A}_0\cap {\bf
A}_1&\to&{\bf A}&\to& {\bf A}/{\bf A}_0\cap {\bf A}_1&\to&0\cr
&&\downarrow&&\downarrow&&\downarrow&&\cr
0&\to&{\bf A}_0\cap {\bf A}_1/{\bf A}(D)\cap
{\bf A}_0\cap {\bf A}_1&\to&{\bf A}/{\bf A}(D)&\to& {\bf A}/{\bf A}(D)+ 
{\bf A}_0\cap
{\bf A}_1&\to&0\cr &&\downarrow&&\downarrow&&\downarrow&&\cr
&&0&&0&&0.&&\cr}$$  for which isomorphisms $${\bf
A}(D)/{\bf A}(D)\cap\big({\bf A}_0\cap {\bf A}_1\big) \simeq {\bf 
A}(D)+{\bf A}_0\cap
{\bf A}_1/{\bf A}_0\cap {\bf A}_1$$ and $${\bf A}_0\cap {\bf A}_1/{\bf 
A}(D)\cap
{\bf A}_0\cap {\bf A}_1\simeq {\bf A}(D)+\big( {\bf A}_0\cap {\bf 
A}_1\big)/{\bf
A}(D)$$ are used;
$$\matrix{&&0&&0&&0&&\cr
&&\downarrow&&\downarrow&&\downarrow&&\cr
0&\to&{\bf A}(D)\cap \big({\bf A}_0\cap {\bf A}_1\big) &\to&{\bf A}(D)\cap {\bf
A}_1&\to& {{{\bf A}(D)\cap {\bf A}_1}\over{{\bf A}(D)\cap \big({\bf 
A}_0\cap {\bf
A}_1\big)}}&\to&0\cr
&&\downarrow&&\downarrow&&\downarrow&&\cr
0&\to& {\bf A}(D)\cap {\bf A}_0&\to&{\bf A}(D)\cap \big({\bf A}_0+ {\bf 
A}_1\big)
&\to& {{{\bf A}(\bf D)\cap \big({\bf A}_0+ {\bf A}_1\big)}\over{{\bf 
A}(D)\cap {\bf
A}_0}}&\to&0\cr
&&\downarrow&&\downarrow&&\downarrow&&\cr
0&\to& {{{\bf A}(D)\cap {\bf
A}_0}\over{{\bf A}(D)\cap\big({\bf A}_0\cap {\bf A}_1\big)}}&\to&
{{{\bf A}(D)\cap \big({\bf A}_0+{\bf A}_1\big)}\over{ {\bf A}(D)\cap {\bf 
A}_1}}
&\to& {{{\bf A}(D)\cap \big({\bf A}_0+ {\bf A}_1\big)}\over{{\bf A}(D)\cap 
{\bf A}_0+{\bf
A}(D)\cap {\bf A}_1}}&\to&0\cr &&\downarrow&&\downarrow&&\downarrow&&\cr
&&0&&0&&0&&\cr}$$
for which the isomorphisms
$$ {\bf A}(D)\cap {\bf A}_1/{\bf A}(D)\cap \big({\bf A}_0\cap {\bf
A}_1\big)\simeq {\bf A}(D)\cap {\bf A}_0+{\bf A}(D)\cap {\bf A}_1/{\bf 
A}(D)\cap {\bf
A}_0$$ and
$${\bf A}(D)\cap {\bf
A}_0/{\bf A}(D)\cap\big({\bf A}_0\cap {\bf A}_1\big)\simeq
  {\bf A}(D)\cap {\bf A}_0+{\bf A}(D)\cap {\bf A}_1/{\bf A}(D)\cap {\bf
A}_1$$ are used; and
$$\matrix{&&0&&0&&0&&\cr
&&\downarrow&&\downarrow&&\downarrow&&\cr
0&\to&{\bf A}(D)\cap \big({\bf A}_0+{\bf A}_1\big) &\to&{\bf A}(D)&\to& 
{\bf A}(D)/{\bf
A}(D)\cap \big({\bf A}_0+{\bf A}_1\big)&\to&0\cr
&&\downarrow&&\downarrow&&\downarrow&&\cr
0&\to& {\bf A}_0+{\bf A}_1&\to&{\bf A}&\to& {\bf A}/\big({\bf A}_0+{\bf
A}_1\big)&\to&0\cr &&\downarrow&&\downarrow&&\downarrow&&\cr
0&\to&\big({\bf A}_0+{\bf A}_1\big)/{\bf A}(D)\cap \big({\bf A}_0+{\bf
A}_1\big)&\to&{\bf A}/{\bf A}(D)&\to& {\bf A}/{\bf A}(D)+ \big({\bf A}_0+{\bf
A}_1\big)&\to&0\cr &&\downarrow&&\downarrow&&\downarrow&&\cr
&&0&&0&&0&&\cr}$$
for which the isomorphisms
$$ {\bf A}(D)/{\bf
A}(D)\cap \big({\bf A}_0+{\bf A}_1\big)\simeq {\bf A}(D)+\big({\bf A}_0+{\bf
A}_1\big)/\big({\bf A}_0+{\bf A}_1\big)$$ and $$\big({\bf A}_0+{\bf 
A}_1\big)/{\bf
A}(D)\cap \big({\bf A}_0+{\bf A}_1\big)\simeq {\bf A}(D)+\big({\bf A}_0+{\bf
A}_1\big)/{\bf A}(D)$$ are used.

Now, define the cohomology groups $H^i$, $i=0,1,2$ by setting
  $$H^0(S,D):={\bf A}(D)\cap
\big({\bf A}_0\cap {\bf A}_1\big);$$
$$H^1(S,D):={\bf A}(D)\cap \big({\bf A}_0+ {\bf A}_1\big)/{\bf A}(D)\cap 
{\bf A}_0+{\bf
A}(D)\cap {\bf A}_1;$$ and
$$H^2(S,D):={\bf A}/{\bf A}(D)+ \big({\bf A}_0+{\bf
A}_1\big).$$
Moreover,  by working with $D+C$ for regular curve $C$ on $S$, we
could get another set of three 9-diagrams, to which there are natural 
morphisms from
the three 9-diagrams for $D$ above. As a direct consequence, by a snake 
chasing, we arrive
at the following  long exact sequence of cohomologies, which as stated 
above, plays
a key role in the induction process:
$$\eqalign{0\to &H^0(S,D)\to H^0(S,D+C)\to H^0(C,D+C|_C)\cr
\to& H^1(S,D)\to H^1(S,D+C)\to
H^1(C,D+C|_C)\to H^2(S,D)\to H^2(S,D+C)\to 0.\cr}$$ In this way, provided that
a good counting is availble, we then are able to give an adelic approach to 
the weak
Riemann-Roch in dimension two by using the Riemann-Roch for curves with the 
help of the
adjunction formula.
\vskip 0.30cm
To understand the analytic structure, motivated by Chevalley's theory on 
linear compacity
for adelic rings of algebraic curves, we need to study the terms
used in the above algebraic discussion. The key is the self-dual property 
of {\bf A}.
In fact, we have the following canonical isomorphisms:

\noindent
(a) $A_0^\perp\simeq {\bf A}_1, A_1^\perp\simeq {\bf A}_0,$ and ${\bf 
A}(D)^\perp={\bf
A}(K_S-D)$, where $K_S$ denotes a canonical divisor of $S$.  In particular,
$$\big({\bf A}_0+{\bf A}_1\big)^\perp\simeq {\bf A}_1\cap {\bf A}_2,\qquad
\big({\bf A}_1\cap {\bf A}_2\big)^\perp\simeq {\bf A}_0+{\bf A}_1.$$

\noindent
(b) Duality between $H^0(S,D)$ and $H^2(S,K_S-D)$
$$\eqalign{&\big({\bf A}(D)\cap{\bf A}_0\cap {\bf
A}_1\big)^\vee\simeq {\bf A}/\big({\bf A}(D)\cap{\bf A}_0\cap {\bf
A}_1\big)^\perp\cr
\simeq& {\bf A}/{\bf A}(D)^\perp+ {\bf
A}_0^\perp+{\bf A}_1^\perp\simeq {\bf A}/{\bf A}(K_S-D)+ \big({\bf A}_0+{\bf
A}_1\big),\cr}$$

\noindent
(c) Note that ${\bf A}_0\cap {\bf A}_1=F$ is nothing but the function field 
of $S$. Moreover,
one checks that algebrically
$$\eqalign{&{\bf A}(D)\cap \big({\bf A}_0+ {\bf A}_1\big)/{\bf A}(D)\cap 
{\bf A}_0+{\bf
A}(D)\cap {\bf A}_1\cr
\simeq &{\bf A}_0\cap \big({\bf A}(D)+{\bf A}_1\big)/\big({\bf A}(D)\cap
{\bf A}_0+{\bf A}_0\cap{\bf A}_1\big)\cr
\simeq&
{\bf A}_1\cap \big({\bf A}(D)+{\bf A}_0\big)/\big({\bf A}(D)\cap
{\bf A}_1+{\bf A}_0\cap{\bf A}_1\big);\cr}$$

\noindent
(d) Duality between $H^1(S,D)$ and $H^1(S,K_S-D)$:
$$\eqalign{&\Big({\bf A}_0\cap \big({\bf A}(D)+{\bf A}_1\big)/\big({\bf 
A}(D)\cap
{\bf A}_0+{\bf A}_0\cap{\bf A}_1\big)\Big)^\vee\cr
\simeq& \big({\bf A}(D)\cap
{\bf A}_0+{\bf A}_0\cap{\bf A}_1\big)^\perp/\big({\bf A}_0\cap \big({\bf 
A}(D)+{\bf
A}_1\big)\big)^\perp\cr
\simeq &\big({\bf A}(K_S-D)+{\bf A}_1\big)\cap\big({\bf A}_0+{\bf
A}_1\big)/\big({\bf A}_1+{\bf A}(K_S-D)\cap {\bf A}_1\big)\cr
\simeq&{\bf A}(K_S-D)\cap \big({\bf A}_0+ {\bf A}_1\big)/{\bf A}(K_S-D)\cap 
{\bf A}_0+{\bf
A}(D)\cap {\bf A}_1.\cr}$$

Thus,  to attack (II), the analytic structure aiming at a reasonable 
counting, we should
introduce a geo-ari theory for surfeces which is compactible with (a), (b), 
(c) and (d) above,
similar to that of linear compactness of Chevalley used in proving the 
Riemann-Roch in
dimension one. For example, we should have a notion of geo-ari compactness 
such that if a
space is both discrete and geo-ari compact, it should be of finite 
dimensional; moreover, the
duality should transform discrete spaces to geo-ari compact speces and vice 
versa.
\vskip 0.30cm
Note that for curves, if the base field is  finite, then we  can equally 
use Fourier
analysis to do the counting. Thus, for two dimensional surfaces, a similar
discussion should work as well.  All this then inevitably leads to a 
geo-ari cohomology
in dimensions one and two over number fields, our primary goal.
\vfill
\eject
\vskip 0.45cm
\centerline {\li REFERENCES}
\vskip 0.25cm
\item{[An]} A.N. Andrianov,  Euler products that correspond to Siegel's 
modular forms of
genus $2$, Russian Math. Surveys {\bf 29}:3(1974), 45-116
\vskip 0.25cm
\item {[Ar]} S. Arakelov: Intersection theory of divisors on an arithmetic
surface, Izv. Akad.  Nauk SSSR Ser. Mat., 38 No. {\bf 6}, (1974)
\vskip 0.25cm
\item{[A]} E. Artin, Quadratische K\"orper im Gebiete der h\"oheren
Kongruenzen, I,II, {\it Math. Zeit}, {\bf 19} 153-246 (1924) (See also
{\it Collected Papers}, pp. 1-94,  Addison-Wesley 1965)
\vskip 0.25cm
\item{[AT]} E. Artin \& J. Tate, {\it Class Field Theory}, Benjamin Inc, 1968
\vskip 0.25cm
\item{[At]} M. Atiyah, Vector Bundles over an elliptic curve, {\it
Proc. LMS, VII}, 414-452 (1957) (See also {\it Collected Works}, Vol. 1,
pp. 105-143, Oxford Science Publications, 1988)
\vskip 0.25cm
\item{[AB]} M. Atiyah \& R. Bott, The Yang-Mills equations over Riemann 
surfaces.  Philos.
Trans. Roy. Soc. London Ser. A  258  (1983),  no. 1505, 523--615.
\vskip 0.25cm
\item{[Ba]} W. Banaszczyk, New bounds in some transference theorems in
the geometry of numbers, Math. Ann., {\bf 296}, (1993), 625-635
\vskip 0.25cm
\item{[Be1]} A. Beilinson, Higher regulators and values of $L$-functions, 
J. Soviet Math. {\bf
30} (1985), 2036-2070
\vskip 0.25cm
\item{[Be2]} A. Beilinson, Height pairings between algebraic cycles, in 
{\it Current
Trends in Arithmetical Algebraic Geometry}, Contemporary Math. {\bf 67} 
(1987), 1-24
\vskip 0.25cm
\item{[BR]} U. Bhosle \& A. Ramanathan, Moduli of parabolic $G$-bundles on 
curves.
Math. Z.  202  (1989),  no. 2, 161--180.
\vskip 0.25cm
\item{[Bl]} S. Bloch, A note on height pairings, Tamagawa numbers, and the 
BSD conjecture,
Inv. Math. {\bf 58} (1980), 65-76
\vskip 0.25cm
\item{[BK]} S. Bloch \& K. Kato, $L$-functions and Tamagawa numbers of 
motives, in {\it The
Grothendieck Festschrift}, Vol. I, Progress in Math. 86 (1990), 333-400
\vskip 0.25cm
\item{[BV]} E. Bombieri \& J. Vaaler, On Siegel's lemma. Invent. Math.
{\bf 73} (1983), no. 1, 11--32.
\vskip 0.25cm
\item{[Bo]} A. Borel, Some finiteness properties of adele groups over
number fields, Publ. Math., IHES, {\bf 16} (1963) 5-25
\vskip 0.25cm
\item{[Bor]} A. Borisov, Convolution structures and arithmetic cohomology, 
preprint
\vskip 0.25cm
\item{[B-PGN]} L. Brambila-Paz, I. Grzegorczyk, \& P.E. Newstead, Geography 
of Brill-Noether Loci
for Small slops, J. Alg. Geo. {\bf 6}  645-669 (1997)
\vskip 0.25cm
\item{[Bu1]} J.-F. Burnol, Weierstrass points on arithmetic surfaces.  Invent.
Math.  107  (1992),  no. 2, 421--432.
\vskip 0.25cm
\item{[Bu2]} J.-F. Burnol, Remarques sur la stabilite' en arithme'tique. 
Internat. Math.
Res. Notices  1992,  no. 6, 117--127
\vskip 0.25cm
\item {[CF]} J.W.S. Cassels \& A. Fr\"ohlich, {\it Algebraic Number 
Theory}, Academic
Press, 1967
\vskip 0.25cm
\item{[Co]} A. Connes, Trace formula in
noncommutative geometry and the zeros of the Riemann zeta function. 
Selecta Math. (N.S.)
5  (1999),  no. 1, 29--106.
\vskip 0.25cm
\item{[CS]} J.H. Conway \& N.J.A. Sloane,  {\it Sphere packings, lattices 
and groups},
Springer-Verlag,  1993
\vskip 0.25cm
\item{[C]} H. Cram\'er,  Studien \"uber die Nullstellen der Riemannschen 
Zetafunktion, Math.
Z. {\bf 4} (1919), 104-130
\vskip 0.25cm
\item{[D1]} P. Deligne, Vari\'et\'es de Shimura, in {\it Automorphic Forms, 
Representations,
and $L$-functions}, Proc. of Sym. in Pure Math., Vol {\bf 33} (2), 
(1977),247-290
\vskip 0.25cm
\item{[D2]} P. Deligne, Cat\'egories Tannakiennes, in  {\it The
Grothendieck Festschrift}, Vol. I, Progress in Math. 86 (1990), 111-195
\vskip 0.25cm
\item{[DM]} P. Deligne \& J.S. Milne, Tannakian categories, in {\it Hodge 
Cycles, Motives and
Shimura Varieties}, LNM {\bf 900}, (1982), 101-228
\vskip 0.25cm
\item{[De1]} Ch. Deninger,  On the $\Gamma$-factors attached to motives.
Invent. Math. {\bf 104} (1991), no. 2, 245--261.
\vskip 0.25cm
\item{[De2]} Ch. Deninger,  Local $L$-factors of
motives and regularized determinants.  Invent. Math.  {\bf 107}  (1992), 
no. 1, 135--150.
\vskip 0.25cm
\item{[De3]} Ch. Deninger, Motivic $L$-functions and regularized determinants,
Motives (Seattle, WA, 1991), 707-743, Peoc. Sympos. Pure Math, {\bf 55} Part 1,
AMS, Providence, RI, 1994
\vskip 0.25cm
\item{[De4]} Ch.  Deninger,  Lefschetz trace formulas and explicit formulas in
analytic number theory. J. Reine Angew. Math. {\bf 441} (1993), 1--15.
\vskip 0.25cm
\item{[De5]} Ch. Deninger, Some analogies between number theory
and dynamical systems on foliated spaces. Proceedings of the International 
Congress of
Mathematicians, Vol. I (Berlin, 1998).  Doc. Math.  1998,  Extra Vol. I, 
163--186
\vskip 0.25cm
\item{[DS]} Ch. Deninger \& M.  Schr\"oter, A distribution-theoretic
proof of Guinand's functional equation for Cram\'er's $V$-function and 
generalizations. J.
London Math. Soc. (2) {\bf 52} (1995), no. 1, 48--60.
\vskip 0.25cm
\item{[DR]} U.V. Desale \& S. Ramanan, Poincar\'e polynomials of the variety of
stable bundles, Maeh. Ann {\bf 216}, 233-244 (1975)
\vskip 0.25cm
\item{[FW]} G. Faltings \& G. W\"ustholz,  Diophantine approximations on 
projective
spaces.  Invent. Math.  116  (1994),  no. 1-3, 109--138.
\vskip 0.25cm
\item{[Fo1]} J.-M.  Fontaine, Repre'sentations $p$-adiques semi-stables.
Aste'risque No. 223 (1994), 113--184.
\vskip 0.25cm
\item{[Fo2]} J.-M.  Fontaine,
Repre'sentations $l$-adiques potentiellement semi-stables.
Aste'risque  No. 223 (1994), 321--347.
\vskip 0.25cm
\item{[FV]} M.D. Fried \& H. V\"olklein,  The embedding problem over a
Hilbertian PAC-field. Ann. of Math. (2) {\bf 135} (1992), no. 3, 469--481
\vskip 0.25cm
\item{[GS]} G. van der Geer \& R. Schoof, Effectivity of Arakelov
Divisors and the Theta Divisor of a Number Field, preprint, 1998
\vskip 0.25cm
\item{[Gr1]} D.R. Grayson, Reduction theory using semistability. Comment. 
Math. Helv. {\bf 59}
(1984), no. 4, 600--634
\vskip 0.25cm
\item{[Gr2]} D.R. Grayson, Reduction theory using semistability. II. 
Comment. Math. Helv. {\bf
61}  (1986),  no. 4, 661--676.
\vskip 0.25cm
\item{[Gu1]} R.C. Gunning, {\it Lectures on Riemann surfaces}, Princeton 
Math. Notes, 1966
\vskip 0.25cm
\item{[Gu2]} R.C. Gunning, {\it Lectures on vector bundles over Riemann 
surfaces},
Princeton Math. Notes, 1967
\vskip 0.25cm
\item{[Ha1]} Sh. Haran,  Index theory, potential theory, and the Riemann 
hypothesis.
$L$-functions and arithmetic (Durham, 1989),  257--270, London Math. Soc. 
Lecture Note
Ser., 153, Cambridge Univ. Press, Cambridge, 1991.
\vskip 0.25cm
\item{[Ha2]} Sh. Haran,  Riesz potentials and explicit sums in arithmetic.
Invent. Math.  101  (1990),  no. 3, 697--703.
\vskip 0.25cm
\item{[HN]} G. Harder \& M.S. Narasimhan, On the cohomology groups of 
moduli spaces
of vector bundles over curves, Math Ann. {\bf 212}, (1975) 215-248
\vskip 0.25cm
\item{[H]} H. Hasse, {\it Mathematische Abhandlungen}, Walter
de Gruyter, Berlin-New York, 1975.
\vskip 0.25cm
\item{[H]} D. Hilbert, \"Uber die Theorie der relativ-Abelschen 
Zahlk\"orper, Nachr. der K.
Ges. der Wiss. G\"ottingen, 377-399 (1898)
\vskip 0.25cm
\item{[Iw1]} K. Iwasawa, {\it Algebraic theory of algebraic functions}, 
Lecture Notes at
Princeton Univ. noted by T. Kimura, 1975
\vskip 0.25cm
\item{[Iw2]} K. Iwasawa, Letter to Dieudonn\'e, April 8, 1952, in
N.Kurokawa and T. Sunuda: {\it Zeta Functions in Geometry}, Advanced
Studies in Pure Math. {\bf 21} (1992), 445-450
\vskip 0.25cm
\item{[JL1]} J. Jorgenson \& S. Lang, On
Crame'r's theorem for general Euler products with functional equation. 
Math. Ann.  297
(1993),  no. 3, 383--416.
\vskip 0.25cm
\item{[JL2]} J. Jorgenson \& S. Lang, {\it Basic analysis of regularized 
series and products}.
Lecture Notes in Mathematics, 1564. Springer-Verlag, Berlin, 1993. viii+122
\vskip 0.25cm
\item{[JL3]} J. Jorgenson \& S. Lang, {\it  Explicit formulas}. Lecture Notes
in Mathematics, 1593. Springer-Verlag, Berlin, 1994. viii+154 pp.
\vskip 0.25cm
\item{[Ka]} K. Kato, Generalized explicit reciprocity laws, Adv. Studies in 
Contemporary
Math., 1 (1999), 57-126
\vskip 0.25cm
\item{[Ka]} Y. Kawada,  Class formulations. VI. Restriction to a
subfamily.  J. Fac. Sci. Univ. Tokyo Sect. I  8  1960 229--262 (1960).
\vskip 0.25cm
\item{[KT]} Y. Kawada \& J. Tate, On the Galois cohomology of unramified 
extensions of
function fields in one variable.  Amer. J. Math.  77,  (1955). 197--217.
\vskip 0.25cm
\item{[Ke]} G.R. Kempf, Instability in invariant theory, Ann. Math. {\bf 
108} (1978),
299-316
\vskip 0.25cm
\item{[L1]} S. Lang, {\it Algebraic Number Theory},
Springer-Verlag, 1986
\vskip 0.25cm
\item{[L2]} S. Lang, {\it Fundamentals on Diophantine Geometry},
Springer-Verlag, 1983
\vskip 0.25cm
\item{[L3]} S. Lang, Introduction to Arakelov theory, Springer Verlag, 1988
\vskip 0.25cm
\item{[La]} R. P.  Langlands,  Automorphic representations, Shimura
varieties, and motives. Ein Ma"rchen.  Automorphic forms, representations and
$L$-functions pp. 205--246, Proc. Sympos. Pure Math., XXXIII, Amer. Math. 
Soc., Providence, R.I., 1979.
\vskip 0.25cm
\item{[Li]} X. Li, A note on the Riemann-Roch theorem for function fields.
{\it Analytic number theory}, Vol. {\bf 2},  567--570, Progr. Math.,
{\bf 139}, Birkha"user, 1996
\vskip 0.25cm
\item{[Ma]} Y. Manin, New dimensions in geometry,  LNM {\bf 1111} (1985)
59-101
\vskip 0.25cm
\item{[M]} B. Mazur, Deforming Galois representations, in {\it Galois 
groups over {\bf Q}},
MSRI Publ. {\bf 16}, (1989), 385-437
\vskip 0.25cm
\item{[MS]} V.B.  Mehta \& C.S.  Seshadri, Moduli of vector bundles on 
curves with
parabolic structures.  Math. Ann.  248  (1980), no. 3, 205--239.
\vskip 0.25cm
\item{[Mi]} H. Minkowski, {\it Geometrie der Zahlen}, Leipzig and Berlin, 1896
\vskip 0.25cm
\item{[Mo]} A. Moriwaki, Private communication
\vskip 0.25cm
\item{[M]} D. Mumford, {\it Geometric Invariant Theory}, Springer-Verlag, 
Berlin
(1965)
\vskip 0.25cm
\item{[Na]} M.  Namba, {\it Branched coverings and algebraic functions}.
Pitman Research Notes in Mathematics Series, 161. Longman Scientific \& 
Technical, Harlow;
John Wiley \& Sons, Inc., New York, 1987. viii+201 pp.
\vskip 0.25cm
\item{[NR]} M.S. Narasimhan \& S. Ramanan, Moduli of vector bundles on a 
compact Riemann surfaces,
Ann. of Math. {\bf 89} 14-51 (1969)
\vskip 0.25cm
\item{[NS]} M.S.  Narasimhan \& C.S.  Seshadri,  Stable
and unitary vector bundles on a compact Riemann surface. Ann. of Math. (2) 
82 1965
\vskip 0.25cm
\item{[Ne]} J. Neukirch, {\it Algebraic Number Theory}, Grundlehren der
Math. Wissenschaften, Vol. {\bf 322}, Springer-Verlag, 1999
\vskip 0.25cm
\item{[Pa]} A.N. Parshin,  On the arithmetic of two-dimensional schemes. I.
Distributions and residues. (Russian)  Izv. Akad. Nauk SSSR Ser. Mat.  40 
(1976), no. 4,
736--773, 949.
\vskip 0.25cm
\item{[RR]} S. Ramanan \& A. Ramanathan,  Some remarks on the instability 
flag.  Tohoku
Math. J. (2)  36  (1984),  no. 2, 269--291.
\vskip 0.25cm
\item{[Ra1]} A. Ramanathan, Stable principal
bundles on a compact Riemann surface.  Math. Ann.  213  (1975), 129--152.
\vskip 0.25cm
\item{[Ra2]} A. Ramanathan, Moduli for principal bundles over algebraic 
curves. I.  Proc.
Indian Acad. Sci. Math. Sci.  106  (1996),  no. 3, 251--328.
\vskip 0.25cm
\item{[Ra3]} A.  Ramanathan, Moduli for
principal bundles over algebraic curves. I. Proc. Indian Acad. Sci. Math. 
Sci. 106 (1996), no.
4, 421--449.
\vskip 0.25cm
\item{[RSS]} M. Rapoport, N. Schappacher, \& P. Schneider, {\it Beilinson's 
conjectures on
special values of $L$-functions}, Perspectives in Math, Vol. {\bf 4.} (1988)
\vskip 0.25cm
\item{[RZ]} M. Rapoport \& Th.  Zink, {\it  Period spaces for $p$-divisible 
groups}. Annals of
Mathematics Studies, 141. Princeton University Press, Princeton, NJ, 1996.
\vskip 0.25cm
\item{[Sch]} A.J. Scholl, An introduction to Kato's Euler systems, in {\it 
Galois
Representations in arithmetic algebraic geometry}, London Math. Soc. 
Lecture Note Series
{\bf 254}, (1998) 379-460
\vskip 0.25cm
\item{[S]} J.-P. Serre, {\it Algebraic Groups and Class Fields}, GTM 117, 
Springer (1988)
\vskip 0.25cm
\item{[Se1]} C.S. Seshadri,  Space of unitary vector bundles on a compact 
Riemann
surface.  Ann. of Math. (2)  85  (1967) 253--336.
\vskip 0.25cm
\item {[Se2]} C. S. Seshadri, {\it Fibr\'es vectoriels sur les courbes 
alg\'ebriques}, Asterisque
{\bf 96}, 1982
\vskip 0.25cm
\item{[Sh]} I. Shafarevich,  {\it Lectures on minimal models and birational
transformations of two dimensional schemes.} Notes by C. P. Ramanujam. Tata 
Institute of
Fundamental Research Lectures on Mathematics and Physics, No. {\bf 37} (1966)
\vskip 0.25cm
\item{[Shi]} G. Shimura, {\it Introduction to the arithmetic theory of 
automorphic
functions}, Iwanami, (1971)
\vskip 0.25cm
\item {[St1]} U. Stuhler,   Eine Bemerkung zur Reduktionstheorie 
quadratischer Formen,
Arch. Math. (Basel) {\bf 27} (1976), no. 6, 604--610
\vskip 0.25cm
\item {[St2]} U. Stuhler, Zur Reduktionstheorie der positiven quadratischen 
Formen. II,
Arch. Math. (Basel)  {\bf 28}  (1977), no. 6, 611--619
\vskip 0.25cm
\item{[T]} J. Tate, Fourier analysis in number fields and Hecke's
zeta functions, Thesis, Princeton University, 1950
\vskip 0.25cm
\item{[Ti]} J. Tilouine, {\it  Deformations of Galois representations and Hecke
algebras}, Published for The Mehta Research Institute of Mathematics and 
Mathematical
Physics, Allahabad; by Narosa Publishing House, New Delhi, 1996.
\vskip 0.25cm
\item{[To]} B. Totaro, Tensor products in $p$-adic Hodge theory.  Duke Math.
J.  83  (1996),  no. 1, 79--104.
\vskip 0.25cm
\item{[W1]} A. Weil, G\'en\'eralisation des fonctions ab\'eliennes, J. Math 
Pures et Appl, {\bf
17}, (1938) 47-87
\vskip 0.25cm
\item{[W2]} A. Weil, {\it Sur les courbes alg\'ebriques et les vari\'et\'es qui
s'en d\'eduisent}, Herman, Paris (1948)
\vskip 0.25cm
\item{[W3]} A. Weil, Sur les formules explicites de la the'orie des 
nombres. (French)  Izv.
Akad. Nauk SSSR Ser. Mat.  36  (1972), 3--18.
\vskip 0.25cm
\item{[W4]} A. Weil, {\it Basic Number Theory}, Springer-Verlag, 1973
\vskip 0.25cm
\item{[W5]} A. Weil,  {\it Adeles and algebraic groups}. With appendices by M.
Demazure and Takashi Ono. Progress in Mathematics, 23. Birkha"user, Boston, 
Mass., 1982.
iii+126 pp.
\vskip 0.25cm
\item{[We1]} L. Weng, $\Omega$-admissible theory II: Deligne pairings over 
moduli spaces
of punctured Riemann surfaces, Math. Ann {\bf 320} (2001), 239-283
\vskip 0.25cm
\item{[We2]} L. Weng,  Riemann-Roch Theorem, Stability and
New Zeta Functions for Number Fields, preprint
\vskip 0.25cm
\item{[We3]} L. Weng, Refined Brill-Noether Locus and Non-Abelian Zeta
Functions  for Elliptic  Curves, preprint
\vskip 0.25cm
\item{[We4]} L. Weng, Constructions of Non-Abelian Zeta Functions for 
Curves, submitted
\vskip 0.25cm
\item{[We5]} L. Weng, Stability and New Non-Abelian Zeta Functions, to appear
\vskip 0.25cm
\item{[We6]} L. Weng, Non-abelian class field theory for Riemann surfaces, 
preprint
\vskip 0.25cm
\item{[Zh]} S. Zhang,  Heights and reductions of semi-stable varieties.
Compositio Math.  {\bf 104} (1996),  no. 1, 77--105.
\vskip 1.30cm
\noindent
Mail Address: Graduate School of Mathematics, Nagoya University, Nagoya
464-8602, Japan

\noindent
E-Mail Address: weng@math.nagoya-u.ac.jp
\end